\documentclass[11pt,a4paper,onecolumn,twosided,rotate]{IEEEtran} %%%%%%%%%%%%%%%%%%%%%%%%%%%
\usepackage{amsmath,amsthm,amssymb,amsfonts,upref,cite,epsf,color}
\usepackage{pst-all,pstricks,pst-node,pst-text,pst-3d}
\usepackage{hyperref}
\usepackage{pstricks-add}
\usepackage{multido}
\usepackage{bm}
\usepackage{graphicx,epsfig,rotating} % Optionen "final", "draft"(= keine bilder)
\usepackage{psfrag,texdraw,bm}
\usepackage{nomencl}
\usepackage{amsmath}
\usepackage{amssymb}
\usepackage{subfigure}
\usepackage{xspace}
\usepackage{dsfont}
\usepackage{pifont}
\usepackage{listings,color}   %% color für Farbmarkierungen
\usepackage[T1]{fontenc}
\usepackage{textcomp}
\newtheorem{theorem}{Theorem}[section]
\newtheorem{definition}[theorem]{Definition}
\newtheorem{lemma}[theorem]{Lemma}

\small\normalsize
\def \expect {\mathsf{E} }
\def \Prob {{\rm P} }
\def \twiddle[#1] {e^{-j \frac{2 \pi}{N}  #1 }}
\def \twiddleneg[#1] {e^{j \frac{2 \pi}{N}  #1 }}

\newcommand{\eq}{\,=\,}
\def\llr{\rho_{\mathcal{E}\rmv,\mathbf{x}_{0}}}

\DeclareMathOperator{\linspan}{span}

\def\ML_est{\hat{\mathbf{x}}_{\text{ML}}}
\newcommand{\CRBfull}{{C}ram\'{e}r--{R}ao bound\xspace}

\newcommand{\be}{\begin{equation}}
\newcommand{\ee}{\end{equation}}
\newcommand{\ist}{\hspace*{.2mm}}
\newcommand{\rmv}{\hspace*{-.2mm}}
\newcommand{\basisfuncRKHS}{u}
\newcommand{\transp}[1]{{#1}^T} % transpose of a matrix

\newcommand{\vecestproblem}{\left(  \mathcal{X}, f(\mathbf{y};\mathbf{x}), \mathbf{g}(\cdot) \right)}
\newcommand{\scalarestproblem}{\big(  \mathcal{X}, f(\mathbf{y};\mathbf{x}),g(\cdot) \big)}

\allowdisplaybreaks

\makeatother

\begin{document}
%%%%%%%%%%%%%%%%%%%%%%%%%%%%%%%%%%%%%%%%%%%%%%%%%%%%%%%%%%

\title{\vspace*{2mm}The RKHS Approach to Minimum Variance\\[-4mm]
Estimation Revisited: Variance Bounds,\\[-4mm]Sufficient Statistics, and Exponential Families
 \vspace*{4mm}}
\author{\emph{ Alexander Jung\textsuperscript{a} (corresponding author), Sebastian Schmutzhard\textsuperscript{$\ist$b}\rmv, and Franz Hlawatsch\textsuperscript{$\ist$a}
\thanks{This work was supported by the FWF under Grants S10602-N13 (Signal and Information Representation) and 
S10603-N13 (Statistical Inference) within the National Research Network SISE 
and by the WWTF under Grant MA 07-004 (SPORTS).}
\thanks{
First revision; submitted to the IEEE Transactions on Information Theory, \today}
%% This paper is mainly based on the first authors PHD thesis \cite{JungPHD}.}
} \\[2mm]
{\small\emph{\textsuperscript{a}}Institute of Telecommunications, Vienna University of Technology;
%% }\\[-2.7mm]
%% {\small Gusshausstrasse\ 25/389, 1040 Vienna, Austria}\\[-2.7mm]
%% {\small Phone: +43 1 58801 38963, Fax: +43 1 58801 38999, e-mail: 
\{ajung,\,fhlawats\}@nt.tuwien.ac.at} \\[-2mm]
{\small\emph{\textsuperscript{b}}NuHAG, Faculty of Mathematics, University of Vienna;
%% }\\[-2.7mm]
%% {\small A-1090 Vienna, Austria; e-mail: 
sebastian.schmutzhard@univie.ac.at} \\[0mm]
}

%\address{\normalsize \textsuperscript{a}Institute of Communications and Radio-Frequency Engineering, Vienna University of Technology\\[-1.2mm]
%\normalsize Gusshausstrasse 25/389, A-1040 Vienna, Austria; e-mail: ajung@nt.tuwien.ac.atÊ\\[1mm]
%\normalsize \textsuperscript{b}Technion---Israel Institute of Technology, %\\[-1.2mm]
%\normalsize Haifa 32000, Israel; e-mail: zvikabh@tx.technion.ac.il}

%% \date{Submitted to IEEE Transactions on Information Theory, \today}

%%%%%%%%%%%%%%%%%%%%%%%%%%%%%%%%%%%%%%%%%%%%%%%%%%%%%%%%%%
\maketitle
%%%%%%%%%%%%%%%%%%%%%%%%%%%%%%%%%%%%%%%%%%%%%%%%%%%%%%%%%%

\vspace*{-10mm}

\begin{abstract}
\vspace*{-1.5mm}
The mathematical theory of reproducing kernel Hilbert spaces (RKHS) provides powerful tools for minimum variance estimation (MVE)
problems. Here, we extend the classical RKHS-based analysis of MVE in several directions. We develop a geometric formulation of  
five known lower bounds on the estimator variance (Barankin bound, Cram\'{e}r--Rao bound, constrained Cram\'{e}r--Rao bound, 
Bhattacharyya bound, and Hammersley-Chapman-Robbins bound) in terms of orthogonal projections onto a subspace of the RKHS associated with a given MVE problem.
%We define the property of \emph{differentiability} of an RKHS and demonstrate its close relation to the subspace associated with the Cram\'{e}r--Rao bound.
We show that, under mild conditions, the Barankin bound (the tightest possible lower bound on the estimator variance) 
%% !!!typically!!! 
is a lower semi-continuous function of the parameter vector. 
We also show that the RKHS associated with an MVE problem remains unchanged if the observation is replaced by a sufficient statistic. 
Finally, for MVE problems conforming to an exponential family of distributions, we derive novel closed-form lower bounds on the estimator variance 
and show that a reduction of the parameter set leaves the minimum achievable variance unchanged.
\vspace*{-1mm}
\end{abstract}

\begin{keywords}
\vspace*{-1.5mm}
Minimum variance estimation, exponential families, RKHS, Cram\'{e}r--Rao bound, Barankin bound, Hammersley--Chapman--Robbins 
bound, Bhattacharyya bound, locally minimum variance unbiased estimator.
%% Parameter estimation
\vspace*{2mm}
\end{keywords}

%% \newpage %%%%%%%%%

%%%%%%%%%%%%%%%%%%%%%%%%%%%%%%%%%%%%%%%%%%%%%%%%%%%%%%%%%%
\section{Introduction}\label{sec.intro}
%%%%%%%%%%%%%%%%%%%%%%%%%%%%%%%%%%%%%%%%%%%%%%%%%%%%%%%%%%

%RKHS approach to MVE; definition of estimation problem: Parameter set, parameter function and statistical model; LMV estimator; UMV estimator; exponential families; Examples: LGM, SLGM, covariance estimation 
%--> show explicitly how these problems are obtained from the generic exponential family 

%% In this paper, 
We consider the problem of estimating the value $\mathbf{g}(\mathbf{x})$ of a known deterministic function $\mathbf{g}(\cdot)$ evaluated at an 
unknown nonrandom parameter vector $\mathbf{x} \rmv\in\rmv \mathcal{X}$,
%%  \rmv\subseteq\rmv \mathbb{R}^{N}\!$.
%% , whose dimension $N$ is assumed known. 
%% The parameter vector $\mathbf{x}$ is unknown except for the fact that it belongs to the known parameter set $\mathcal{X}$.
where 
%% its dimension $N$ and 
the parameter set $\mathcal{X}$ is known.
%%  and it is also known that $\mathbf{x}$ belongs to a parameter set $\mathcal{X} \subseteq \mathbb{R}^{N}\!$. 
%% One contribution of this paper is an analysis of the effect of restricting???reducing the parameter set $\mathcal{X}$, which corresponds to an 
%% %% can be interpreted as in 
%% increased \emph{a priori} knowledge. 
The estimation of $\mathbf{g}(\mathbf{x})$ is based on an observed vector $\mathbf{y}$,
%%  \rmv\in\rmv \mathbb{R}^{M}\!$,
%% , with known dimension $M$ is known. 
%% The observation $\mathbf{y}$ 
which is modeled as a random vector with an associated probability measure
\cite{BillingsleyProbMeasure} $\mu^{\mathbf{y}}_{\mathbf{x}}$ 
or, as a special case, an associated probability density function (pdf) $f(\mathbf{y}; \mathbf{x})$, both parametrized by $\mathbf{x} \rmv\in\rmv \mathcal{X}$. 
%% Thus, the observed random vector $\mathbf{y}$ is related to the parameter vector $\mathbf{x}$ via its probabilistic properties. 
%% We will refer to both the set of measures ${\{ \mu^{\mathbf{y}}_{\mathbf{x}} \}}_{\mathbf{x} \in \mathcal{X}}$ and the set of pdf's 
%% ${\{ f(\mathbf{y}; \mathbf{x}) \}}_{\mathbf{x} \in \mathcal{X}}$
%% %% , which are obtained for every parameter vector $\mathbf{x} \in \mathcal{X}$, 
%% as the \emph{statistical model}. 
%% 
%% In general, we will assume that the statistical model is dominated by a reference measure on $\mathbb{R}^{N}\!$, 
%% e.g., the Lebesgue measure on $\mathbb{R}^{M}$ \cite{BillingsleyProbMeasure,LC,HalmosRadonNikodymSuffStat,IbragimovBook}, 
%% so that we can define for each $\mathbf{x} \rmv\in\rmv \mathcal{X}$ the \emph{probability density function} (pdf) $f(\mathbf{y}; \mathbf{x})$ of $\mathbf{y}$.
%% %% the observation vector. 
%% In this case, we will refer also to the set of pdfs ${\{ f(\mathbf{y}; \mathbf{x}) \}}_{\mathbf{x} \in \mathcal{X}}$ as a \emph{statistical model}. 
%% %We will discuss the concept of probability measures and pdfs in more detail below. 
%% 
%% Based on the reproducing kernel Hilbert space (RKHS) framework, 
%% 
More specifically, we study the problem of \emph{minimum variance estimation} (MVE), where one aims at finding estimators with minimum variance
under the constraint of a prescribed bias. Our treatment of MVE will be based on the mathematical framework and methodology
of \emph{reproducing kernel Hilbert spaces} (RKHS).

\vspace{-7.5mm}

\textcolor{red}{
\subsection{State of the Art and Motivation}\label{sec.state}
%%%%%%%%%%%%%%%%%%%%%%%%%%%%%%%%%%%%%%%%%%%%%%%%%%%%%%%%%%
%% \paragraph{State of the Art}
}
\vspace{1mm}
The RKHS approach to MVE was introduced in the seminal papers \cite{Parzen59} and \cite{Duttweiler73b}. 
\textcolor{red}{On a general
%% high 
level, the theory of RKHS yields
%% allows for 
efficient methods for
%% to deal with 
high-dimensional optimization problems. These methods are popular, e.g., 
%%  for coping with high-dimensional optimization problems occurring 
in machine learning \cite{SmaleZhou2007,Cucker02onthe}.
For the MVE problem considered here, the optimization problem 
%% consists in minimizing 
is the minimization of the estimator variance subject to a bias constraint.} 
 The RKHS approach to MVE enables a consistent and intuitive geometric treatment of the MVE problem. In particular, the determination of the minimum achievable variance (or Barankin bound) and of the
%% accompanying 
locally minimum variance estimator reduces to the computation of the squared norm and isometric image of a specific vector---representing the prescribed estimator bias---that belongs to the RKHS associated with the estimation problem. \textcolor{red}{This reformulation is interesting from a theoretical perspective; in addition, it may also be the basis
for an efficient computational evaluation.}
Furthermore, a wide class of lower bounds on the minimum achievable variance (and, in turn, on the variance of any estimator) is obtained by performing projections onto subspaces of the RKHS. \textcolor{red}{Again, this enables an efficient computational evaluation of these bounds.}
%%  associated with the estimation problem. 

\textcolor{red}{A specialization to estimation problems involving sparsity constraints was presented in \cite{RKHSAsilomar2010,RKHSISIT2011,AlexSebICASSP2011}. For certain special cases of these sparse estimation problems, the RKHS approach allows the derivation of closed-form expressions of 
the minimum achievable variance and the corresponding locally minimum variance estimators.
The RKHS approach has also proven to be a valuable tool for the analysis of estimation problems involving continuous-time random processes \cite{Parzen59,Duttweiler73b,Kailath71}. 
}

%% \vspace{-1mm}

\vspace{-3mm}

\textcolor{red}{
\subsection{Contribution and Outline}\label{sec.contr}
%%%%%%%%%%%%%%%%%%%%%%%%%%%%%%%%%%%%%%%%%%%%%%%%%%%%%%%%%%
}

%% \paragraph{Contribution}
The main contributions of this paper 
%% can be summarized as follows.
concern an RKHS-theoretic analysis of the performance of MVE, with a focus on questions related to lower variance bounds,
sufficient statistics, and observations conforming to an exponential family of distributions.
%% exponential family-based estimation problems.
First, 
%% based on the introduction of the concept of a \emph{differentiable RKHS}, 
we give a 
%% purely 
geometric 
%% RKHS-based 
interpretation of some well-known lower bounds on the estimator variance. 
The tightest of these bounds, i.e., the Barankin bound, is proven to be a lower semi-continuous function of the parameter vector $\mathbf{x}$ under mild conditions. 
We then analyze the role of a sufficient statistic from the RKHS viewpoint.
%of the effect of two specific modification of an estimation problem: The first modification is the reduction of the parameter set and the second modification is the replacement of the observation by a sufficient statistic. 
In particular, we prove that the RKHS associated with an estimation problem remains unchanged if the observation $\mathbf{y}$ is replaced by any sufficient statistic. 
Furthermore, we characterize the RKHS for estimation problems with observations conforming to an exponential family of distributions.
%% exponential family-based estimation problems.
%%  whose observation follows an exponential family. 
It is found that this RKHS has a strong structural property, 
%% is a differentiable RKHS, which constitutes a strong structural property, 
%% Moreover, the RKHS 
and that it is explicitly related to the moment-generating function of the 
%% respective 
exponential family. Inspired by this relation, we derive novel lower bounds on the estimator variance, and we analyze the effect of parameter set reductions. 
The lower bounds have a particularly simple form.

%% \paragraph{Outline}
The remainder of this paper is organized as follows. 
\textcolor{red}{
In Section \ref{SecMVE}, basic elements of MVE are reviewed and the RKHS approach to MVE is summarized.
%%  and introduce the concept of a differentiable RKHS. 
In Section \ref{sec_rkhs_interpret_well_known_bounds}, we present an RKHS-based geometric interpretation of known variance bounds 
and demonstrate the lower semi-continuity of the Barankin bound.}
%%  !!! under mild conditions!!!. 
The effect of replacing the observation by a sufficient statistic is studied in Section \ref{sec_suff_stat_RKHS}. 
In Section \ref{sec_rkhs_exp_fam}, the RKHS for exponential family-based estimation problems is investigated,
%% estimation problems whose observation follows an exponentially family. We prove that 
%% the RKHS associated to an estimation problem involving an exponential family is differentiable and explicitly 
%% tied to the moment generating and cumulant function of the exponential family. 
%% In Section \ref{sec_param_set_reduction_exp_fam}, still considering exponential family-based problems, 
novel lower bounds on the estimator variance are derived, and the effect of a parameter set reduction is analyzed.
%The variance bounds are then illustrated by means of numerical experiments in Section \ref{sec_num_results}.
We note that the proofs of most of the new results presented can be found in the doctoral dissertation \cite{JungPHD} and will be referenced in each case. 

\vspace{-2mm}

\textcolor{red}{
\subsection{Notation and Basic Definitions}\label{sec.notation}
%%%%%%%%%%%%%%%%%%%%%%%%%%%%%%%%%%%%%%%%%%%%%%%%%%%%%%%%%%
}

%% \paragraph{Notation and basic definitions} 
%% \paragraph{Notation}
We will use the shorthand notations
%% The set of natural numbers is denoted by 
$\mathbb{N} \triangleq \{1,2,3, \ldotsÊ\}$, 
%% and the set of nonnegative integers by 
$\mathbb{Z}_{+} \triangleq \{0,1,2,\ldots\}$, and 
$[N] \triangleq \{1,2,\ldots,N\}$.
%% $[N] \triangleq \{0,1,\ldots,N\!-\!1\}$.
%Given a nonnegative integer $L \in \mathbb{Z}_{+}$, we denote by $[L]$ the set of all numbers up to $L$, i.e., $[L] \triangleq \{1,\ldots,L\}$ if $L>0$ and for the special case $L=0$ we set $[0] \triangleq \emptyset$ where $\emptyset$ denotes the empty set. 
The open ball in $\mathbb{R}^{N}$ with radius $r > 0$ and centered at $\mathbf{x}_{c}$ is defined as 
$\mathcal{B}(\mathbf{x}_{c},r) \triangleq \big\{ \mathbf{x} \!\in\! \mathbb{R}^{N} \ist\big|\ist {\| \mathbf{x} \rmv-\rmv \mathbf{x}_{c} \|}_{2} \!<\rmv r \big\}$.
We call 
%% a point 
$\mathbf{x} \in \mathcal{X}\subseteq \mathbb{R}^{N}$ an \emph{interior point} if $\mathcal{B}(\mathbf{x},r) \subseteq \mathcal{X}$ for some $r \!>\! 0$. 
The set of all interior points of $\mathcal{X}$ is called the \emph{interior} of $\mathcal{X}$ and denoted $\mathcal{X}^{\text{o}}$. A set $\mathcal{X}$ is called \emph{open} if 
%% it coincides with its interior, i.e., 
$\mathcal{X} = \mathcal{X}^{\text{o}}$. 
%The boundary $\boundary(\mathcal{X})$ of the set $\mathcal{X} \subseteq \mathbb{R}^{N}$ is defined 
%by the set points $\mathbf{x} \in \mathcal{X}$ such that for all $r \!>\!0$ the ball $\mathcal{B}(\mathbf{x},r)$ contains at least one point $\mathbf{x}'$ outside $\mathcal{X}$, i.e., 
%$\mathcal{B}(\mathbf{x},r) \cap \mathcal{X} \neq \mathcal{B}(\mathbf{x},r)$. 

Boldface lowercase (uppercase) letters denote vectors (matrices). The superscript $^{T}$ stands for transposition. 
The $k$th entry of a vector $\mathbf{x}$ and the entry in the $k$th row and $l$th column of a matrix $\mathbf{A}$ are denoted 
by ${(\mathbf{x})}_{k} = x_k$ and ${(\mathbf{A})}_{k,l} = A_{k,l}$, respectively.
The $k$th unit vector is denoted by $\mathbf{e}_{k}$, and the identity matrix of size $N \!\times\! N$ by $\mathbf{I}_{N}$.
The Moore-Penrose pseudoinverse \cite{golub96} of a rectangular matrix $\mathbf{F} \in \mathbb{R}^{M \times N}$ is denoted by $\mathbf{F}^{\dagger}$. 
%% The set of positive definite Hermitian $NÊ\!\times\! N$ matrices is denoted $\mathcal{S}_{++}^{N}$.
%and by $\Ker(\mathbf{F}) \triangleq \{ \mathbf{x} \rmv\in\rmv \mathbb{R}^{N} |\ist  \mathbf{F} \mathbf{x} \rmv=\rmv \mathbf{0} \}$ its kernel (or null space). 

A function $f(\cdot) \!: \mathcal{D} \!\rightarrow\! \mathbb{R}$, with $\mathcal{D} \!\subseteq\! \mathbb{R}^{N}\rmv$, is said to be
\emph{lower semi-continuous} at 
%% the point 
$\mathbf{x}_{0} \!\in\! \mathcal{D}$ if for every $\varepsilon \!>\! 0$ there is a radius $r \!>\! 0$ such that 
$f(\mathbf{x}) \geq f(\mathbf{x}_{0}) - \varepsilon$ for all $\mathbf{x} \rmv\in\rmv \mathcal{B}(\mathbf{x}_{0},r)$. (This definition is equivalent to
$\liminf_{\mathbf{x} \rightarrow \mathbf{x}_{0}} f(\mathbf{x}) \geq f(\mathbf{x}_{0})$, where $\liminf_{\mathbf{x} \rightarrow \mathbf{x}_{0}} f(\mathbf{x}) 
\triangleq \sup_{r >0} \big\{ \rmv\rmv \inf_{\mathbf{x} \,\in\, \mathcal{D} \,\cap\, [\mathcal{B}(\mathbf{x}_{0},r) \setminus \{ \mathbf{x}_{0}\}]} f(\mathbf{x}) \big\}$
\cite{CounterExamplesAnalysis,RudinBookPrinciplesMatheAnalysis}.)
%Given a real-valued function $R(\cdot,\cdot): \mathcal{D} \times \mathcal{D} \rightarrow \mathbb{R}$ and 
%an element $\mathbf{x} \in \mathcal{D}$, we denote by $R(\cdot, \mathbf{x})$ the specific function $f(\cdot): \mathcal{D} \rightarrow \mathbb{R}$ given by $f(\mathbf{x}') = R(\mathbf{x}', \mathbf{x})$ for every $\mathbf{x}' \in \mathcal{D}$. 
%Given two functions $f(\cdot): \mathcal{A} \rightarrow \mathcal{B}$, $g(\cdot): \mathcal{C} \rightarrow \mathcal{D}$, we denote by $f(\cdot)=g(\cdot)$ the fact that $\mathcal{A} = \mathcal{C}$, $\mathcal{B} = \mathcal{D}$ and 
%for every $\mathbf{x} \inÊ\mathcal{A}$ we have $f(\mathbf{x}) = g(\mathbf{x})$. 
The restriction of a function $f(\cdot) \!: \mathcal{D} \rightarrow \mathbb{R}$ 
%% with domain $\mathcal{D}$
to a subdomain $\mathcal{D}' \!\subseteq\! \mathcal{D}$ is denoted by $f(\cdot)\big|_{\mathcal{D}'}$. 
Given a multi-index $\mathbf{p} = (p_1 \cdots p_N)^T \!\in \mathbb{Z}_{+}^{N}$, we define the partial derivative of order $\mathbf{p}$ 
of a real-valued function $f(\cdot): \mathcal{D} \rightarrow \mathbb{R}$, with $\mathcal{D} \subseteq \mathbb{R}^{N}\rmv$, as 
$\frac{\partial^{\mathbf{p}} \! f(\mathbf{x})}{\partial \mathbf{x}^{\mathbf{p}}} \triangleq \frac{\partial^{p_{1}}}{\partial x_{k}^{p_{1}}} \cdots \frac{\partial^{p_{N}}}{\partial x_{N}^{p_{N}}}  f(\mathbf{x})$ (if it exists) \cite{KranzPrimerAnalytic,RudinBookPrinciplesMatheAnalysis}.
Similarly, for a function $f(\cdot\ist\ist,\cdot) : \mathcal{D} \times \mathcal{D} \rightarrow \mathbb{R}$ and two multi-indices $\mathbf{p}_{1}, \mathbf{p}_{2} \in \mathbb{Z}_{+}^{N}$, we denote   
by $\frac{\partial^{\mathbf{p}_{1}} \partial^{\mathbf{p}_{2}} f(\mathbf{x}_{1},\mathbf{x}_{2})}{\partial \mathbf{x}_{1}^{\mathbf{p}_{1}}\mathbf{x}_{2}^{\mathbf{p}_{2}}}$ the partial derivative 
of order $(\mathbf{p}_{1},\mathbf{p}_{2})$, where $f(\mathbf{x}_{1}, \mathbf{x}_{2})$ is considered as a 
%% real-valued 
function of the ``super-vector'' $(\mathbf{x}_{1}^T \,\mathbf{x}_{2}^T)^T\rmv$ of length $2N$.
Given a vector-valued function ${\bm \phi}(\cdot):\mathbb{R}^{M} \!\rightarrow \mathbb{R}^{N}\rmv$ 
and $\mathbf{p} 
%% = (p_1 \cdots p_N)^T \!
\in \mathbb{Z}_{+}^{N}$, we denote the product $\prod_{k=1}^{N} \rmv\big( \phi_{k}(\mathbf{y}) \big)^{p_{k}}$ by ${\bm \phi}^{\mathbf{p}}(\mathbf{y})$.

The probability measure of a random vector $\mathbf{y}$ taking on values in $\mathbb{R}^{M}$ 
%% for some $M \rmv\in\rmv \mathbb{N}$,
is denoted by $\mu^{\mathbf{y}}$ \cite{BillingsleyProbMeasure,LC,HalmosRadonNikodymSuffStat,IbragimovBook}. 
We 
%% only 
consider probability measures that are defined on the measure space 
given by all $M$-dimensional Borel sets on $\mathbb{R}^{M}$ \cite[Sec. 10]{BillingsleyProbMeasure}. The probability measure assigns to a measureable set $\mathcal{A} \subseteq \mathbb{R}^{M}$ the probability 
\[ 
%% \label{equ_def_prob}
\Prob \{ \mathbf{y} \rmv\in\rmv \mathcal{A} \} \,\triangleq \int_{\mathbb{R}^{M}} \!I_{\mathcal{A}}(\mathbf{y}') \,d \mu^{\mathbf{y}}(\mathbf{y}') 
\,=  \int_{\mathcal{A}}d \mu^{\mathbf{y}}(\mathbf{y}') \,,
%% \vspace{-.5mm}
\] 
where $I_{\mathcal{A}}(\cdot) \!: \mathbb{R}^{M} \!\rightarrow\rmv \{0,1\}$ denotes the indicator function of the set $\mathcal{A}$.
We will also consider a family of probability measures ${\{ \mu^{\mathbf{y}}_{\mathbf{x}} \}}_{\mathbf{x} \in \mathcal{X}}$ 
parametrized by a nonrandom parameter vector $\mathbf{x} \rmv\in\rmv \mathcal{X}$. We assume that there exists a dominating measure 
$\mu_{\mathcal{E}}$, so that we can define the pdf $f(\mathbf{y}; \mathbf{x})$ (again parametrized by $\mathbf{x}$) 
as the Radon-Nikodym derivative of the measure 
$\mu^{\mathbf{y}}_{\mathbf{x}}$ with respect to the measure $ \mu_{\mathcal{E}}$ \cite{BillingsleyProbMeasure,LC,HalmosRadonNikodymSuffStat,IbragimovBook}. 
(In general, we will choose for $\mu_{\mathcal{E}}$ the Lebesgue measure on $\mathbb{R}^{M}\rmv$.)
We refer to both the set of measures ${\{ \mu^{\mathbf{y}}_{\mathbf{x}} \}}_{\mathbf{x} \in \mathcal{X}}$ and the set of pdfs 
${\{ f(\mathbf{y}; \mathbf{x}) \}}_{\mathbf{x} \in \mathcal{X}}$
%% , which are obtained for every parameter vector $\mathbf{x} \in \mathcal{X}$, 
as the \emph{statistical model}. 
Given a (possibly vector-valued) deterministic function $\mathbf{t}(\mathbf{y})$, the expectation operation 
%% based on the underlying probability measure $\mu^{\mathbf{y}}_{\mathbf{x}}$ \cite{BillingsleyProbMeasure}
is 
%% denoted by $\expect_{\mathbf{x}}Ê\{ \cdot \}$ and 
defined by \cite{BillingsleyProbMeasure} 
\begin{equation}
\label{equ_def_expect}
\expect_{\mathbf{x}}Ê\{ \mathbf{t}(\mathbf{y}) \} \,\triangleq \int_{\mathbb{R}^{M}} \! \mathbf{t}(\mathbf{y}') \,d \mu^{\mathbf{y}}_{\mathbf{x}}(\mathbf{y}') 
\,= \int_{\mathbb{R}^{M}} \! \mathbf{t}(\mathbf{y}') \, f(\mathbf{y}';\mathbf{x}) \,d \mathbf{y}', 
%% \vspace{-.3mm}
\end{equation} 
where the subscript 
%% $_{\mathbf{x}$ 
in $\expect_{\mathbf{x}}$ indicates the dependence on the parameter vector $\mathbf{x}$ parametrizing $\mu^{\mathbf{y}}_{\mathbf{x}}(\mathbf{y})$ and $f(\mathbf{y};\mathbf{x})$.

\vspace{-1mm}

%%%%%%%%%%%%%%%%%%%%%%%%%%%%%%%%%%%%%%%%%%%%%%%%%%%%%%%%%%
\textcolor{red}{
\section{Fundamentals}
\label{SecMVE}
%%%%%%%%%%%%%%%%%%%%%%%%%%%%%%%%%%%%%%%%%%%%%%%%%%%%%%%%%%
\subsection{Review of MVE}}
\setcounter{paragraph}{0}
%% \subsection{Basic Concepts} 
%% \label{sec_basic_concepts}
%Short discussion about the formal framework of MVE ---> Estimation problems, parameter functions, bias functions, equivalence bias and parameter function, 
%vector valued parameter functions can be reduced to scalar parameter functions (since the variances add up) 
%% \vspace*{5mm}
%% \paragraph{The Notion of an Estimation Problem}

%We consider the 
%% following ``classical'' or ``frequentist'' estimation problem. We wish to estimate 
%estimation of a function value $\mathbf{g}(\mathbf{x})$ from an observed vector $\mathbf{y}$,
%where the deterministic \emph{parameter vector} $\mathbf{x} \!\in\! \mathcal{X} \!\subseteq\! \mathbb{R}^{N}$ is unknown
%except for the fact that it belongs to a known \emph{parameter set} $\mathcal{X}$, and the deterministic \emph{parameter function} 
%$\mathbf{g}(\cdot) \!: \mathcal{X} \!\rightarrow \mathbb{R}^{P}$ 
%% (where the dimension $P$ of the parameter function is assumed deterministic and known) is assumed known and deterministic. .
%is known. Furthermore, the random \emph{observation} $\mathbf{y} \rmv\in\rmv \mathbb{R}^{M}$ 
%is distributed according to the parametrized set of pdfs (the \emph{statistical model}) ${\{ f(\mathbf{y}; \mathbf{x}) \}}_{\mathbf{x} \in \mathcal{X}}$.
\textcolor{red}{
It will be convenient to denote a classical (frequentist) estimation problem by the triple 
$\mathcal{E} = \big( Ê\mathcal{X}, f(\mathbf{y}; \mathbf{x}), \mathbf{g}(\cdot) \big)$, consisting of the parameter set $\mathcal{X}$, 
the statistical model ${\{ f(\mathbf{y}; \mathbf{x}) \}}_{\mathbf{x} \in \mathcal{X}}$, and the parameter function $\mathbf{g}(\cdot) \!: \mathcal{X} \!\rightarrow \mathbb{R}^{P}$.}
Note that our setting includes estimation of the parameter vector $\mathbf{x}$ itself, which is obtained when $\mathbf{g}(\mathbf{x}) = \mathbf{x}$.
The result of estimating $\mathbf{g}(\mathbf{x})$ from $\mathbf{y}$ is an \emph{estimate} $\hat{\mathbf{g}} \rmv\in\rmv \mathbb{R}^{P}\!$,
which is derived from $\mathbf{y}$ via a deterministic \emph{estimator} $\hat{\mathbf{g}}(\cdot) \!: \mathbb{R}^{M} \!\rightarrow \mathbb{R}^{P}\!$,
i.e., $\hat{\mathbf{g}} = \hat{\mathbf{g}}(\mathbf{y})$. We assume that any estimator is a measurable mapping from $\mathbb{R}^{M}$ to $\mathbb{R}^{P}$ \cite[Sec.\ 13]{BillingsleyProbMeasure}. 
%% The general goal in the design of
%% when constructing the 
%% an estimator $\hat{\mathbf{g}}(\cdot)$ is that $\hat{\mathbf{g}}(\mathbf{y})$ be close to the true value $\mathbf{g}(\mathbf{x})$. 
A 
%% reasonable and 
convenient
%% popular 
characterization of the performance of an estimator $\hat{\mathbf{g}}(\cdot)$
%% for assessing the quality of a specific estimator $\hat{\mathbf{g}}(\mathbf{y})$ 
is the \emph{mean squared error} (MSE) 
%% $\varepsilon$ 
defined as 
\[
%% \label{equ_def_MSE}
 \varepsilon \,\triangleq\, \expect_{\mathbf{x}} \big\{ {\|  \hat{\mathbf{g}}(\mathbf{y}) \rmv-\rmv \mathbf{g}(\mathbf{x})\|}^{2}_{2} \big\} 
 \,= \int_{\mathbb{R}^{M}} \! {\|  \hat{\mathbf{g}}(\mathbf{y}) \rmv-\rmv \mathbf{g}(\mathbf{x})\|}^{2}_{2} 
 \,\ist f(\mathbf{y}; \mathbf{x}) \, d \mathbf{y} \,. 
\]
We will write $\varepsilon(\hat{\mathbf{g}}(\cdot); \mathbf{x})$ to explicitly indicate the dependence of the MSE on the estimator $\hat{\mathbf{g}}(\cdot)$ and the parameter vector 
$\mathbf{x}$. Unfortunately, for a general estimation problem $\mathcal{E} = \big( Ê\mathcal{X}, f(\mathbf{y}; \mathbf{x}), \mathbf{g}(\cdot) \big)$, 
there does not exist an estimator $\hat{\mathbf{g}}(\cdot)$ that minimizes the MSE simultaneously for all parameter vectors $\mathbf{x} \rmv\in\rmv \mathcal{X}$ 
\cite{kay,RethinkingBiasedEldar}. 
This follows from the fact that minimizing the MSE at a given parameter vector $\mathbf{x}_{0}$ always yields zero MSE; this is achieved by the 
%% dumb 
estimator 
%% \begin{equation}
%% \label{equ_dumb_estimator}
$\hat{\mathbf{g}}_{0}(\mathbf{y}) = \mathbf{g}(\mathbf{x}_{0})$,
%% \end{equation}
which completely ignores the observation $\mathbf{y}$.

A popular rationale for the design of good estimators is MVE.
%% that of \emph{minimum variance estimation} (MVE). 
This approach is based on the MSE decomposition 
%% \emph{bias-variance decomposition} of the MSE which refers to the fact that for any estimator $\hat{\mathbf{g}}(\cdot)$, we have 
\vspace{-2mm}
\begin{equation} 
\label{equ_bias_variance_tradeoff}
\varepsilon(\hat{\mathbf{g}}(\cdot); \mathbf{x}) \,=\, {\| \mathbf{b}(\hat{\mathbf{g}}(\cdot); \mathbf{x}) \|}_2^{2} \ist\ist+\ist\ist v(\hat{\mathbf{g}}(\cdot); \mathbf{x}) \,, 
\end{equation}
with the \emph{estimator bias} $\mathbf{b}(\hat{\mathbf{g}}(\cdot); \mathbf{x}) \triangleq \expect_{\mathbf{x}} \{ \hat{\mathbf{g}}(\mathbf{y}) \} - \mathbf{g}(\mathbf{x})$ and 
the \emph{estimator variance} $v(\hat{\mathbf{g}}(\cdot); \mathbf{x}) \triangleq \expect_{\mathbf{x}} \big\{ \| \hat{\mathbf{g}}(\mathbf{y})$\linebreak %%%%%%%
$-\, {\expect_{\mathbf{x}} \{ \hat{\mathbf{g}}(\mathbf{y}) \} \|}^{2}_{2} \big\}$.
In MVE, one fixes the bias for all parameter vectors, i.e., 
$\mathbf{b}(\hat{\mathbf{g}}(\cdot); \mathbf{x}) \stackrel{!}{=} \mathbf{c}(\mathbf{x})$ for all $\mathbf{x} \!\in\! \mathcal{X}$, with a \emph{prescribed bias function} 
$\mathbf{c}(\cdot): \mathcal{X} \!\rightarrow \mathbb{R}^{P}\!$, and considers only estimators with the given bias.
%%  function. 
Note that fixing the estimator bias is equivalent to fixing the estimator mean, i.e., 
$\expect_{\mathbf{x}}\big\{ \hat{\mathbf{g}}(\mathbf{y}) \big\} \stackrel{!}{=} {\bm \gamma}(\mathbf{x})$ for all $\mathbf{x} \!\in\! \mathcal{X}$,
with the \emph{prescribed mean function} ${\bm \gamma}(\mathbf{x}) \triangleq \mathbf{c}(\mathbf{x}) + \mathbf{g}(\mathbf{x})$. 
The important special case of \emph{unbiased estimation} is obtained for $\mathbf{c}(\mathbf{x}) \equiv \mathbf{0}$ or equivalently 
${\bm \gamma}(\mathbf{x}) \equiv \mathbf{g}(\mathbf{x})$ for all $\mathbf{x} \!\in\! \mathcal{X}$.
Fixing the bias can be viewed
%% motivated conceptually 
as a kind of regularization of the set of considered estimators \cite{LC,RethinkingBiasedEldar}, 
because useless estimators like the estimator $\hat{\mathbf{g}}_{0}(\mathbf{y}) = \mathbf{g}(\mathbf{x}_{0})$ are excluded. 
Another justification for 
%% considering only estimators with a fixed bias 
fixing the bias is the fact that, if a large number 
%% $L$ 
of independent and identically distributed (i.i.d.)\ 
%% observations 
realizations ${\{ \mathbf{y}_{i} \}}_{i=1}^{L}$ of the vector $\mathbf{y}$ are observed,
%% realizations of the observation $\mathbf{y}$, 
then, under certain technical conditions, the bias term dominates in the decomposition \eqref{equ_bias_variance_tradeoff}. Thus, in that case, the MSE is small if and only if the 
%% prescribing a small or zero bias will lead to a small MSE; 
bias
%% function $\mathbf{c}(\mathbf{x})$ 
is small; this means that the estimator has to be effectively unbiased, i.e., $\mathbf{b}(\hat{\mathbf{g}}(\cdot); \mathbf{x}) \approx \mathbf{0}$ 
for all $\mathbf{x} \!\in\! \mathcal{X}$.

%% \newpage %%%%%%%%%

%% In what follows, we consider 
For a fixed ``reference'' parameter vector $\mathbf{x}_{0} \!\in\! \mathcal{X}$ and a prescribed bias function $\mathbf{c}(\cdot)$, 
%% For $\mathbf{x}_{0}$ and $\mathbf{c}(\cdot)$, 
we define the \emph{set of allowed estimators} 
\vspace{-1mm}
by 
\[
%% \label{equ_est_finite_var_prescr_bias}
\mathcal{A}(\mathbf{c}(\cdot),\mathbf{x}_{0}) \,\triangleq\, \big\{ \hat{\mathbf{g}}(\cdot) \,\big|\, 
v(\hat{\mathbf{g}}(\cdot);\mathbf{x}_{0}) < \infty \,, \, \mathbf{b}(\hat{\mathbf{g}}(\cdot);\mathbf{x}) = \mathbf{c}(\mathbf{x}) \,\, \forall \mathbf{x} \!\in\! \mathcal{X} \big\} \,. 
\]
We call a bias function $\mathbf{c}(\cdot)$ \emph{valid} for the estimation problem $\mathcal{E}=\vecestproblem$ at $\mathbf{x}_{0} \!\in\! \mathcal{X}$ 
if the set $\mathcal{A}(\mathbf{c}(\cdot),\mathbf{x}_{0})$ is nonempty. This means that there is at least one estimator $\hat{\mathbf{g}}(\cdot)$ 
with finite variance at $\mathbf{x}_{0}$ and whose bias 
%% function 
equals $\mathbf{c}(\cdot)$, i.e., 
$\mathbf{b}(\hat{\mathbf{g}}(\cdot);\mathbf{x}) = \mathbf{c}(\mathbf{x})$ for all $\mathbf{x} \!\in\! \mathcal{X}$. 
From \eqref{equ_bias_variance_tradeoff}, it follows that for a fixed bias $\mathbf{c}(\cdot)$, minimizing the MSE $\varepsilon(\hat{\mathbf{g}}(\cdot); \mathbf{x}_{0})$ 
is equivalent to minimizing the variance $v(\hat{\mathbf{g}}(\cdot); \mathbf{x}_{0})$. Therefore, in MVE, one attempts to find estimators that minimize the variance under
the constraint of a prescribed bias $\mathbf{c}(\cdot)$ function. Let
%% The minimum variance is
%% the \emph{minimum variance problem} (MVP)
\vspace*{-2mm}
\begin{equation}
\label{equ_minimum_variance_problem}
M(\mathbf{c}(\cdot),\mathbf{x}_{0}) \,\triangleq \inf_{\hat{\mathbf{g}}(\cdot) \ist\in\ist \mathcal{A}(\mathbf{c}(\cdot),\mathbf{x}_{0})} \rmv v(\hat{\mathbf{g}}(\cdot); \mathbf{x}_{0})
\vspace{1mm}
\end{equation}
denote the minimum (strictly speaking, infimum) variance at $\mathbf{x}_{0}$ for bias function $\mathbf{c}(\cdot)$.
If $\mathcal{A}(\mathbf{c}(\cdot),\mathbf{x}_{0})$ is empty, i.e., if $\mathbf{c}(\cdot)$ is not valid, we set $M(\mathbf{c}(\cdot),\mathbf{x}_{0}) \triangleq \infty$.
Any estimator $\hat{\mathbf{g}}^{(\mathbf{x}_{0})}(\cdot)  \in \mathcal{A}(\mathbf{c}(\cdot),\mathbf{x}_{0})$ that achieves the infimum in \eqref{equ_minimum_variance_problem}, i.e., 
for which 
%% \[
%% %% \label{equ_def_LMV_estimator}
$v\big(\hat{\mathbf{g}}^{(\mathbf{x}_{0})}(\cdot);\mathbf{x}_{0}\big) = M(\mathbf{c}(\cdot),\mathbf{x}_{0})$,
%% \]
is called a \emph{locally minimum variance} (LMV) estimator at $\mathbf{x}_{0}$
%%  \rmv\in\rmv \mathcal{X}$ 
for bias function $\mathbf{c}(\cdot)$ \cite{LC,Parzen59,Duttweiler73b}. 
The corresponding minimum variance $M(\mathbf{c}(\cdot),\mathbf{x}_{0})$ is called the \emph{minimum achievable variance} at $\mathbf{x}_{0}$
%%  \rmv\in\rmv \mathcal{X}$  
for bias function $\mathbf{c}(\cdot)$. The minimization problem \eqref{equ_minimum_variance_problem} is referred to as a \emph{minimum variance problem} (MVP).
By its definition in \eqref{equ_minimum_variance_problem}, 
%% the minimum achievable variance 
$M(\mathbf{c}(\cdot),\mathbf{x}_{0})$ is a lower bound on the variance at $\mathbf{x}_{0}$ of any estimator with bias function $\mathbf{c}(\cdot)$, i.e., 
\begin{equation} 
\label{equ_lower_vound_variance_trivial_min_achiev_var}
\hat{\mathbf{g}}(\cdot) \in \mathcal{A}(\mathbf{c}(\cdot),\mathbf{x}_{0}) \;\;\Rightarrow\;\; v(\hat{\mathbf{g}}(\cdot);\mathbf{x}_{0}) \ist\geq\ist M(\mathbf{c}(\cdot),\mathbf{x}_{0}) \,.
\end{equation}
In fact, $M(\mathbf{c}(\cdot),\mathbf{x}_{0})$ is the tightest lower bound,
%%  on the variance of estimators with bias $\mathbf{c}(\cdot)$, 
which is sometimes referred to as the \emph{Barankin bound}.
%for any estimator $\hat{\mathbf{g}}(\cdot)$ such that $\mathbf{b}(\hat{\mathbf{g}}(\cdot),\mathbf{x}) = \mathbf{c}(\mathbf{x})$ for every $\mathbf{x} \in \mathcal{X}$. 

If, for a prescribed bias function $\mathbf{c}(\cdot)$, there exists an estimator that is the LMV estimator \emph{simultaneously} at all $\mathbf{x}_0 \rmv\in\rmv \mathcal{X}$, 
then that estimator is called the \emph{uniformly minimum variance} (UMV) estimator for bias function $\mathbf{c}(\cdot)$ \cite{LC,Parzen59,Duttweiler73b}. 
%Note that in this paper, we consider only estimation problems where the existence of a LMV is guaranteed, while there does not exist a UMV estimator in general. 
For many estimation problems, a UMV estimator does not exist. However, it always exists
%% a UMV estimator can always be found 
if there exists a \emph{complete sufficient statistic} 
%% (see, e.g., 
\cite[Theorem 1.11 and Corollary 1.12]{LC}, \cite[Theorem 6.2.25]{CasellaBergerStatInf}. Under mild conditions, this includes the case where the statistical model 
corresponds to an exponential family.
%% , one can show that under mild conditions there exists an UMV estimator.

The variance to be minimized can be decomposed as
\[
%% \label{equ_sum_var_scalar}
v(\hat{\mathbf{g}}(\cdot); \mathbf{x}_{0}) \,= \sum_{l \in [P]} \! v(\hat{g}_{l}(\cdot); \mathbf{x}_{0}) \,,
\] 
where $\hat{g}_{l}(\cdot) \triangleq \big(\hat{\mathbf{g}}(\cdot)\big)_l$
and $v(\hat{g}_{l}(\cdot); \mathbf{x}_{0}) \triangleq \expect_{\mathbf{x}} \big\{ \big[ \hat{g}_{l}(\mathbf{y}) - \expect_{\mathbf{x}} \{ \hat{g}_{l}(\mathbf{y}) \} \big]^{2} \big\}$
for $l \rmv\in\rmv [P]$. Moreover, 
%% \[
%%  \hat{\mathbf{g}}(\cdot) \in \mathcal{A}(\mathbf{c}(\cdot),\mathbf{x}_{0}) \quad\;Ê\Leftrightarrow \quad\; \hat{g}_{l}(\cdot) \in \mathcal{A}(c_{l}(\cdot),\mathbf{x}_{0})  
%%    \;\mbox{ for all }  k \rmv\in\rmv [P] \,.
%% \]
$\hat{\mathbf{g}}(\cdot) \in \mathcal{A}(\mathbf{c}(\cdot),\mathbf{x}_{0})$ if and only if $\hat{g}_{l}(\cdot) \in \mathcal{A}(c_{l}(\cdot),\mathbf{x}_{0})$ for all 
$l \rmv\in\rmv [P]$, where $c_{l}(\cdot) \triangleq \big(\mathbf{c}(\cdot)\big)_l$.
It follows that the minimization of $v(\hat{\mathbf{g}}(\cdot); \mathbf{x}_{0})$ can be reduced to $P$ separate problems of minimizing the component variances
$v(\hat{g}_{l}(\cdot); \mathbf{x}_{0})$, each involving the optimization of a single scalar component $\hat{g}_{l}(\cdot)$ of $\hat{\mathbf{g}}(\cdot)$
subject to the scalar bias constraint $b(\hat{g}_{l}(\cdot); \mathbf{x}) = c_{l}(\mathbf{x})$ for all $\mathbf{x} \!\in\! \mathcal{X}$.  
Therefore, without loss of generality, we will hereafter assume that the parameter function $\mathbf{g}(\mathbf{x})$ is scalar-valued, i.e., $P \!=\! 1$. 

\vspace{-2mm}

\textcolor{red}{
%%%%%%%%%%%%%%%%%%%%%%%%%%%%%%%%%%%%%%%%%%%%%%%%%%%%%%%%%%
\subsection{Review of the RKHS Approach to MVE}}
\label{SecRKHSMVE}
%%%%%%%%%%%%%%%%%%%%%%%%%%%%%%%%%%%%%%%%%%%%%%%%%%%%%%%%%%

\vspace{.5mm}

A powerful mathematical toolbox for MVE 
%% and the MVP \eqref{equ_minimum_variance_problem} 
is provided by RKHS theory \cite{Parzen59,Duttweiler73b,aronszajn1950}. 
In this subsection, we review basic definitions and results of RKHS theory and its application to MVE, and we discuss a differentiability property that will be relevant to
the 
%% lower 
variance bounds considered in Section \ref{sec_rkhs_interpret_well_known_bounds}. 
%% Furthermore, we introduce the concepts of a regular estimation problem and of a differentiable RKHS and show their equivalence.

An RKHS is associated with a \emph{kernel function}, which is a function $R(\cdot\ist\ist,\cdot) \!: \mathcal{X} \rmv\times\rmv \mathcal{X} \rightarrow \mathbb{R}$ 
with the following two properties \cite{aronszajn1950}:
%% \label{def_kernel_function}

\begin{itemize}

\item It is symmetric, i.e., $R(\mathbf{x}_{1}, \mathbf{x}_{2}) = R(\mathbf{x}_{2}, \mathbf{x}_{1})$ for all $\mathbf{x}_{1}, \mathbf{x}_{2} \in \mathcal{X}$.

\item For every finite set 
%% $\mathcal{D}=
$\{\mathbf{x}_{1},\ldots,\mathbf{x}_{D}\} \subseteq \mathcal{X}$, the matrix $\mathbf{R} \in \mathbb{R}^{D \times D}$ with entries
%% defined elementwise as 
$R_{m,n} = R(\mathbf{x}_{m},\mathbf{x}_{n})$ is positive semidefinite.

\end{itemize}

\noindent 
%% In what follows, we will use the shorthand notation $R(\cdot\ist,\mathbf{x})$ for the function $f_{\mathbf{x}}(\mathbf{x}') = R(\mathbf{x}'\!, \mathbf{x})$ 
%% with a fixed $\mathbf{x} \rmv\in\rmv \mathcal{X}$.
%% It is a fundamental result \cite{aronszajn1950} that there 
There exists an RKHS for any kernel function $R(\cdot\ist\ist,\cdot)\!: \mathcal{X} \rmv\times\rmv \mathcal{X} \rightarrow \mathbb{R}$ \cite{aronszajn1950}.
%% More specifically, this 
This RKHS, denoted $\mathcal{H}(R)$, is a Hilbert space equipped with an inner product 
${\langle \cdot\ist\ist, \cdot \rangle}_{\mathcal{H}(R)}$ such that, for any $\mathbf{x} \!\in\! \mathcal{X}$, 

\begin{itemize}

\item  
%% the function $f_{\mathbf{x}}(\mathbf{x}') = R(\mathbf{x}', \mathbf{x})$, for which we will use the shorthand $R(\cdot\ist\ist,\mathbf{x})$ in the following, belongs to $\mathcal{H}(R)$, i.e., 
%% \[
$R(\cdot\ist,\mathbf{x}) \rmv\in\rmv \mathcal{H}(R)$ (here, $R(\cdot\ist,\mathbf{x})$ denotes the function $f_{\mathbf{x}}(\mathbf{x}') = R(\mathbf{x}'\!, \mathbf{x})$ 
with a fixed $\mathbf{x} \rmv\in\rmv \mathcal{X}$);
%% \] 

\item for any function 
\vspace{-2mm}
$f(\cdot) \rmv\in\rmv \mathcal{H}(R)$,
%%  and $\mathbf{x} \in \mathcal{X}$, we have 
\begin{equation} 
\label{equ_reproducing_property}
\big\langle f(\cdot), R(\cdot\ist\ist,\mathbf{x}) \big\rangle_{\mathcal{H}(R)} = f(\mathbf{x}) \,.
\end{equation}
%% (This is known as the \emph{reproducing property}.)

\end{itemize}

\noindent Relation \eqref{equ_reproducing_property}, which is known as the \emph{reproducing property}, 
defines the inner product ${\langle f, g \rangle}_{\mathcal{H}(R)}$ 
for all $f(\cdot), g(\cdot) \rmv\in \mathcal{H}(R)$ because (in a certain sense) any $f(\cdot) \rmv\in\rmv \mathcal{H}(R)$ can be expanded into the 
set of functions ${\{ R(\cdot\ist,\mathbf{x}) \}}_{\mathbf{x} \in \mathcal{X}}$. 
\textcolor{red}{In particular, consider two functions $f(\cdot), g(\cdot) \in \mathcal{H}(R)$ that are given as 
$f(\cdot) = \sum_{\mathbf{x}_{k} \in \mathcal{D}} a_{k}  R(\cdot\ist,\mathbf{x}_{k})$ and 
$g(\cdot) = \sum_{\mathbf{x}'_{l} \in \mathcal{D}'} b_{l}  R(\cdot\ist,\mathbf{x}'_{l})$ 
with coefficients $a_{k}, b_{l} \in \mathbb{R}$ and (possibly infinite) sets $\mathcal{D}, \mathcal{D}' \subseteq \mathcal{X}$. Then, by the linearity of inner products and \eqref{equ_reproducing_property}, 
\begin{equation}
\big\langle f(\cdot), g(\cdot) \big\rangle_{\mathcal{H}(R)} 
= \sum_{\mathbf{x}_{k} \in \mathcal{D}} \sum_{\mathbf{x}'_{l} \in \mathcal{D}'} \!a_{k} b_{l} R(\mathbf{x}_{k},\mathbf{x}'_{l}) \,. \nonumber 
\end{equation} }

%% \newpage %%%%%%%

\vspace{-2mm}

\subsubsection{\textcolor{red}{The RKHS Associated with an MVP}}
\label{SecRKHS-basics}
%%%%%%%%%%%%%%%%%%%%%%%%%%%%%%%%%%%%%%%%%%%%%%%%%%%%%%%%%%

Consider the class of MVPs that is defined by an estimation problem $\mathcal{E}=\big( \mathcal{X}, f(\mathbf{y}; \mathbf{x}), g(\cdot) \big)$,
a reference parameter vector $\mathbf{x}_{0}\!\in\! \mathcal{X}$, and all possible prescribed bias functions $c(\cdot) \!: \mathcal{X} \rightarrow \mathbb{R}$.
With this class of MVPs, we can associate a kernel function $R_{\mathcal{E}\rmv,\mathbf{x}_{0}}(\cdot\ist\ist,\cdot) \!: \mathcal{X} \rmv\times\rmv \mathcal{X} \rightarrow \mathbb{R}$
and, in turn, an RKHS $\mathcal{H}(R_{\mathcal{E}\rmv,\mathbf{x}_{0}})$ \cite{Parzen59,Duttweiler73b}. (Note that, as our notation indicates, 
$R_{\mathcal{E}\rmv,\mathbf{x}_{0}}(\cdot\ist\ist,\cdot)$ and $\mathcal{H}(R_{\mathcal{E}\rmv,\mathbf{x}_{0}})$ depend on $\mathcal{E}$ and $\mathbf{x}_{0}$ but not on $c(\cdot)$.)
%% For this, we 
\textcolor{red}{
We assume that 
%% for all reference parameters $\mathbf{x}_{0} \!\in\! \mathcal{X}$ for which the MVP \eqref{equ_minimum_variance_problem} is considered, 
%% the pdf $f(\mathbf{y};\mathbf{x}_{0})$ does not vanish anywhere, i.e., 
\be
\label{equ_f_y_x_0_not_vanish}
\Prob \{ f(\mathbf{y}; \mathbf{x}_{0}) \neq 0 \} = 1 \,, 
\ee
where the probability is evaluated for the underlying dominating measure $\mu_{\mathcal{E}}$.
We can then define the \emph{likelihood ratio} 
%% \vspace{-1mm}
\pagebreak %%%%%%%
as
\begin{equation} 
\label{equ_def_likelihood}
\llr (\mathbf{y},\mathbf{x}) \,\triangleq\, \begin{cases} \displaystyle\frac{f(\mathbf{y};\mathbf{x})}{f(\mathbf{y};\mathbf{x}_{0})} \,, &\mbox{if}\; f(\mathbf{y};\mathbf{x}_{0}) \neq 0 \\ 
0 \,, &\mbox{else.} \end{cases}
\vspace{-3mm}
\end{equation}}

\noindent
We consider $\llr (\mathbf{y},\mathbf{x})$ as a random variable (since it is a function of the random vector $\mathbf{y}$) that is parametrized by $\mathbf{x} \rmv\in\rmv \mathcal{X}$. 
Furthermore, we define the Hilbert space $\mathcal{L}_{\mathcal{E}\rmv,\mathbf{x}_{0}}$ as the closure of the linear 
span\footnote{A %%%%%%%%%
detailed discussion of the concepts of closure, inner product, orthonormal basis, and linear span in the context of abstract 
Hilbert space theory can be found in 
\cite{Parzen59,RudinBook}.} %%%%%%%%%% 
of the set of random variables $\big\{ \llr (\mathbf{y},\mathbf{x}) \big\}_{\mathbf{x} 
\in \mathcal{X}}$. The topology of $\mathcal{L}_{\mathcal{E}\rmv,\mathbf{x}_{0}}$ is determined by the inner product 
${\langle \cdot \ist\ist, \cdot \rangle}_{\text{RV}} \!: \mathcal{L}_{\mathcal{E}\rmv,\mathbf{x}_{0}} \times \mathcal{L}_{\mathcal{E}\rmv,\mathbf{x}_{0}} \!\rightarrow \mathbb{R}$
defined by 
\begin{equation}
\label{equ_def_inn_prod_RV}
\big\langle  \llr (\mathbf{y}, \mathbf{x}_{1}) \ist, \llr (\mathbf{y}, \mathbf{x}_{2}) \big\rangle_{\text{RV}} 
\,\triangleq\, \expect_{\mathbf{x}_{0}} \big\{ \llr (\mathbf{y}, \mathbf{x}_{1}) \, \llr (\mathbf{y}, \mathbf{x}_{2}) \big\}.
\vspace{1mm}
\end{equation}
It can be shown that it is sufficient to define the inner product only for the random variables $\big\{ \llr(\mathbf{y}, \mathbf{x}) \big\}_{\mathbf{x} 
\in \mathcal{X}}$ \cite{Parzen59}. \textcolor{red}{We will assume that 
\begin{equation}
\label{equ_corr_likelihood_finite}
\big\langle  \llr (\mathbf{y}, \mathbf{x}_{1}) \ist, \llr (\mathbf{y}, \mathbf{x}_{2}) \big\rangle_{\text{RV}} < \infty \,, \quad \text{for all} \;\, 
\mathbf{x}_{1}, \mathbf{x}_{2} \in \mathcal{X} \,.
\vspace{1mm}
\end{equation}
The assumptions \eqref{equ_f_y_x_0_not_vanish} and \eqref{equ_corr_likelihood_finite} (or variants thereof)
%% the condition \eqref{equ_corr_likelihood_finite} for the analysis of MVE 
are standard in the literature on MVE 
%% (cf.\ 
\cite{Duttweiler73b,Barankin49,Parzen59,stein50}. 
%% Moreover, the conditions \eqref{equ_f_y_x_0_not_vanish}, \eqref{equ_corr_likelihood_finite} 
They are typically satisfied for the important and large class of estimation problems arising from exponential families (cf.\ Section \ref{sec_rkhs_exp_fam}).}

The inner product ${\langle \cdot \ist\ist, \cdot \rangle}_{\text{RV}} \!: \mathcal{L}_{\mathcal{E}\rmv,\mathbf{x}_{0}} \!\times\rmv \mathcal{L}_{\mathcal{E}\rmv,\mathbf{x}_{0}} \!\rightarrow \mathbb{R}$ 
can now be interpreted as a kernel function 
$R_{\mathcal{E}\rmv,\mathbf{x}_{0}}(\cdot\ist\ist,\cdot) \!: \mathcal{X} \!\times\! \mathcal{X} \rightarrow \mathbb{R}$:
\vspace*{-2mm}
\be
\label{equ_def_kernel_est_problem}
R_{\mathcal{E}\rmv,\mathbf{x}_{0}}(\mathbf{x}_{1},\mathbf{x}_{2}) \,\triangleq\, \big\langle \llr(\mathbf{y},\mathbf{x}_{1}) \ist, \llr (\mathbf{y}, \mathbf{x}_{2}) \big\rangle_{\text{RV}} 
\vspace{1mm}
\ee
The RKHS induced by $R_{\mathcal{E}\rmv,\mathbf{x}_{0}}(\cdot\ist\ist,\cdot)$ will be denoted by $\mathcal{H}_{\mathcal{E}\rmv,\mathbf{x}_{0}}$, i.e.,
$\mathcal{H}_{\mathcal{E}\rmv,\mathbf{x}_{0}} \triangleq \mathcal{H}(R_{\mathcal{E}\rmv,\mathbf{x}_{0}})$. This is the RKHS associated with the 
estimation problem $\mathcal{E}=\big( Ê\mathcal{X}, f(\mathbf{y}; \mathbf{x}), g(\cdot) \big)$ and the corresponding class of MVPs at $\mathbf{x}_{0}\rmv\in\rmv \mathcal{X}$.

We note that assumption \eqref{equ_f_y_x_0_not_vanish} implies that the likelihood ratio \textcolor{red}{$\llr (\mathbf{y},\mathbf{x})$} 
is measurable with respect to the underlying dominating measure $\mu_{\mathcal{E}}$.
%which in our case is the Lebesgue measure on $\mathbb{R}^{M}$ (cf.\ \cite[Theorem 13.3]{BillingsleyProbMeasure}). 
Furthermore, the likelihood ratio $\llr (\mathbf{y},\mathbf{x})$
%% $\frac{f(\mathbf{y};\mathbf{x})}{f(\mathbf{y};\mathbf{x}_{0})}$ 
is 
%% trivially 
the 
%% \emph{Radon-Nikodym derivative} 
Radon-Nikodym derivative \cite{BillingsleyProbMeasure,HalmosRadonNikodymSuffStat} of the probability measure $\mu_{\mathbf{x}}^{\mathbf{y}}$ 
induced by $f(\mathbf{y}; \mathbf{x})$ with respect to
the probability measure $\mu_{\mathbf{x}_{0}}^{\mathbf{y}}$ induced by 
$f(\mathbf{y}; \mathbf{x}_{0})$ (cf. \cite{RudinBook,BillingsleyProbMeasure,HalmosMeasure}).
It is also important to observe that $\llr (\mathbf{y},\mathbf{x})$ does not depend 
on the dominating measure $\mu_{\mathcal{E}}$ underlying the definition of the pdfs $f(\mathbf{y}; \mathbf{x})$.
Thus, the kernel $R_{\mathcal{E}\rmv,\mathbf{x}_{0}}(\cdot\ist\ist,\cdot)$ given by \eqref{equ_def_kernel_est_problem} 
does not depend on $\mu_{\mathcal{E}}$ either. Moreover, under assumption \eqref{equ_f_y_x_0_not_vanish}, 
we can always use the measure $\mu^{\mathbf{y}}_{\mathbf{x}_{0}}$ as the base
measure $\mu_{\mathcal{E}}$ for the estimation problem $\mathcal{E}$, since the Radon-Nikodym derivative 
\textcolor{red}{$\llr(\mathbf{y}, \mathbf{x})$} is well defined. Note that, trivially, this also implies that the measure $\mu^{\mathbf{y}}_{\mathbf{x}_{0}}$ dominates the measures ${\{ \mu^{\mathbf{y}}_{\mathbf{x}}Ê\}}_{\mathbf{x} \in \mathcal{X}}$ \cite[p. 443]{BillingsleyProbMeasure}.
%!!! Note that the existence of the Radon-Nikodym derivative of the measure $\mu^{\mathbf{y}}_{\mathbf{x}}$ w.r.t. 
%the measure $\mu^{\mathbf{y}}_{\mathbf{x}_{0}}$ trivially implies that the measure $\mu^{\mathbf{y}}_{\mathbf{x}_{0}}$ dominates $\mu^{\mathbf{y}}_{\mathbf{x}}$ (cf.\ \cite[pp. 216]{BillingsleyProbMeasure}).!!!

\textcolor{red}{ 
The two Hilbert spaces $\mathcal{L}_{\mathcal{E}\rmv,\mathbf{x}_{0}}$ and $\mathcal{H}_{\mathcal{E}\rmv,\mathbf{x}_{0}}$ are 
%% closely related. Indeed, consider an estimation problem $\mathcal{E}=\scalarestproblem$ and a fixed reference parameter vector 
%% $\mathbf{x}_{0} \rmv\in\rmv \mathcal{X}$. The  
%% associated 
%% Hilbert spaces $\mathcal{L}_{\mathcal{E}\rmv,\mathbf{x}_{0}}$ and $\mathcal{H}_{\mathcal{E}\rmv,\mathbf{x}_{0}}$ 
isometric. In fact, as proven in \cite{Parzen59}, 
a specific congruence (i.e., isometric mapping of functions in $\mathcal{H}_{\mathcal{E}\rmv,\mathbf{x}_{0}}$ to functions in $\mathcal{L}_{\mathcal{E}\rmv,\mathbf{x}_{0}}$) 
%% is constituted by the map  
$\mathsf{J}[\cdot] \!: \mathcal{H}_{\mathcal{E}\rmv,\mathbf{x}_{0}} \!\rightarrow \mathcal{L}_{\mathcal{E}\rmv,\mathbf{x}_{0}}$ is given by 
\[
%% \mathsf{J}[\cdot]: \, \mathcal{H}_{\mathcal{E}\rmv,\mathbf{x}_{0}} \rightarrow \mathcal{L}_{\mathcal{E}\rmv,\mathbf{x}_{0}} \!:\ist 
\mathsf{J}[ R_{\mathcal{E}\rmv,\mathbf{x}_{0}}(\cdot\ist,\mathbf{x})] = \llr (\cdot\ist,\mathbf{x}) \,.  
\vspace{-1mm}
\] 
%% A specific approach to evaluate 
The isometry $\mathsf{J}[f(\cdot)]$ can be evaluated for an arbitrary function $f(\cdot) \in \mathcal{H}_{\mathcal{E}\rmv,\mathbf{x}_{0}}$ by expanding 
$f(\cdot)$ into the elementary functions ${\{ R_{\mathcal{E}\rmv,\mathbf{x}_{0}}(\cdot,\mathbf{x}) \}}_{\mathbf{x} \in \mathcal{X}}$ (cf.\ \cite{Parzen59}). 
Given the expansion $f(\cdot) = \sum_{\mathbf{x}_{k} \in \mathcal{D}} 
a_{k}  R_{\mathcal{E}\rmv,\mathbf{x}_{0}}(\cdot\ist,\mathbf{x}_{k})$ with coefficients $a_{k} \in \mathbb{R}$ and a (possibly infinite) set $\mathcal{D} \subseteq \mathcal{X}$, the isometric image of $f(\cdot)$ is obtained as $\mathsf{J}[ f(\cdot)] =  \sum_{\mathbf{x}_{k} \in \mathcal{D}} a_{k} \, \llr(\cdot,\mathbf{x}_{k})$.}
%\end{theorem}
%% \begin{proof}
%% \cite{Parzen59} 
%% \end{proof}

\vspace{3.5mm}

\subsubsection{\textcolor{red}{RKHS-based Analysis of MVE}}
\label{SecRKHS-mve}
%%%%%%%%%%%%%%%%%%%%%%%%%%%%%%%%%%%%%%%%%%%%%%%%%%%%%%%%%%

%% \paragraph{RKHS based analysis of MVE}

%% PHDAJ -- p. 52 (page number appearing on the page)\\
\textcolor{red}{
An RKHS-based analysis of MVE is enabled by
%% based on the RKHS introduced in the previous subsection and 
the following central result.}
Consider an estimation problem $\mathcal{E}= \big( \mathcal{X}, f(\mathbf{y}; \mathbf{x}), g(\cdot) \big)$, 
a fixed reference parameter vector $\mathbf{x}_{0}Ê\rmv\in\rmv \mathcal{X}$, and 
a prescribed bias function $c(\cdot): \mathcal{X} \!\rightarrow \mathbb{R}$, corresponding to the prescribed mean function $\gamma(\cdot) \triangleq c(\cdot) + g(\cdot)$.
Then, \textcolor{red}{as shown in \cite{Parzen59,Duttweiler73b},} the following holds:
\begin{itemize} 
\item The bias function $c(\cdot)$ is valid for $\mathcal{E}$ at $\mathbf{x}_{0}$ if and only if 
%% the prescribed mean function 
$\gamma(\cdot)$ belongs to the RKHS $\mathcal{H}_{\mathcal{E}\rmv,\mathbf{x}_{0}}$, \textcolor{red}{i.e., 
\begin{equation}
\label{equ_nec_suff_cond_validity_bias_function}
\mathcal{A}(c(\cdot),\mathbf{x}_{0}) \neq \emptyset \;\;\Longleftrightarrow\;\; \gamma(\cdot) 
%% \triangleq c(\cdot)+g(\cdot)  
\in \mathcal{H}_{\mathcal{E}\rmv,\mathbf{x}_{0}} \,.
\vspace{-7.5mm}
\end{equation}
}
%% , i.e., 
%% \begin{equation}
%% c(\cdot) \mbox{ valid at } \mathbf{x}_{0} \mbox{ for } \mathcal{E} \Longleftrightarrow \gamma(\cdot) \in  \mathcal{H}_{\mathcal{E}\rmv,\mathbf{x}_{0}}.
%% \end{equation} 
\item If the bias function $c(\cdot)$ is valid, the corresponding minimum achievable variance at $\mathbf{x}_{0}$ is given by
\begin{equation}
\label{equ_min_achiev_var_sqared_norm}
M(c(\cdot),\mathbf{x}_{0}) \eq {\| \gamma(\cdot) \|}^{2}_{\mathcal{H}_{\mathcal{E}\rmv,\mathbf{x}_{0}}} \!- \gamma^{2}(\mathbf{x}_{0}) \,,
\end{equation}
and the LMV estimator at $\mathbf{x}_{0}$ is given 
\vspace{-3mm}
by 
\[
%% \label{equ_lmv_rkhs}
\hat{g}^{(\mathbf{x}_{0})}(\cdot) \ist=\ist \mathsf{J}[\gamma(\cdot)] \,.
\vspace{-1mm}
\]
\end{itemize}

%% \begin{proof}
%% \cite{Parzen59,Duttweiler73b} 
%% \end{proof}

This result
%% theorem 
shows that the RKHS $\mathcal{H}_{\mathcal{E},\mathbf{x}_{0}}$ is equal to 
the set of the mean functions $\gamma(\mathbf{x}) = \expect_{\mathbf{x}}Ê\{ \hat{g}(\mathbf{y}) \}$ of 
all estimators $\hat{g}(\cdot)$ with a finite variance at $\mathbf{x}_{0}$, i.e., $v(\hat{g}(\cdot);\mathbf{x}_{0}) < \infty$.
Furthermore, the problem of solving the MVP \eqref{equ_minimum_variance_problem} can be reduced to the computation 
of the squared norm ${\| \gamma(\cdot) \|}^{2}_{\mathcal{H}_{\mathcal{E}\rmv,\mathbf{x}_{0}}}$ and the isometric image $\mathsf{J}[\gamma(\cdot)]$ of the prescribed mean 
function $\gamma(\cdot)$, viewed as an element of the RKHS $\mathcal{H}_{\mathcal{E}\rmv,\mathbf{x}_{0}}$. This 
%% result 
is especially helpful if a simple characterization of 
%% the RKHS 
$\mathcal{H}_{\mathcal{E}\rmv,\mathbf{x}_{0}}$ is available. Here, following the terminology of \cite{Duttweiler73b}, what is meant by ``simple characterization'' is
%% the knowledge of 
the availability of an orthonormal basis (ONB) for $\mathcal{H}_{\mathcal{E}\rmv,\mathbf{x}_{0}}$ such that the inner products of 
%% the mean function 
$\gamma(\cdot)$ with
%% and 
the 
%% elements of the 
ONB functions can be computed easily.
%%  or in an efficient manner.

%% Even if 
If such an ONB of $\mathcal{H}_{\mathcal{E}\rmv,\mathbf{x}_{0}}$ cannot be found, the relation \eqref{equ_min_achiev_var_sqared_norm} can still be used 
to derive lower bounds on the minimum achievable variance $M(c(\cdot),\mathbf{x}_{0})$. Indeed, because of \eqref{equ_min_achiev_var_sqared_norm}, 
any lower bound on 
%% the squared norm 
${\| \gamma(\cdot) \|}^{2}_{\mathcal{H}_{\mathcal{E}\rmv,\mathbf{x}_{0}}}$ induces a lower bound on $M(c(\cdot),\mathbf{x}_{0})$. A large 
class of lower bounds on ${\| \gamma(\cdot) \|}^{2}_{\mathcal{H}_{\mathcal{E}\rmv,\mathbf{x}_{0}}}$ can be obtained via projections of 
%% the mean function 
$\gamma(\cdot)$ onto a
%% n arbitrary 
subspace $\mathcal{U} \subseteq \mathcal{H}_{\mathcal{E}\rmv,\mathbf{x}_{0}}$. 
Denoting the orthogonal projection of $\gamma(\cdot)$ onto 
%% the subspace 
$\mathcal{U}$ by 
$\gamma_{\mathcal{U}}(\cdot)$, we have 
${\| \gamma_{\mathcal{U}}(\cdot) \|}^{2}_{\mathcal{H}_{\mathcal{E}\rmv,\mathbf{x}_{0}}} \le {\| \gamma(\cdot) \|}^{2}_{\mathcal{H}_{\mathcal{E}\rmv,\mathbf{x}_{0}}}$ 
\cite[Chapter 4]{RudinBook} and thus, from \eqref{equ_min_achiev_var_sqared_norm},
\begin{equation}
\label{equ_lower_bound_min_achiev_var_projection}
 M(c(\cdot),\mathbf{x}_{0}) 
 %% = \| \gamma(\cdot) \|^{2}_{\mathcal{H}_{\mathcal{E}\rmv,\mathbf{x}_{0}}}- \gamma^{2}(\mathbf{x}_{0})  
 \,\geq\, {\| \gamma_{\mathcal{U}}(\cdot) \|}^{2}_{\mathcal{H}_{\mathcal{E}\rmv,\mathbf{x}_{0}}} \!- \gamma^{2}(\mathbf{x}_{0}) \,,
\end{equation}
for an arbitrary subspace $\mathcal{U} \subseteq \mathcal{H}_{\mathcal{E}\rmv,\mathbf{x}_{0}}$. 
%% The inequality in \eqref{equ_lower_bound_min_achiev_var_projection} is due to the Hilbert space projection theorem \cite[Chapter 4]{RudinBook}.
In particular, let us consider the special case of a finite-dimensional subspace $\mathcal{U} \subseteq \mathcal{H}_{\mathcal{E}\rmv,\mathbf{x}_{0}}$ 
that is spanned by a given set of functions $u_{l}(\cdot) \in \mathcal{H}_{\mathcal{E}\rmv,\mathbf{x}_{0}}$, i.e., 
\begin{equation} 
\label{equ_def_finit_dim_subspace_span} 
\mathcal{U} \ist=\, \linspan \rmv{\{ u_{l}(\cdot) \}}_{l \in [L]} \,\triangleq\ist \Bigg\{ f(\cdot) =\rmv \sum_{l \in [L]} a_{l} u_{l}(\cdot) \,\Bigg|\, a_{l} \!\in\rmv \mathbb{R} \Bigg\} \,.
\end{equation}
Here, 
%% the squared norm 
${\| \gamma_{\mathcal{U}}(\cdot) \|}^{2}_{\mathcal{H}_{\mathcal{E}\rmv,\mathbf{x}_{0}}}\!$ can be evaluated very easily 
due to the following expression
%% , which holds for any $f(\cdot) \in \mathcal{H}_{\mathcal{E}\rmv,\mathbf{x}_{0}}$ 
\cite[Theorem 3.1.8]{JungPHD}:
%% , which holds for arbitrary Hilbert spaces (not necessarily being a RKHS).
\begin{equation}
%% \label{equ_proof_lower_bound_exp_fam_squared_norm_projection_finite_dim_subspace}
\label{equ_projection_finite_dim_subspace}
{\| \gamma_{\mathcal{U}}(\cdot) \|}^{2}_{\mathcal{H}_{\mathcal{E}\rmv,\mathbf{x}_{0}}} \!=\, \bm{\gamma}^{T} \mathbf{G}^{\dagger} \bm{\gamma} \,,
\end{equation}
where the vector $\bm{\gamma} \rmv\in\rmv \mathbb{R}^{L}$ and the matrix $\mathbf{G} \rmv\in\rmv \mathbb{R}^{L \times L}$ are given elementwise by
\be
\label{equ_n_S_0}
\gamma_{l} \ist=\ist {\langle \gamma(\cdot), u_l(\cdot) \rangle}_{\mathcal{H}_{\mathcal{E}\rmv,\mathbf{x}_{0}}} \,, \qquad
G_{l,l'} \ist=\ist {\langle u_l(\cdot), u_{l'}(\cdot) \rangle}_{\mathcal{H}_{\mathcal{E}\rmv,\mathbf{x}_{0}}} \,.
\ee
If all $u_{l}(\cdot)$ are linearly independent, then a larger number $L$ of basis functions $u_{l}(\cdot)$ entails a higher dimension of $\mathcal{U}$ and, thus, a larger
${\| \gamma_{\mathcal{U}}(\cdot) \|}^{2}_{\mathcal{H}_{\mathcal{E}\rmv,\mathbf{x}_{0}}}\!$; this implies that the lower bound \eqref{equ_lower_bound_min_achiev_var_projection}
will be higher (i.e., tighter).
In Section \ref{sec_rkhs_interpret_well_known_bounds}, we will show that some well-known lower bounds on the estimator variance are obtained from \eqref{equ_lower_bound_min_achiev_var_projection} and \eqref{equ_projection_finite_dim_subspace}, 
using a subspace $\mathcal{U}$ of the form \eqref{equ_def_finit_dim_subspace_span} and specific choices for the functions $u_{l}(\cdot)$ spanning $\mathcal{U}$.

%% \paragraph{Estimation problems with Differentiable RKHS}

\vspace{3.5mm}

\subsubsection{\textcolor{red}{Regular Estimation Problems and Differentiable RKHS}}
\label{SecRKHS-diff}
%%%%%%%%%%%%%%%%%%%%%%%%%%%%%%%%%%%%%%%%%%%%%%%%%%%%%%%%%%

%% 
Some of the lower bounds to be considered in Section \ref{sec_rkhs_interpret_well_known_bounds} require the estimation problem to satisfy certain 
regularity conditions.
%%  summarized in 
%The validity of the CRB requires some weak technical conditions to be placed on the estimation problem $\mathcal{E}=\scalarestproblem$ and associated 
%MVP $\mathcal{M}=\left(\mathcal{E},c(\cdot)\equiv 0,\mathbf{x}_{0}\right)$, which are summarized in 

\vspace{1mm}

\begin{definition}
\label{post_regular_cond_CRB}

An estimation problem $\mathcal{E}=\scalarestproblem$ satisfying \eqref{equ_corr_likelihood_finite} is said to be 
\emph{regular up to order $m \in \mathbb{N}$ at an interior point $\mathbf{x}_{0} \rmv\in\rmv \mathcal{X}^{\text{o}}$} %???if there exists a 
%% sufficiently small 
%radius $r$ such that $\mathcal{B}(\mathbf{x}_{0},r) \subseteq \mathcal{X}$ 
if the following holds:

\vspace{.5mm}

\begin{itemize} 
%\item There exists a positive radius $r_{0} >0$ such that the parameter set $\mathcal{X} \subseteq \mathbb{R}^{N}$ contains the $N$-dimensional ball $\mathcal{B}(\mathbf{x}_{0},r_{0})$ centered at $\mathbf{x}_{0} \in \mathcal{X}$, i.e.,
%\begin{equation}
%\label{equ_unconstr_crb_param_set_contains_ball} 
%\mathcal{B}(\mathbf{x}_{0},r_{0}) \subseteq \mathcal{X}. 
%\end{equation} 

\vspace{1mm}

\item For every multi-index $\mathbf{p}Ê \in \mathbb{Z}_{+}^{N}$ with 
%% $\| \mathbf{p} \|_{\infty} \leq m$, 
entries $p_k \rmv\leq\rmv m$,
%%  for all $k \in \{1,\ldots,N\}$, 
the partial derivatives $ \frac{\partial^{\mathbf{p}} f(\mathbf{y}; \mathbf{x})}{\partial \mathbf{x}^{\mathbf{p}}}$ exist
and 
%% \vspace{-1mm}
%% moreover 
satisfy
\begin{equation} 
\label{equ_crb_regul_cond_finite_part_der}
\expect_{\mathbf{x}_{0}} \bigg\{ \! \left( \frac{1}{f(\mathbf{y};\mathbf{x}_{0})} \frac{\partial^{\mathbf{p}} f(\mathbf{y}; \mathbf{x})}{\partial \mathbf{x}^{\mathbf{p}}}Ê\right)^{\!\!2} \bigg\} < \infty \,,
\quad\; \text{for all} \;\; \mathbf{x} \in \mathcal{B}(\mathbf{x}_{0}, r) \,, 
%% \vspace{2.5mm}
\end{equation} 
where $r >0$ is a suitably chosen radius such that $\mathcal{B}(\mathbf{x}_{0},r) \subseteq \mathcal{X}$.
%% for all $\mathbf{x} \in \mathcal{B}(\mathbf{x}_{0}, r)$.
%and is continuous w.r.t. to $(\mathbf{y},\mathbf{x})$ where $(\mathbf{y}, \mathbf{x}) \in \mathbb{R}^{M} \times \mathcal{X}$. 

\vspace{1mm}

\item For any function $h(\cdot) \!: \mathbb{R}^{M} \!\rightarrow \mathbb{R}$ such that 
$\expect_{\mathbf{x}} \{ h(\mathbf{y}) \}$ exists,
%% \footnote{By definition, this implies that $\expect_{\mathbf{x}} \{ | h(\mathbf{y})| \} < \infty$ \cite{AshProbMeasure,HalmosMeasure, BillingsleyProbMeasure,Lapidoth09}.} 
%% we have that for all $\mathbf{x} \in \mathcal{B}(\mathbf{x}_{0},r)$ 
the expectation operation commutes with partial differentiation in the sense that, for every multi-index $\mathbf{p}Ê \in \mathbb{Z}_{+}^{N}$ with 
%% $\| \mathbf{p} \|_{\infty} \leq m$, 
\vspace{.5mm}
$p_k \rmv\leq\rmv m$,
\begin{equation} 
\label{equ_crb_regul_cond_intercah_diff_integr_intral_repr}
\frac{\partial^{\mathbf{p}}}{\partial \mathbf{x}^{\mathbf{p}}} \rmv\rmv \int_{\mathbb{R}^M} \! h(\mathbf{y}) \, f(\mathbf{y}; \mathbf{x}) \, d\mathbf{y}  
\,= \int_{\mathbb{R}^M} \! h(\mathbf{y}) \, \frac{\partial^{\mathbf{p}} f(\mathbf{y}; \mathbf{x})}{\partial \mathbf{x}^{\mathbf{p}}} \, d\mathbf{y}\,,
\quad\; \text{for all} \;\; \mathbf{x} \in \mathcal{B}(\mathbf{x}_{0}, r) \,,
\vspace{1mm}
\end{equation} 
or equivalently %expressed in terms of expectations, we have 
\begin{equation}
\label{equ_crb_regul_cond_intercah_diff_integr}
\frac{\partial^{\mathbf{p}} \ist\ist \expect_{\mathbf{x}} \{  h(\mathbf{y}) \}}{\partial \mathbf{x}^{\mathbf{p}}}
\eq \expect_{\mathbf{x}} \bigg\{  h(\mathbf{y}) \ist \frac{1}{f(\mathbf{y};\mathbf{x})} \ist \frac{\partial^{\mathbf{p}} f(\mathbf{y}; \mathbf{x})}{\partial \mathbf{x}^{\mathbf{p}}}  \bigg\}\,,
\quad\; \text{for all} \;\; \mathbf{x} \in \mathcal{B}(\mathbf{x}_{0}, r) \,,
\vspace{.5mm}
\end{equation}  
%%  for all $k \in \{1,\ldots,N\}$, 
provided that the right hand side of \eqref{equ_crb_regul_cond_intercah_diff_integr_intral_repr} and \eqref{equ_crb_regul_cond_intercah_diff_integr} is 
\vspace{1mm}
finite. 
%\item The parameter function $g(\cdot): \mathcal{X} \rightarrow \mathbb{R}$ is such that the partial derivatives $\frac{\partial^{\mathbf{p}} g(\mathbf{x})}{\partial \mathbf{x}^{\mathbf{p}}} \big|_{\mathbf{x} = \mathbf{x}_{0}}$ exist 
%at $\mathbf{x}_{0}$ for every multi-index $\mathbf{p}Ê \in \mathbb{Z}_{+}^{N}$ with $\| \mathbf{p} \|_{\infty} \leq m$.

\item For every pair of multi-indices $\mathbf{p}_{1}, \mathbf{p}_{2}Ê \in \mathbb{Z}_{+}^{N}$ with $p_{1,k} \rmv\leq\rmv m$ and $p_{2,k} \rmv\leq\rmv m$,
%% $\| \mathbf{p}_{1} \|_{\infty},\| \mathbf{p}_{2} \|_{\infty} \leq m$
the 
\vspace{.5mm}
expectation 
\begin{equation}
\label{equ_CRB_regul_cond_inner_prod_cont_x_1_x_2}
 \expect_{\mathbf{x}_{0}} \!\left\{ \frac{1}{f^2(\mathbf{y}; \mathbf{x}_{0})}  
 \ist\frac{\partial^{\mathbf{p}_{1}}f(\mathbf{y};\mathbf{x}_{1})} {\partial \mathbf{x}_{1}^{\mathbf{p}_{1}} } 
 %% \,\frac{1}{f(\mathbf{y}; \mathbf{x}_{0})} 
 \ist\frac{\partial^{\mathbf{p}_{2}} f(\mathbf{y};\mathbf{x}_{2})} {\partial \mathbf{x}_{2}^{\mathbf{p}_{2}} }  \right\}
\vspace{.5mm}
\end{equation} 
depends continuously on the parameter vectors 
$\mathbf{x}_1, \mathbf{x}_{2} \in \mathcal{B}(\mathbf{x}_{0}, r)$.

\end{itemize}

\end{definition}

\vspace{1mm}

%\noindent ???We note that condition \eqref{equ_crb_regul_cond_finite_part_der} and 
%condition \eqref{equ_crb_regul_cond_intercah_diff_integr_intral_repr}, \eqref{equ_crb_regul_cond_intercah_diff_integr} have been previously considered in \cite{} and 
%in \cite{}, respectively.

%% \newpage %%%%%%%

%% The use of Definition \ref{post_regular_cond_CRB} is to ensure that the RKHS associated to 
We remark 
that the notion of a \emph{regular estimation problem} according to Definition \ref{post_regular_cond_CRB} is somewhat similar to the notion of a \emph{regular statistical experiment} introduced in \cite[Section I.7]{IbragimovBook}.

\textcolor{red}{
As shown in \cite[Thm. 4.4.3.]{JungPHD}, the RKHS associated with a regular estimation problem has an important
%% appealing 
structural property, which we will term \emph{differentiable}. More precisely, we call an RKHS $\mathcal{H}(R)$ 
\emph{differentiable up to order $m$} if it is associated with a kernel $R(\cdot\ist\ist,\cdot): \mathcal{X} \rmv\times\rmv \mathcal{X} \!\rightarrow \mathbb{R}$ 
that is differentiable up to a given order $m$. The properties of differentiable RKHSs have been previously studied, e.g., in \cite{sun_jfaa,zhou_jfaa,ZhouTIT03}. 
}

%% \newpage %%%%%%%%

It will be seen that,
%% turn out that, 
under certain
%% specific 
conditions, the functions belonging to an RKHS $\mathcal{H}(R)$ that is differentiable are characterized completely by their partial derivatives 
\textcolor{red}{at any point $\mathbf{x}_0 \rmv\in\rmv \mathcal{X}^{\text{o}}$. This implies via \eqref{equ_nec_suff_cond_validity_bias_function} together with 
identity \eqref{equ_der_reproduction_prop} below}
%% Theorem \ref{thm_regul_cond_CRB_imply_diff} (see below) 
that, for a regular estimation problem, the mean function $\gamma(\mathbf{x})=\expect_{\mathbf{x}} \{ \hat{g}(\mathbf{y}) \}$ of any estimator $\hat{g}(\cdot)$ 
with finite variance at $\mathbf{x}_{0}$
%% , i.e., $v(\hat{g}(\cdot), \mathbf{x}_{0}) < \infty$, 
is completely specified by the partial derivatives $\big\{ \frac{\partial^{\mathbf{p}} \ist \gamma(\mathbf{x})}{\partial \mathbf{x}^{\mathbf{p}}} \big|_{\mathbf{x} = \mathbf{x}_{0}}Ê\big\}_{\mathbf{p} \in \mathbb{Z}_{+}^{N}}$ (cf.\ Lemma \ref{thm_any_valid_bias_func_analytic} in Section \ref{sec_param_set_reduction_exp_fam}).

\textcolor{red}{
Further important properties of a differentiable RKHS have been reported in \cite{zhou_jfaa,Kailath71}. In particular, for an RKHS $\mathcal{H}(R)$ that is differentiable up to order $m$, and 
for any $\mathbf{x}_0 \rmv\in\rmv \mathcal{X}^{\text{o}}\rmv$ and any $\mathbf{p} \in \mathbb{Z}_{+}^{N}$ with $p_k \rmv\leq\rmv m$,}
%% $\|\mathbf{p}\|_{\infty}\leq m$, 
the following 
\vspace{1mm}
holds:
%% we have that: 
\begin{itemize}

\item
The function $r^{(\mathbf{p})}_{\mathbf{x}_0}(\cdot) \!: \mathcal{X} \!\rightarrow\rmv \mathbb{R}$ defined by 
\begin{equation} 
\label{equ_def_part_der_func}
r^{(\mathbf{p})}_{\mathbf{x}_0} (\mathbf{x}) \,\triangleq\, \frac{\partial^{\mathbf{p}}  R(\mathbf{x}, \mathbf{x}_{2})}{\partial \mathbf{x}_{2}^{\mathbf{p}}} \bigg|_{\mathbf{x}_{2} = \mathbf{x}_0}
\vspace{-1mm}
\end{equation}
%where $\mathbf{x}_0 \in \mathcal{X}$ is such that $\mathcal{N}^{\supp(\mathbf{p})}_{\mathbf{x}_0}(\varepsilon) \subseteq \mathcal{X}$
%(for a suitable $\varepsilon >0$), 
is an element  of $\mathcal{H}(R)$, i.e., 
\vspace{1.5mm}
$r^{(\mathbf{p})}_{\mathbf{x}_0} (\cdot) \rmv\in \mathcal{H}(R)$.
%% \[ 
%% r^{(\mathbf{p})}_{\mathbf{x}_0} (\cdot) \in \mathcal{H}(R). 
%% \] 

\item
%% Moreover, for 
For any function $f(\cdot) \in \mathcal{H}(R)$, the partial derivative $\frac{\partial^{\mathbf{p}}  f(\mathbf{x})}{\partial \mathbf{x}^{\mathbf{p}}} \big|_{\mathbf{x} = \mathbf{x}_0}\!$ 
\vspace{1.5mm}
exists. 

\item
%% Finally, the 
The inner product of 
%% the function 
$r^{(\mathbf{p})}_{\mathbf{x}_0} (\cdot)$ with an arbitrary function $f(\cdot) \in \mathcal{H}(R)$ is given 
\vspace{.5mm}
by 
\begin{equation}
\label{equ_der_reproduction_prop} 
\big\langle r^{(\mathbf{p})}_{\mathbf{x}_0}(\cdot) , f(\cdot) \big\rangle_{\mathcal{H}(R)} 
  \ist=\, \frac{\partial^{\mathbf{p}}  f(\mathbf{x})}{\partial \mathbf{x}^{\mathbf{p}}} \bigg|_{\mathbf{x} = \mathbf{x}_0} \,.
\vspace{.5mm}
\end{equation} 

\end{itemize}

%\end{theorem} 

%% \begin{proof}
%% \cite{zhou_jfaa}
%% \end{proof}

Thus, an RKHS $\mathcal{H}(R)$ that is differentiable up to order $m$
%% \textcolor{red}{, for any $\mathbf{x}_0 \in  \mathcal{X}^{\text{o}}$,} 
%%  \in \mathbb{N}$ 
contains the functions $\big\{ r^{(\mathbf{p})}_{\mathbf{x}_0} (\mathbf{x}) \big\}_{p_k \leq m}$, 
%% _{\| \mathbf{p} \|_{\infty} \leq m}$. Moreover, 
and the inner products of any function $f(\cdot) \in \mathcal{H}(R)$ with the $r^{(\mathbf{p})}_{\mathbf{x}_0} (\mathbf{x})$
%% these functions with can be, according to \eqref{equ_der_reproduction_prop}, 
can be computed easily via differentiation of $f(\cdot)$.
%% that function. 
This makes function sets $\big\{ r^{(\mathbf{p})}_{\mathbf{x}_0} (\mathbf{x}) \big\}$ appear as interesting
%%  reasonable 
candidates for a simple characterization of the RKHS $\mathcal{H}(R)$. 
However, in general, these function sets are not guaranteed to be complete or orthonormal, i.e., they do not 
constitute an ONB. An important exception is constituted by certain estimation problems $\mathcal{E}$ involving an exponential family of distributions,
which will be studied in Section \ref{sec_rkhs_exp_fam}.
%% However, in what follows we will see that for certain estimation problems $\mathcal{E}$ involving exponential families, 
%% the associated kernel $R_{\mathcal{E}\rmv,\mathbf{x}_{0}}$ is differentiable to any order, and therefore contains the 
%% countable infinite set $\{ r^{(\mathbf{p})}_{\mathbf{x}_0} (\mathbf{x})\}_{\mathbf{p}  \in \mathbb{Z}_{+}^{N}} $. 
%% Moreover, for these estimation problems it will turn out that this function set is complete, i.e., the 
%% inner products of any function $f(\cdot) \in \mathcal{H}_{\mathcal{E}\rmv,\mathbf{x}_{0}}$ with the functions $r^{(\mathbf{p})}_{\mathbf{x}_0} (\cdot)$ 
%% completely determine the function $f(\cdot)$. 

Consider an estimation problem $\mathcal{E}=\scalarestproblem$ that is regular up to order $m \!\in\! \mathbb{N}$ at $\mathbf{x}_{0} \!\in\! \mathcal{X}^{\text{o}}\rmv$.
%% and the associated MVP at $\mathbf{x}_{0}$ for prescribed bias function $c(\cdot)$. 
\textcolor{red}{According to \eqref{equ_nec_suff_cond_validity_bias_function}}, the mean function $\gamma(\cdot)$ of any estimator 
with finite variance at $\mathbf{x}_{0}$ belongs to the RKHS $\mathcal{H}_{\mathcal{E}\rmv,\mathbf{x}_{0}}$. 
%By Theorem \ref{thm_regul_cond_CRB_imply_diff}, 
\textcolor{red}{Since $\mathcal{E}$ is assumed} regular up to order $m$, $\mathcal{H}_{\mathcal{E}\rmv,\mathbf{x}_{0}}$ is differentiable up to order $m$. 
This, in turn, 
implies\footnote{\textcolor{red}{Indeed, %%%%%%%%
it follows from \eqref{equ_nec_suff_cond_validity_bias_function} that the mean function $\gamma(\cdot)$ belongs to the 
RKHS $\mathcal{H}_{\mathcal{E}\rmv,\mathbf{x}_{0}}$. Therefore, by \eqref{equ_der_reproduction_prop}, the partial derivatives of 
$\gamma(\cdot)$ at $\mathbf{x}_{0}$ coincide with well-defined inner products of functions in 
$\mathcal{H}_{\mathcal{E}\rmv,\mathbf{x}_{0}}$.}} %%%%%%%%%
via \textcolor{red}{\eqref{equ_nec_suff_cond_validity_bias_function} and \eqref{equ_der_reproduction_prop}} that the partial derivatives of $\gamma(\cdot)$ at $\mathbf{x}_{0}$ exist up to order $m$. Therefore, for the derivation of lower bounds on the minimum achievable variance at $\mathbf{x}_{0}$ 
in the case of an estimation problem 
%% $\mathcal{E}=\scalarestproblem$ 
that is regular up to order $m$ at $\mathbf{x}_{0}$, we can always tacitly 
assume that the partial derivatives of $\gamma(\cdot)$ at $\mathbf{x}_{0}$ exist up to order $m$; otherwise the corresponding bias function
$c(\cdot)=\gamma(\cdot)-g(\cdot)$ cannot be valid, i.e., there would not exist any estimator with mean function $\gamma(\cdot)$ 
(or, equivalently, bias function $c(\cdot)$) and finite variance at 
\pagebreak %%%%%%%%
$\mathbf{x}_{0}$.

%% \vspace{-2mm}

%%%%%%%%%%%%%%%%%%%%%%%%%%%%%%%%%%%%%%%%%%%%%%%%%%%%%%%%%%
\section{RKHS Formulation
%% Interpretation 
of Known Variance Bounds} 
\label{sec_rkhs_interpret_well_known_bounds}
%%%%%%%%%%%%%%%%%%%%%%%%%%%%%%%%%%%%%%%%%%%%%%%%%%%%%%%%%%

\vspace{1mm}

Consider an estimation problem $\mathcal{E} = \scalarestproblem$ and an estimator $\hat{g}(\cdot)$ with 
%% finite variance at $\mathbf{x}_{0}$, 
mean function $\gamma(\mathbf{x}) = \expect_{\mathbf{x}} \{ \hat{g}(\mathbf{y}) \}$ and bias function 
$c(\mathbf{x}) = \gamma(\mathbf{x}) - g(\mathbf{x})$. We assume that $\hat{g}(\cdot)$ has a finite variance at $\mathbf{x}_{0}$,
which implies that the bias function $c(\cdot)$ is valid and $\hat{g}(\cdot)$ is an element of $\mathcal{A}({c}(\cdot),\mathbf{x}_{0})$, 
the set of allowed estimators at $\mathbf{x}_{0}$ for prescribed bias function ${c}(\cdot)$, which therefore is nonempty. 
Then, $\gamma(\cdot) \in  \mathcal{H}_{\mathcal{E}\rmv,\mathbf{x}_{0}}$ \textcolor{red}{according to \eqref{equ_nec_suff_cond_validity_bias_function}.}
We also recall from our discussion further above that if the estimation problem $\mathcal{E}$ is regular at $\mathbf{x}_{0}$ up to order $m$, 
%% it follows by \eqref{equ_crb_regul_cond_finite_part_der} and \eqref{equ_crb_regul_cond_intercah_diff_integr_intral_repr} that 
then the partial derivatives $\frac{\partial^{\mathbf{p}} \ist \gamma(\mathbf{x})}{\partial \mathbf{x}^{\mathbf{p}}} \big|_{\mathbf{x} = \mathbf{x}_{0}}\!$ exist
for all $\mathbf{p} \rmv\in\rmv \mathbb{Z}_{+}^{N}$ with $p_k \rmv\leq\rmv m$.
%% Thus, as already discussed above, we will tacitly assume that the partial derivatives up to a sufficient order exists for any mean function considered.
%The existence of the derivatives of $\gamma_{0}(\cdot)$ will be required by some of the bounds that we will now consider. 

%% \newpage %%%%%%%

In this section, we will demonstrate how five known lower bounds on the variance---Barankin bound, Cram\'{e}r--Rao bound, constrained Cram\'{e}r--Rao bound, 
Bhattacharyya bound, and Hammersley-Chapman-Robbins bound---can be formulated in a unified manner within the RKHS framework. More specifically,
%% these bounds can be obtained by applying the projection theorem to the RKHS $\mathcal{H}_{\mathcal{E}\rmv,\mathbf{x}_{0}}$.
%% %% associated with the estimation problem. 
by combining \eqref{equ_lower_vound_variance_trivial_min_achiev_var} with \eqref{equ_lower_bound_min_achiev_var_projection}, it follows that the variance of 
$\hat{g}(\cdot)$ at $\mathbf{x}_{0}$ is lower bounded as 
\begin{equation}
\label{equ_lower_bound_variance_projection} 
v(\hat{g}(\cdot);\mathbf{x}_{0})  \,\geq\, {\| \gamma_{\mathcal{U}}(\cdot) \|}^{2}_{\mathcal{H}_{\mathcal{E}\rmv,\mathbf{x}_{0}}} \!- \gamma^{2}(\mathbf{x}_{0}) \,, 
\end{equation}
where $\mathcal{U}$ is any subspace of $\mathcal{H}_{\mathcal{E}\rmv,\mathbf{x}_{0}}$.
The five variance bounds to be considered are obtained via specific choices of $\mathcal{U}$.

\vspace{-1mm}

\subsection{Barankin Bound} 
\label{SecBarankin}
%%%%%%%%%%%%%%%%%%%%%%%%%%%%%%%%%%%%%%%%%%%%%%%%%%%%%%%%%%

For a (valid) prescribed bias function $c(\cdot)$, the Barankin bound \cite{Barankin49,mcaulay71} is the minimum achievable variance at $\mathbf{x}_{0}$, 
i.e., the variance of the LMV estimator at $\mathbf{x}_{0}$,
which we denoted $M(c(\cdot),\mathbf{x}_{0})$. This is the tightest lower bound on the variance, cf.\ \eqref{equ_lower_vound_variance_trivial_min_achiev_var}. 
Using the RKHS expression of $M(c(\cdot),\mathbf{x}_{0})$ in \eqref{equ_min_achiev_var_sqared_norm}, the Barankin bound can be written 
\vspace{-2mm}
as
\begin{equation} 
\label{equ_lower_vound_variance_trivial_min_achiev_var_barankin}
v(\hat{g}(\cdot);\mathbf{x}_{0}) \,\geq\, M(c(\cdot),\mathbf{x}_{0}) 
\eq {\| \gamma(\cdot) \|}^{2}_{\mathcal{H}_{\mathcal{E}\rmv,\mathbf{x}_{0}}} \!- \gamma^{2}(\mathbf{x}_{0})\,,
\end{equation}
with $\gamma(\cdot) = c(\cdot) + g(\cdot)$, for any estimator $\hat{g}(\cdot)$ with bias function $c(\cdot)$.
Comparing with \eqref{equ_lower_bound_variance_projection}, we see that the Barankin bound is obtained for the special choice 
$\mathcal{U} =  \mathcal{H}_{\mathcal{E}\rmv,\mathbf{x}_{0}}$, in which case $\gamma_{\mathcal{U}}(\cdot) = \gamma(\cdot)$ 
%%reduces to the identity operator
and \eqref{equ_lower_bound_variance_projection} reduces to \eqref{equ_lower_vound_variance_trivial_min_achiev_var_barankin}. 
%% It can be shown \cite[Section 2.3.5]{JungPHD} that

In the literature \cite{Barankin49,mcaulay71}, the following special expression of the Barankin bound is usually considered.
Let $\mathcal{D} \triangleq \{\mathbf{x}_{1},\ldots,\mathbf{x}_{L}\} \subseteq \mathcal{X}$
be a subset of $\mathcal{X}$, with finite size $L = |\mathcal{D}| \in \mathbb{N}$ and elements $\mathbf{x}_{l}  \in \mathcal{X}$, 
and let $\mathbf{a} \triangleq ( a_{1} \cdots\ist a_{L})^T\!$ with $a_l \in \mathbb{R}$.
Then the Barankin bound 
%% $v(\hat{g}(\cdot);\mathbf{x}_{0}) \geq M(c(\cdot),\mathbf{x}_{0})$ 
can be written as \cite[Theorem 4]{Barankin49}
\begin{equation}
\label{equ_barankin_bound1}
v(\hat{g}(\cdot);  \mathbf{x}_{0}) 
\,\geq\, M(c(\cdot),\mathbf{x}_{0}) 
\,= \sup_{\mathcal{D} \subseteq \mathcal{X}, 
\ist\ist L\in \mathbb{N}, 
\ist\ist \mathbf{a} \in \mathcal{A}_{\mathcal{D}}}
%% \sup_{\substack{\mathcal{D} \in \mathcal{X}^L \\ L\in \mathbb{N}\mbox{,} \,\, \mathbf{x}_{l} \in \mathcal{X} \\ \mathbf{a} \in \mathcal{B}_{\mathcal{D}}}} 
\,\frac{ \Big( \rmv \sum_{l \in [L]}  a_{l} \ist\ist [ \gamma(\mathbf{x}_l) \rmv-\rmv \gamma(\mathbf{x}_{0}) ] 
%% h(\mathbf{x}_{l}) 
\Big)^{\!2}}{ \expect_{\mathbf{x}_{0}} \rmv\Big\{ \Big( \rmv \sum_{l \in [L]} a_{l} \ist\ist \llr(\mathbf{y},\mathbf{x}_{l}) \Big)^{\! 2} \Big\}} \,,
\end{equation}
where $\llr(\mathbf{y},\mathbf{x}_{l})$ is the likelihood ratio as defined in \eqref{equ_def_likelihood}
and $\mathcal{A}_{\mathcal{D}}$ is defined as the set of all $\mathbf{a} \rmv\in\rmv \mathbb{R}^{L}$ for which the denominator 
$\expect_{\mathbf{x}_{0}} \rmv\big\{ \big( \sum_{l \in [L]} a_{l} \ist\ist \llr(\mathbf{y},\mathbf{x}_{l}) \big)^{\rmv 2} \big\}$ does not vanish.
%% $\mathcal{B}_{\mathcal{D}} \triangleq \Big\{ \mathbf{a} \rmv\in\rmv \mathbb{R}^{L} \ist \Big|\, 
%% \expect_{\mathbf{x}_{0}} \rmv\Big\{ \Big( \sum_{l \in [L]} a_{l} \ist\ist \llr(\mathbf{y},\mathbf{x}_{l}) \Big)^{\!2} \Big\} \neq 0 \Big\}$
Note that our notation $\sup_{\mathcal{D} \subseteq \mathcal{X}, \ist\ist L\in \mathbb{N}, \ist\ist \mathbf{a} \in \mathcal{A}_{\mathcal{D}}}$ 
is intended to indicate that the supremum is taken not only with respect to the elements $\mathbf{x}_{l} $ of $\mathcal{D}$ but also with respect to
the size of $\mathcal{D}$ (number of elements), $L$.
We will now verify that the bound in \eqref{equ_barankin_bound1} can be obtained from our RKHS expression in \eqref{equ_lower_vound_variance_trivial_min_achiev_var_barankin}.
We will use the following result that we reported in \cite[Theorem 3.1.2]{JungPHD}.

\vspace{1mm}

\begin{lemma}
\label{thm_approx_norm_dense_set}
Consider an RKHS $\mathcal{H}(R)$ with kernel $R(\cdot\ist\ist,\cdot) \!: \mathcal{X} \rmv\times\rmv \mathcal{X} \!\rightarrow \mathbb{R}$. 
Let $\mathcal{D} \triangleq \{\mathbf{x}_{1},\ldots,\mathbf{x}_{L}\} \subseteq \mathcal{X}$
%%  \in \mathcal{X}^L$ 
with some $L = |\mathcal{D}| \in \mathbb{N}$ and $\mathbf{x}_{l} \!\in\! \mathcal{X}$, and let $\mathbf{a} \triangleq ( a_{1} \cdots\ist a_{L})^T\!$ with $a_l \!\in\! \mathbb{R}$.
Then the norm ${\| f(\cdot) \|}_{\mathcal{H}(R)}$ of any function $f(\cdot) \in \mathcal{H}(R)$ can be expressed as
\begin{equation}
\label{equ_approx_norm_HS_sup_dense_set}
{\| f(\cdot) \|}_{\mathcal{H}(R)} \,= \sup_{\mathcal{D} \subseteq \mathcal{X}, \ist\ist L\in \mathbb{N}, \ist\ist \mathbf{a} \in \mathcal{A}'_{\mathcal{D}}}
\frac{\sum_{l \in [L]} a_{l} \ist f(\mathbf{x}_{l})}{\sqrt{\sum_{l,l' \in [L]} a_{l} \ist a_{l'} R(\mathbf{x}_{l}, \mathbf{x}_{l'})}} \,,
\end{equation} 
where $\mathcal{A}'_{\mathcal{D}}$
%%  \triangleq \big\{ \sum_{k,k' \in [L]} a_{k} a_{k'} R(\mathbf{x}_{k}, \mathbf{x}_{k'})> 0 \big\}$.
is the set of all $\mathbf{a} \rmv\in\rmv \mathbb{R}^{L}$ for which 
%% the denominator 
$\sum_{l,l' \in [L]} a_{l} \ist a_{l'} R(\mathbf{x}_{l}, \mathbf{x}_{l'})$ does not vanish.
\end{lemma}

\vspace{1mm}

%% \begin{proof}
%% \cite{JungPHD}
%% \end{proof} 

%% \noindent
We will furthermore use the fact---shown in \cite[Section 2.3.5]{JungPHD}---that the minimum achievable variance at $\mathbf{x}_{0}$, $M(c(\cdot),\mathbf{x}_{0})$ 
(i.e., the Barankin bound)
remains unchanged when the prescribed mean function $\gamma(\mathbf{x})$ is replaced by $\tilde\gamma(\mathbf{x}) \triangleq \gamma(\mathbf{x}) + c$
with an arbitrary constant $c$. Setting in particular $c = -\gamma(\mathbf{x}_0)$, we have $\tilde\gamma(\mathbf{x}) = \gamma(\mathbf{x}) - \gamma(\mathbf{x}_0)$
and $\tilde\gamma(\mathbf{x}_0) = 0$,
and thus \eqref{equ_lower_vound_variance_trivial_min_achiev_var_barankin} simplifies to
\begin{equation} 
\label{equ_lower_vound_variance_trivial_min_achiev_var_barankin_0}
v(\hat{g}(\cdot);\mathbf{x}_{0}) \,\geq\, M(c(\cdot),\mathbf{x}_{0}) 
\eq {\| \tilde\gamma(\cdot) \|}^{2}_{\mathcal{H}_{\mathcal{E}\rmv,\mathbf{x}_{0}}} \,.
\vspace{-1.5mm}
\end{equation}
Using \eqref{equ_approx_norm_HS_sup_dense_set} in \eqref{equ_lower_vound_variance_trivial_min_achiev_var_barankin_0}, we obtain
\be
M(c(\cdot),\mathbf{x}_{0}) \,= \sup_{\mathcal{D} \subseteq \mathcal{X}, \ist\ist L\in \mathbb{N}, \ist\ist \mathbf{a} \in \mathcal{A}'_{\mathcal{D}}}
\frac{ \Big( \rmv \sum_{l \in [L]} a_{l} \ist\ist \tilde\gamma(\mathbf{x}_{l}) \Big)^{\!2} }{ \sum_{l,l' \in [L]} a_{l} \ist a_{l'} R_{\mathcal{E}\rmv,\mathbf{x}_{0}}(\mathbf{x}_{l}, \mathbf{x}_{l'})} 
\,= \sup_{\mathcal{D} \subseteq \mathcal{X}, \ist\ist L\in \mathbb{N}, \ist\ist \mathbf{a} \in \mathcal{A}'_{\mathcal{D}}}
\!\rmv\frac{ \Big( \rmv \sum_{l \in [L]} a_{l} \ist\ist [ \gamma(\mathbf{x}_l) \rmv-\rmv \gamma(\mathbf{x}_{0}) ]  \Big)^{\!2} }{ \sum_{l,l' \in [L]} a_{l} \ist a_{l'} 
  R_{\mathcal{E}\rmv,\mathbf{x}_{0}}(\mathbf{x}_{l}, \mathbf{x}_{l'})} \,.
\label{equ_barankin_bound_der2}
\vspace{-1mm}
\ee
From \eqref{equ_def_kernel_est_problem} and \eqref{equ_def_inn_prod_RV}, we have 
$R_{\mathcal{E}\rmv,\mathbf{x}_{0}}(\mathbf{x}_{1},\mathbf{x}_{2}) = \expect_{\mathbf{x}_{0}} \big\{  \llr(\mathbf{y},\mathbf{x}_{1}) \ist\ist \llr (\mathbf{y}, \mathbf{x}_{2}) \big\}$, 
and thus the denominator in \eqref{equ_barankin_bound_der2} becomes
\vspace{-1.5mm}
\[
\sum_{l,l' \in [L]} \!a_{l} \ist a_{l'} R_{\mathcal{E}\rmv,\mathbf{x}_{0}}(\mathbf{x}_{l}, \mathbf{x}_{l'}) 
\eq \expect_{\mathbf{x}_{0}} \rmv\Bigg\{ \sum_{l,l' \in [L]} \!a_{l} \ist a_{l'} \llr(\mathbf{y},\mathbf{x}_{l}) \ist\ist \llr (\mathbf{y}, \mathbf{x}_{l'}) \Bigg\}
\eq \expect_{\mathbf{x}_{0}} \rmv\Bigg\{ \! \Bigg( \sum_{l \in [L]} \rmv a_{l} \ist\ist \llr(\mathbf{y},\mathbf{x}_{l}) \rmv\rmv \Bigg)^{\!\! 2} \Bigg\} \,,
\vspace{.5mm}
\]
whence it also follows that $\mathcal{A}'_{\mathcal{D}} = \mathcal{A}_{\mathcal{D}}$. Therefore, \eqref{equ_barankin_bound_der2} is equivalent to \eqref{equ_barankin_bound1}.
Hence, we have shown that our RKHS expression \eqref{equ_lower_vound_variance_trivial_min_achiev_var_barankin} is equivalent to \eqref{equ_barankin_bound1}.

\vspace{-1mm}

\subsection{Cram\'{e}r--Rao Bound}
\label{sec_CRB}
%%%%%%%%%%%%%%%%%%%%%%%%%%%%%%%%%%%%%%%%%%%%%%%%%%%%%%%%%%

The \CRBfull (CRB) \cite{cramer45,rao45,kay} is the most popular
%% best known 
lower variance bound.
%%  on the variance of estimators for a given estimation problem $\mathcal{E}=\scalarestproblem$ which satisfies some weak regularity conditions as discussed presently. 
Since the CRB applies to any estimator with a 
%% given 
prescribed bias function $c(\cdot)$, it yields also a lower bound on the minimum achievable variance $M(c(\cdot),\mathbf{x}_{0})$ 
(cf.\ \eqref{equ_lower_vound_variance_trivial_min_achiev_var}).

%% We will consider the CRB for two different settings, which are defined by the structure of the parameter set $\mathcal{X}$ 
%% associated with the estimation problem $\mathcal{E}=\scalarestproblem$. These two instances of the CRB are termed the \emph{unconstrained CRB} and the 
%% \emph{constrained CRB}, respectively. 
%% %and the moments of $\Phi(\mathbf{y})$ (cf. \eqref{equ_def_pdf_exp_fam}) exist up to a sufficient order which depends on $m$. 

%% \newpage %%%%%%%%

Consider an estimation problem $\mathcal{E} = \scalarestproblem$ that is regular \textcolor{red}{
up to order 1 at $\mathbf{x}_{0} \rmv\in\rmv \mathcal{X}^{\text{o}}$ in the sense of Definition \ref{post_regular_cond_CRB}.
Let $\hat{g}(\cdot)$ denote an estimator} with mean function $\gamma(\mathbf{x}) = \expect_{\mathbf{x}} \{ \hat{g}(\mathbf{y}) \}$ 
and finite variance at $\mathbf{x}_{0}$ (i.e., $v(\hat{g}(\cdot);\mathbf{x}_{0}) < \infty$). Then, this variance is lower bounded by the 
%% unconstrained 
CRB 
%% (cf.\ \cite{HeroRecursiveCRB,HeroUniformCRB})
\begin{equation}
\label{equ_unconstrained_CRB}
v(\hat{g}(\cdot);\mathbf{x}_{0}) \,\geq\, \mathbf{b}^T\rmv(\mathbf{x}_{0}) \,
%%  \bigg( \frac{\partial \gamma(\mathbf{x})} {\partial \mathbf{x}}\bigg|_{\mathbf{x}_{0}} \bigg)^{\!\rmv T} 
  \mathbf{J}^{\dagger}(\mathbf{x}_{0}) \ist\ist \mathbf{b}(\mathbf{x}_{0}) \,,
  %% \frac{\partial \gamma(\mathbf{x})} {\partial \mathbf{x}}\bigg|_{\mathbf{x}_{0}} ,
\end{equation} 
where $\mathbf{b}(\mathbf{x}_{0}) \triangleq \frac{\partial \gamma(\mathbf{x})} {\partial \mathbf{x}}\big|_{\mathbf{x}_{0}}$ and 
$\mathbf{J}(\mathbf{x}_{0}) \rmv\in\rmv \mathbb{R}^{N \times N}\rmv$, known as the \emph{Fisher information matrix} 
%% (FIM) 
associated with 
%% the estimation problem 
$\mathcal{E}$, is given elementwise by 
\begin{equation} 
\label{equ_def_FIM_CRB}
\big( \mathbf{J}(\mathbf{x}_{0}) \big)_{k,l} \ist\triangleq\, \expect_{\mathbf{x}_{0}} \bigg\{  \frac{\partial \log f(\mathbf{y}; \mathbf{x})}{\partial x_{k}} 
%% \bigg|_{\mathbf{x} = \mathbf{x}_{0}}
\, \frac{\partial \log f(\mathbf{y}; \mathbf{x})}{\partial x_{l}} \bigg|_{\mathbf{x} = \mathbf{x}_{0}} \bigg \} \,.
\vspace{.5mm}
\end{equation}

\textcolor{red}{Since the estimation problem $\mathcal{E}$ is assumed regular up to order 1 at $\mathbf{x}_{0}$, 
the associated RKHS $\mathcal{H}_{\mathcal{E}\rmv,\mathbf{x}_{0}}$ is differentiable up to order 1.} %, according to
%Theorem \ref{thm_regul_cond_CRB_imply_diff}. 
This differentiability is used in the proof of the following result \cite[Section 4.4.2]{JungPHD}.
%% our subsequent development.

%% \subsubsection{Unconstrained CRB}

\vspace{1mm}

\begin{theorem} 
\label{th_crb}
Consider an estimation problem that is regular up to order 1 in the sense of Definition \ref{post_regular_cond_CRB}. Then, for a reference parameter vector 
\emph{$\mathbf{x}_{0} \rmv\in\rmv \mathcal{X}^{\text{o}}$}, the 
%% unconstrained 
CRB in \eqref{equ_unconstrained_CRB} is obtained from \eqref{equ_lower_bound_variance_projection} by using the subspace 
\[
\mathcal{U}_{\text{\emph{CR}}} \ist\triangleq\, \linspan \rmv\big\{ \{ \basisfuncRKHS_{0} (\cdot) \} \cup {\{ \basisfuncRKHS_{l}(\cdot) \}}_{l \in [N]} \ist \big\} \,, 
%\{ f(\cdot) \in \mathcal{H}(\mathcal{M}) \big| f(\cdot) = a_{0} v_{0}(\cdot) + \sum_{l=1}^{N} a_{l} v_{l}(\cdot) \,\, \mbox{,} \,\, a_{0}, a_{l} \in \mathbb{R} \} 
\vspace{-3mm}
\] 
with the 
\vspace{-.5mm}
functions 
%% $v_{0} (\cdot ) \triangleq R_{\mathcal{E}\rmv,\mathbf{x}_{0}}(\cdot, \mathbf{x}_{0}) \in \mathcal{H}_{\mathcal{E}\rmv,\mathbf{x}_{0}}$ and 
\[
%% \label{equ_subspace_CR_RKHS_vectors}
\basisfuncRKHS_{0} (\cdot ) \ist\triangleq\ist R_{\mathcal{E}\rmv,\mathbf{x}_{0}}(\cdot\ist\ist, \mathbf{x}_{0}) \in \mathcal{H}_{\mathcal{E}\rmv,\mathbf{x}_{0}} \,, \qquad\;
\basisfuncRKHS_{l}(\cdot) \ist\triangleq\ist \frac{\partial R_{\mathcal{E}\rmv,\mathbf{x}_{0}}(\cdot\ist\ist,\mathbf{x})}{\partial x_{l}} \bigg|_{\mathbf{x} = \mathbf{x}_{0}} 
  \!\!\in \mathcal{H}_{\mathcal{E}\rmv,\mathbf{x}_{0}} \,, \quad l \rmv\in\rmv [N] \,.  
\vspace{1mm}
\] 
\end{theorem}

%% \begin{proof}
%% \cite[Section 4.4.2]{JungPHD}
%% \end{proof}

\subsection{Constrained Cram\'{e}r--Rao Bound} 
\label{sec_ConstrCRB}
%%%%%%%%%%%%%%%%%%%%%%%%%%%%%%%%%%%%%%%%%%%%%%%%%%%%%%%%%%

%% Consider an estimation problem $\mathcal{E}=\scalarestproblem$, with the parameter set % $\mathcal{X}$ which is defined via a set of equality constraints, i.e., 
The constrained CRB \cite{StoicaNgCCRB,ZvikaCCRB,MooreCCRB} is an evolution
%%  generalization 
of the CRB in \eqref{equ_unconstrained_CRB} for estimation problems $\mathcal{E}=\scalarestproblem$ with a parameter set 
%% $\mathcal{X}$ 
of the 
\vspace{-1mm}
form 
%% \eqref{equ_def_paramter_set_CCRB}.
\begin{equation} 
\label{equ_def_parameter_set_CCRB_12}
%% \label{equ_def_paramter_set_CCRB} 
\mathcal{X} = \big\{ÊÊ\mathbf{x} \rmv\in\rmv \mathbb{R}^{N} \big|\ist \mathbf{f}(\mathbf{x}) \rmv=\rmv \mathbf{0} \big\} \,,
\end{equation} 
where $\mathbf{f}(\cdot) \!: \mathbb{R}^{N} \!\rightarrow \mathbb{R}^{Q}$ with $Q \!\le\! N$ is a continuously differentiable
%% \footnote{cf.\ \cite[Definition 9.20]{RudinBookPrinciplesMatheAnalysis}} 
%% vector-valued 
function. We assume that the set $\mathcal{X}$ has a nonempty interior.
Moreover, we require the Jacobian matrix $\mathbf{F}(\mathbf{x}) \triangleq \frac{\partial \ist \mathbf{f}(\mathbf{x})}{ \partial \mathbf{x}} \in \mathbb{R}^{Q \times N}$ 
%% given elementwise by $\left( \mathbf{F}(\mathbf{x}) \right)_{m,n} \triangleq \frac{\partial  f_{m}(\mathbf{x})}{\partial x_{n}}$, 
to have rank $Q$ whenever $\mathbf{f}(\mathbf{x})=\mathbf{0}$, i.e., for every $\mathbf{x} \in \mathcal{X}$. 
This full-rank requirement implies that the constraints represented by $\mathbf{f}(\mathbf{x}) = \mathbf{0}$ are nonredundant \cite{ZvikaCCRB}. 
Such parameter sets are considered, e.g., in \cite{StoicaNgCCRB,ZvikaCCRB,MooreCCRB}. Under these conditions, 
%% we have according to 
the implicit function theorem 
%% (cf., e.g., 
\cite[Theorem 3.3]{MooreCCRB}, \cite[Theorem 9.28]{RudinBookPrinciplesMatheAnalysis} states
that for any $\mathbf{x}_{0} \rmv\in\rmv \mathcal{X}$, with $\mathcal{X}$ given by \eqref{equ_def_parameter_set_CCRB_12}, there exists a continuously differentiable map $\mathbf{r}(\cdot)$ from an open set $\mathcal{O} \subseteq \mathbb{R}^{N-Q}$ into a set $\mathcal{P} \subseteq \mathcal{X}$ 
%% also 
containing $\mathbf{x}_{0}$, 
\vspace{-2mm}
i.e., 
\begin{equation} 
\label{equ_ccrb_reparam} 
\mathbf{r}(\cdot) \rmv :\ist  \mathcal{O} \subseteq \mathbb{R}^{N-Q} \ist\rightarrow\ist\ist \mathcal{P} \rmv\subseteq\rmv \mathcal{X} \ist, \quad 
%% \mathcal{O}, \mathcal{P} \mbox{ open, and } 
\text{with} \;\ist \mathbf{x}_{0} \!\in\rmv \mathcal{P}.
\end{equation} 
%and moreover $\mathbf{r}(\cdot)$ has a continuously differentiable inverse
%% \footnote{Such a function is called a diffeomorphism \cite{MooreCCRB}.} 
%$\mathbf{r}^{-1}(\cdot)$.
%% where $\mathbf{x}_{0} \in \mathcal{P}$. 
%% Note that any function value $\mathbf{r}({\bm \theta})$ is an element of the parameter set, i.e., $\mathbf{r}({\bm \theta}) \in \mathcal{X}$. 

%% \newpage %%%%%%%%

The constrained CRB in the form presented in \cite{ZvikaCCRB} reads
\begin{equation} 
\label{equ_CCRB} 
v(\hat{g}(\cdot); \mathbf{x}_{0}) \,\geq\, \mathbf{b}^T\rmv(\mathbf{x}_{0}) \ist
\mathbf{U}(\mathbf{x}_{0}) \,\big( \mathbf{U}^{T}\rmv(\mathbf{x}_{0}) \,\mathbf{J}(\mathbf{x}_{0}) \ist\mathbf{U}(\mathbf{x}_{0}) \big)^{\rmv\dagger} \,
\mathbf{U}^{T}\rmv(\mathbf{x}_{0}) \ist\ist \mathbf{b}(\mathbf{x}_{0}) \,, 
\end{equation} 
where $\mathbf{b}(\mathbf{x}_{0}) = \frac{\partial \gamma(\mathbf{x})} {\partial \mathbf{x}}\big|_{\mathbf{x}_{0}}$,
$\mathbf{J}(\mathbf{x}_{0})$ is again the Fisher information matrix defined in \eqref{equ_def_FIM_CRB}, 
and $\mathbf{U}(\mathbf{x}_{0}) \in \mathbb{R}^{N \times (N-Q)}$ is any matrix whose columns
%%  vectors 
form an ONB for the null space 
%% $\Ker(\mathbf{F}(\mathbf{x}_{0}))$ 
of the Jacobian matrix $\mathbf{F}(\mathbf{x}_{0})$, i.e., 
\[
\mathbf{F}(\mathbf{x}_{0}) \ist\mathbf{U}(\mathbf{x}_{0}) = \mathbf{0}\,, \quad\;\; \mathbf{U}^{T}\rmv(\mathbf{x}_{0}) \ist\mathbf{U}(\mathbf{x}_{0}) = \mathbf{I}_{N-Q} \,.  
\] 
%Note that if 
%% $\mathbf{f}(\mathbf{x})$ in \eqref{equ_def_paramter_set_CCRB} is given as 
%$\mathbf{f}(\mathbf{x}) \equiv \mathbf{0}$ (unconstrained case), 
%% one can verify that 
%% we have $Q=???$ and 
%the matrix $\mathbf{U}(\mathbf{x}_{0})$ 
%% used in \eqref{equ_CCRB} 
%can be chosen as $\mathbf{U}(\mathbf{x}_{0}) = \mathbf{I}_{???}$, in which case 
%the constrained CRB \eqref{equ_CCRB} reduces to the unconstrained CRB \eqref{equ_unconstrained_CRB}. 
The next result is proved in \cite[Section 4.4.2]{JungPHD}.

\vspace{1mm}

\begin{theorem} 
Consider an estimation problem that is regular up to order 1 in the sense of Definition \ref{post_regular_cond_CRB}. Then, for a reference parameter vector 
\emph{$\mathbf{x}_{0} \rmv\in\rmv \mathcal{X}^{\text{o}}$}, 
the constrained CRB in \eqref{equ_CCRB} is obtained from \eqref{equ_lower_bound_variance_projection} by using the 
%% \vspace{-2mm}
\pagebreak %%%%%%
subspace 
\[
%% \label{equ_def_U_CCRB}
\mathcal{U}_{\text{\emph{CCR}}} \ist\triangleq\, \linspan \rmv\big\{ \{ \basisfuncRKHS_{0} (\cdot) \} \cup {\{ \basisfuncRKHS_{l}(\cdot) \}}_{l \in [N-Q]} \big\} \,,
\vspace{-3mm}
\] 
with the functions 
%% $v_{0} (\cdot ) \triangleq R_{\mathcal{E}\rmv,\mathbf{x}_{0}}(\cdot, \mathbf{x}_{0}) \in \mathcal{H}(\mathcal{M})$ and 
\[
%% \label{equ_def_vectors_ccrb}
\basisfuncRKHS_{0} (\cdot ) \ist\triangleq\ist R_{\mathcal{E}\rmv,\mathbf{x}_{0}}(\cdot\ist\ist, \mathbf{x}_{0}) \in \mathcal{H}_{\mathcal{E}\rmv,\mathbf{x}_{0}} \,, \qquad\;
\basisfuncRKHS_{l}(\cdot) \ist\triangleq\ist \frac{\partial R_{\mathcal{E}\rmv,\mathbf{x}_{0}}(\cdot\ist\ist,\mathbf{r}({\bm \theta}))}{\partial \theta_{l}} \bigg|_{{\bm \theta} = \mathbf{r}^{-1}(\mathbf{x}_{0})} 
  \!\! \in \mathcal{H}_{\mathcal{E}\rmv,\mathbf{x}_{0}} \,, \quad l \rmv\in\rmv [N \!-\! Q] \,,  
\]
where $\mathbf{r}(\cdot)$ is any continuously differentiable function of the form \eqref{equ_ccrb_reparam}. 
%% The existence of such a function is guaranteed by the implicit function theorem \cite[Theorem 9.28]{RudinBookPrinciplesMatheAnalysis}. 
\vspace{-1mm}
\end{theorem}

%% \begin{proof}
%% \cite[Section 4.4.2]{JungPHD}
%% \end{proof}

\subsection{Bhattacharyya Bound} 
\label{sec_Bhatt}
%%%%%%%%%%%%%%%%%%%%%%%%%%%%%%%%%%%%%%%%%%%%%%%%%%%%%%%%%%

Whereas the CRB depends only on the first-order partial derivatives of 
%% the statistical model 
$f(\mathbf{y}; \mathbf{x})$ with respect to 
%% the parameter 
$\mathbf{x}$, the Bhattacharyya bound \cite{LowerBoundAbel,bhattacharyya} involves also higher-order 
%% partial 
derivatives. 
For an estimation problem $\mathcal{E}=\scalarestproblem$ that is regular at 
%% a specific 
$\mathbf{x}_{0}Ê\rmv\in\rmv \mathcal{X}^{\text{o}}\rmv$ up to order $m \in \mathbb{N}$, 
%% where in general $m> 1$. 
the Bhattacharyya bound states that 
%Like the CRB the Bhattacharyya bound, is a lower bound on the variance 
%of any unbiased estimator for a given estimation problem (which is required to satisfy some regularity conditions). The Bhattacharyya bound may be formulated in our setting as follows.
\begin{equation}
\label{equ_bhattacharyya_bound}
v(\hat{g}(\cdot);\mathbf{x}_{0}) \,\geq\, \mathbf{a}^T\rmv(\mathbf{x}_{0}) \, \mathbf{B}^{\dagger}(\mathbf{x}_{0}) \ist\ist \mathbf{a}(\mathbf{x}_{0}) \,,
\end{equation} 
where the vector $\mathbf{a}(\mathbf{x}_{0}) \rmv\in\rmv \mathbb{R}^{L}\rmv$ and the matrix $\mathbf{B}(\mathbf{x}_{0}) \rmv\in\rmv \mathbb{R}^{L \times L}\rmv$ 
are given elementwise by 
$\big( \mathbf{a}(\mathbf{x}_{0}) \big)_l \triangleq \frac{\partial^{{\mathbf{p}}_{l}} \gamma(\mathbf{x})} {\partial \mathbf{x}^{{\mathbf{p}}_{l}} }\big|_{\mathbf{x}_{0}}$
%% \begin{equation} 
%% b_{l} \triangleq \frac{\partial^{\tilde{\mathbf{p}}_{l}} \gamma_{0}(\mathbf{x})} {\partial \mathbf{x}^{\tilde{\mathbf{p}}_{l}} }\bigg|_{\mathbf{x}_{0}}
%% \end{equation}
and 
\vspace{-2mm}
\[ 
%% \label{equ_def_bhattacharyya_matrix}
\big( \mathbf{B}(\mathbf{x}_{0}) \big)_{l,l'} \ist\triangleq\, \expect_{\mathbf{x}_{0}}
 \bigg\{  \frac{1}{f^{2}(\mathbf{y}; \mathbf{x}_{0})} \ist \frac{\partial^{{\mathbf{p}}_{l}} \rmv f(\mathbf{y}; \mathbf{x})}{\partial \mathbf{x}^{{\mathbf{p}}_{l}} } \ist
 \frac{\partial^{{\mathbf{p}}_{l'}} \rmv f(\mathbf{y}; \mathbf{x})}{\partial \mathbf{x}^{{\mathbf{p}}_{l'}} }  \bigg|_{\mathbf{x} = \mathbf{x}_{0}} \bigg \} \,,
\vspace{.5mm}
\] 
respectively. Here, the ${\mathbf{p}}_{l}$, $l \rmv\in\rmv [L]$ are $L$ distinct multi-indices 
%% $\{ \tilde{\mathbf{p}}_{l} \}_{l \in [L]}$ 
with ${({\mathbf{p}}_{l})}_k \leq m$.

The following result is proved in \cite[Section 4.4.3]{JungPHD}.

\vspace{1mm}

\begin{theorem} 
\label{th_bhat}
Consider an estimation problem that is regular up to order $m$ in the sense of Definition \ref{post_regular_cond_CRB}. Then, for a reference parameter vector 
\emph{$\mathbf{x}_{0} \rmv\in\rmv \mathcal{X}^{\text{o}}$},
the Bhattacharyya bound in \eqref{equ_bhattacharyya_bound} is obtained from \eqref{equ_lower_bound_variance_projection} by using the 
\vspace*{-2mm}
subspace
\[
%% \label{equ_def_bhatt}
\mathcal{U}_{\text{\emph{B}}} \ist\triangleq\, \linspan \rmv\big\{ \{ \basisfuncRKHS_{0} (\cdot) \} \cup {\{ \basisfuncRKHS_{l}(\cdot) \}}_{l \in [L]} \big\} \,,
\vspace{-3.5mm}
\] 
with the functions 
%% $v_{0} (\cdot ) \triangleq R_{\mathcal{E}\rmv,\mathbf{x}_{0}}(\cdot, \mathbf{x}_{0}) \in \mathcal{H}(\mathcal{M})$ and 
\be
\label{equ_def_vectors_ccrb}
\basisfuncRKHS_{0} (\cdot ) \ist\triangleq\ist R_{\mathcal{E}\rmv,\mathbf{x}_{0}}(\cdot\ist\ist, \mathbf{x}_{0}) \in \mathcal{H}_{\mathcal{E}\rmv,\mathbf{x}_{0}} \,, \qquad\;
\basisfuncRKHS_{l}(\cdot) \ist\triangleq\ist \frac{\partial^{{\mathbf{p}}_{l}} \rmv R_{\mathcal{E}\rmv,\mathbf{x}_{0}}(\cdot\ist\ist,\mathbf{x})}{\partial \mathbf{x}^{{\mathbf{p}}_{l}}} \bigg|_{\mathbf{x} = \mathbf{x}_{0}} 
  \!\!\in \mathcal{H}_{\mathcal{E}\rmv,\mathbf{x}_{0}} \,, \quad l \rmv\in\rmv [L] \,.  
\vspace{3mm}
\ee
\end{theorem}

%% \begin{proof}
%% \cite[Section 4.4.3]{JungPHD}
%% \end{proof}

While the RKHS interpretation of the Bhattacharyya bound has been presented previously in \cite{Duttweiler73b} for a specific estimation problem, 
the above result holds for general estimation problems. We note that the bound tends to become higher (tighter) if $L$ is increased in the sense that
additional functions $\basisfuncRKHS_{l}(\cdot)$ are used (i.e., in addition to the functions already used).
Finally, we note that the CRB subspace $\mathcal{U}_{\text{CR}}$ in Theorem \ref{th_crb} is obtained as a special case of the Bhattacharyya bound 
subspace $\mathcal{U}_{\text{B}}$ by setting $L \rmv=\rmv N$, $m \rmv=\rmv 1$, and 
%% multi-indices 
$\mathbf{p}_{l} = \mathbf{e}_{l}$ in \eqref{equ_def_vectors_ccrb}. 

\vspace{-1mm}

\subsection{Hammersley-Chapman-Robbins Bound} 
\label{sec_HCRB}
%%%%%%%%%%%%%%%%%%%%%%%%%%%%%%%%%%%%%%%%%%%%%%%%%%%%%%%%%%

%% Despite its popularity, a 
A drawback of the CRB and the Bhattacharyya bound is that they exploit only the local structure of an estimation problem $\mathcal{E}$ 
around a specific point $\mathbf{x}_{0} \in \mathcal{X}^{\text{o}}$ \cite{LowerBoundAbel}. 
%% \begin{example}
As an illustrative example, consider two different estimation problems $\mathcal{E}_{1} = \big(\mathcal{X}_{1},f(\mathbf{y}; \mathbf{x}),g(\cdot)\big)$ and 
$\mathcal{E}_{2} = \big(\mathcal{X}_{2},f(\mathbf{y}; \mathbf{x}),g(\cdot)\big)$ 
with
%% that share 
the same statistical model $f(\mathbf{y}; \mathbf{x})$ and parameter function $g(\cdot)$ but 
%% are defined for two 
different parameter sets $\mathcal{X}_{1}$ and $\mathcal{X}_{2}$. These parameter sets are assumed to be open balls 
centered at $\mathbf{x}_{0}$ with different radii $r_{1}$ and $r_{2}$, i.e., $\mathcal{X}_{1} = \mathcal{B}(\mathbf{x}_{0},r_{1})$ and $\mathcal{X}_{2} = \mathcal{B}(\mathbf{x}_{0},r_{2})$ with $r_{1} \neq r_{2}$. Then the CRB at $\mathbf{x}_{0}$ for both estimation problems will be identical, 
irrespective of the 
%% precise 
values of $r_{1}$ and $r_{2}$, and similarly for the Bhattacharyya bound. Thus, these bounds 
%% ignore in some sense 
do not take into account a part of the information contained in the parameter set $\mathcal{X}$. 
%% \end{example}
The Barankin bound, on the other hand, exploits the full information carried by the parameter set $\mathcal{X}$ since it is the tightest possible lower bound on the estimator variance. 
However, the Barankin bound is difficult to evaluate in general.
%%  the exact Barankin bound in \eqref{equ_barankin_bound}. 

The Hammersley-Chapman-Robbins bound (HCRB) \cite{GormanHero,ChapmanRobbins51,Hammersley50} is a lower bound on the estimator variance that 
takes into account the global structure of the estimation problem associated with the entire parameter set $\mathcal{X}$.
%% , i.e., it depends on the parameter set $\mathcal{X}$.
%%  not only through its local structure as do the CRB and the Bhattacharyya bound. 
It can be evaluated much more easily than the Barankin bound, and it does not require the estimation problem to be regular. 
Based on a suitably chosen set of ``test points'' $\{{\mathbf{x}}_{1},\ldots,{\mathbf{x}}_{L}\} \subseteq \mathcal{X}$, the HCRB
%% $\{{\mathbf{x}}_{1},\ldots,\tilde{\mathbf{x}}_{L}\} \subseteq \mathcal{X}$, with a finite number $L \in \mathbb{N}$ of test-points. 
states that \cite{GormanHero}
\begin{equation}
\label{equ_HCRB} 
v(\hat{g}(\cdot);\mathbf{x}_{0}) \,\geq\, \mathbf{m}^T\rmv(\mathbf{x}_{0}) \ist
  \mathbf{V}^{\dagger}(\mathbf{x}_{0}) \, \mathbf{m}(\mathbf{x}_{0}) \,,
\end{equation}  
where the vector $\mathbf{m}(\mathbf{x}_{0}) \rmv\in\rmv \mathbb{R}^{L}\rmv$ and the matrix $\mathbf{V}(\mathbf{x}_{0}) \rmv\in\rmv \mathbb{R}^{L \times L}\rmv$ 
are given elementwise by 
$\big( \mathbf{m}(\mathbf{x}_{0}) \big)_l \triangleq \gamma({\mathbf{x}}_{l}) - \gamma(\mathbf{x}_{0})$ 
\vspace{-1mm}
and 
\[ 
\big( \mathbf{V}(\mathbf{x}_{0}) \big)_{l,l'} \ist\triangleq\, \expect_{\mathbf{x}_{0}}
 \bigg\{  \frac{ [ f(\mathbf{y};{\mathbf{x}}_{l}) \rmv\rmv-\! f(\mathbf{y};\mathbf{x}_{0}) ] \ist 
   [ f(\mathbf{y};{\mathbf{x}}_{l'}) \rmv\rmv-\! f(\mathbf{y};\mathbf{x}_{0}) ] }{f^2(\mathbf{y}; \mathbf{x}_{0})} \bigg \} \,,
\vspace{-2mm}
\] 
respectively. 

The following result is proved in \cite[Section 4.4.4]{JungPHD}.

\vspace{1mm}

\begin{theorem} 
\label{th_hcrb}
The HCRB in \eqref{equ_HCRB}, with test points ${\{ {\mathbf{x}}_{l} \}}_{l \in [L]} \subseteq \mathcal{X}$, 
is obtained from \eqref{equ_lower_bound_variance_projection} by using the subspace 
\[
%% \label{equ_def_subspace_HCRB}
\mathcal{U}_{\text{\emph{HCR}}}\ist\triangleq\, \linspan \rmv\big\{ \{ \basisfuncRKHS_{0} (\cdot) \} \cup {\{ \basisfuncRKHS_{l}(\cdot) \}}_{l \in [L]} \big\} \,,
\vspace{-3mm}
\] 
with the 
\vspace*{-.5mm}
functions 
\[
\basisfuncRKHS_{0} (\cdot ) \ist\triangleq\ist R_{\mathcal{E}\rmv,\mathbf{x}_{0}}(\cdot\ist\ist, \mathbf{x}_{0}) \in \mathcal{H}_{\mathcal{E}\rmv,\mathbf{x}_{0}} \,, \qquad\;
\basisfuncRKHS_{l}(\cdot) \ist\triangleq\ist R_{\mathcal{E}\rmv,\mathbf{x}_{0}}(\cdot\ist\ist, {\mathbf{x}}_{l}) \rmv-\rmv R_{\mathcal{E}\rmv,\mathbf{x}_{0}}(\cdot\ist\ist, \mathbf{x}_{0}) \,, \quad l \rmv\in\rmv [L] \,.  
\vspace{2.5mm}
\]
\end{theorem}

The HCRB tends to become higher (tighter) if $L$ is increased in the sense that test points $\mathbf{x}_{l}$ or, equivalently, functions $\basisfuncRKHS_{l}(\cdot)$
are added to those already used.

%% \begin{proof}
%% \cite[Section 4.4.4]{JungPHD}
%% \end{proof}

\subsection{Lower Semi-continuity of the Barankin Bound} 
\label{sec_lower_semi_cont_bar_bound}
%%%%%%%%%%%%%%%%%%%%%%%%%%%%%%%%%%%%%%%%%%%%%%%%%%%%%%%%%%
%% PHDAJ -- p. 57 

For a given estimation problem $\mathcal{E} = \scalarestproblem$ and a prescribed bias function $c(\cdot)$, it is sometimes of interest to characterize not only 
the minimum achievable variance $M(c(\cdot), \mathbf{x}_{0})$ at a single parameter vector $\mathbf{x}_{0} \rmv\in\rmv \mathcal{X}$ but also how 
%% the minimum achievable variance 
$M(c(\cdot), \mathbf{x}_{0})$ changes if $\mathbf{x}_{0}$ is varied. 
%% For each value of $\mathbf{x}_{0}$, we can define the kernel $R_{\mathcal{E}\rmv,\mathbf{x}_{0}}(\cdot\ist\ist,\cdot)$ 
%% according to \eqref{equ_def_kernel_est_problem} and the associated RKHS $\mathcal{H}_{\mathcal{E}\rmv,\mathbf{x}_{0}}$. 
The following result 
%% on the dependence of the minimum achievable variance $M(c(\cdot),\mathbf{x}_{0})$ on $\mathbf{x}_{0} \in \mathcal{X}$ 
is proved in Appendix \ref{app_proof_thm_lower_semi_cont_RKHS}.

\vspace{1mm}

\begin{theorem} 
\label{thm_lower_semi_cont_varying_kernel}
Consider an estimation problem $\mathcal{E}= \scalarestproblem$ with parameter set $\mathcal{X} \rmv\!\subseteq\! \mathbb{R}^{N}\rmv$ and a prescribed bias function 
$c(\cdot) \!: \mathcal{X} \!\rightarrow \mathbb{R}$ that is valid at all $\mathbf{x}_{0} \!\in\rmv \mathcal{C}$ for some open set $\mathcal{C} \!\subseteq\! \mathcal{X}$
and for which the associated prescribed mean function $\gamma(\cdot) = c(\cdot) + g(\cdot)$
%% : \mathcal{X} \!\rightarrow \mathbb{R}$ 
is a continuous function 
%% of $\mathbf{x}$ 
on $\mathcal{C}$. 
%% We denote by $\gamma(\cdot): \mathcal{X} \rightarrow \mathbb{R} : \gamma(\mathbf{x}) = c(\mathbf{x}) + g(\mathbf{x})$ the prescribed mean function. 
Furthermore assume that
%% the prescribed mean function $\gamma(\mathbf{x}) = c(\mathbf{x}) + g(\mathbf{x})$ is a continuous function of $\mathbf{x}$ over $\mathcal{C}$ and 
%% In what follows, we will assume that the estimation problem $\mathcal{E}$, 
%% in particular its statistical model ${\{ f(\mathbf{y}; \mathbf{x}) \}}_{\mathbf{x} \in \mathcal{X}}$, is such that 
for any fixed $\mathbf{x}_{1}, \mathbf{x}_{2} \rmv\in\rmv \mathcal{X}$, $R_{\mathcal{E}\rmv,\mathbf{x}_{0}}(\mathbf{x}_{1}, \mathbf{x}_{2})$
%% : \mathcal{X} \times \mathcal{X} \rightarrow \mathbb{R}$ 
is continuous with respect to
%%  the parameter 
$\mathbf{x}_{0}$ on $\mathcal{C}$, i.e., 
\begin{equation} 
\label{equ_def_cont_varying_kernel}
\lim\limits_{\mathbf{x}'_{0} \rightarrow \mathbf{x}_{0}} \! R_{\mathcal{E}\rmv,\mathbf{x}'_{0}}(\mathbf{x}_{1}, \mathbf{x}_{2}) 
  = R_{\mathcal{E}\rmv,\mathbf{x}_{0}}(\mathbf{x}_{1}, \mathbf{x}_{2})\,, \quad\;  \forall\, \mathbf{x}_{0}  \in \mathcal{C} \,, \; \forall\, \mathbf{x}_{1}, \mathbf{x}_{2}  \in \mathcal{X} \ist.
\vspace{.5mm}
\end{equation}
Then, the minimum achievable variance $M(c(\cdot), \mathbf{x})$, 
%% exists, i.e.\ is finite, for every $\mathbf{x} \in \mathcal{C}$. Furthermore, 
viewed as a function of $\mathbf{x}$, 
%% $M(c(\cdot), \mathbf{x})$ 
is lower semi-continuous on $\mathcal{C}$.
\end{theorem}

\vspace{1mm}

\begin{figure}
\vspace{-1mm}
\centering
\psfrag{x0}[c][c][.9]{\uput{0mm}[270]{0}{$\!\rmv\mathbf{x}_{0}$}}
\psfrag{x}[c][c][.9]{\uput{0mm}[270]{0}{$\!\!\mathbf{x}$}}
\psfrag{fx}[c][c][.9]{\uput{0mm}[0]{0}{\raisebox{-9.5mm}{\hspace{-10.5mm}$f(\mathbf{x})$}}}
\centering
\hspace*{1.5mm}\includegraphics[height=4cm,width=7.3cm]{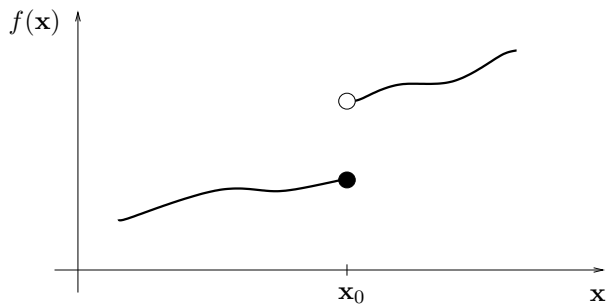}
%\vspace{-3.5mm}
\caption{Graph of a function that is lower semi-continuous at $\mathbf{x}_{0}$. The solid dot indicates the function value $f(\mathbf{x}_{0})$.} %%  and $\sigma^{2}=1$
%\label{fig_bounds_1}
\label{fig_lower_upper_semi_cont}
%% \vspace*{-1.5mm}
\end{figure}

A schematic illustration of a lower semi-continuous function is given in Fig.\ \ref{fig_lower_upper_semi_cont}.
The application of Theorem \ref{thm_lower_semi_cont_varying_kernel} to the estimation problems considered in \cite{ZvikaCRB}---corresponding to the
linear/Gaussian model with a sparse parameter vector---allows us to conclude that the ``sparse CRB'' introduced in \cite{ZvikaCRB} cannot be maximally tight, 
i.e., it is not equal to the minimum achievable variance. Indeed, the sparse CRB 
derived in \cite{ZvikaCRB} is in general a strictly upper 
semi-continuous\footnote{A %%%%%%%
function is said to be \emph{strictly upper semi-continuous} if it is upper semi-continuous but not 
continuous.} %%%%%%%%
function of the parameter vector $\mathbf{x}$, whereas the minimum achievable variance $M(c(\cdot),Ê\mathbf{x})$ is 
%% necessarily 
lower semi-continuous according to Theorem \ref{thm_lower_semi_cont_varying_kernel}. Since a function cannot be simultaneously strictly upper semi-continuous
and lower semi-continuous, the sparse CRB cannot be equal to $M(c(\cdot),Ê\mathbf{x})$.

%%  \begin{proof} 
%% Appendix \ref{app_proof_thm_lower_semi_cont_RKHS}
%% \end{proof} 

%%%%%%%%%%%%%%%%%%%%%%%%%%%%%%%%%%%%%%%%%%%%%%%%%%%%%%%%%%
\section{Sufficient Statistics}
\label{sec_suff_stat_RKHS}
%%%%%%%%%%%%%%%%%%%%%%%%%%%%%%%%%%%%%%%%%%%%%%%%%%%%%%%%%%

\vspace{.5mm}

For some estimation problems $\mathcal{E} = \scalarestproblem$, the observation $\mathbf{y} \in \mathbb{R}^{M}$ contains 
information that is irrelevant to $\mathcal{E}$,
%% the problem of estimating $g(\mathbf{x})$, 
and thus $\mathbf{y}$ can be compressed in some sense. 
Accordingly, let us replace $\mathbf{y}$ by a transformed observation 
%% (known as a \emph{statistic}) 
$\mathbf{z} = \mathbf{t}(\mathbf{y}) \rmv\in\rmv \mathbb{R}^{K}\rmv$, 
with a deterministic mapping $\mathbf{t}(\cdot)\rmv: \mathbb{R}^{M} \!\rightarrow \mathbb{R}^{K}\rmv$.
A compression is achieved if $K \!<\! M$. 
%We will make the weak 
%technical assumption that the mapping $\mathbf{T}(\cdot)$ is such that the modified observation $\mathbf{z}= \mathbf{T}(\mathbf{y})$, which is a random vector, 
%possesses a pdf \cite{papoulis,BillingsleyProbMeasure,AshProbMeasure}. 
%For clarity, we will denote the evaluation of the pdf of the random vector $\mathbf{z} = \mathbf{T}(\mathbf{y})$ at the specific vector $\mathbf{z}'  \in \mathbb{R}^{K}$ by $f_{\mathbf{z}}(\mathbf{z}'; \mathbf{x})$.
Any transformed observation $\mathbf{z} = \mathbf{t}(\mathbf{y})$ is termed a \emph{statistic}, and in particular it is said to be a \emph{sufficient statistic}
if it preserves all the information that is relevant to $\mathcal{E}$
%% the problem of estimating $g(\mathbf{x})$
\cite{Kullback68informationtheory, BillingsleyProbMeasure,LC,HalmosRadonNikodymSuffStat,IbragimovBook,kay}.
In particular, a sufficient statistic 
%% $\mathbf{z} = \mathbf{t}(\mathbf{y})$ 
preserves the minimum achievable variance (Barankin bound) $M(c(\cdot),\mathbf{x}_{0})$.
In the following, the mapping $\mathbf{t}(\cdot)$ will be assumed to be measurable.

%% In what follows we will 
For a given reference parameter vector $\mathbf{x}_{0} \rmv\in\rmv \mathcal{X}$, we consider estimation problems $\mathcal{E}= \scalarestproblem$ 
%% satisfying the assumptions of Section \ref{SecRKHS-basics}, i.e., in particular we assume the existence of 
for which there exists a 
%% base 
dominating measure $\mu_{\mathcal{E}}$ such that the % $\mu_{\mathcal{E}}$ \cite{HalmosRadonNikodymSuffStat} with respect to which 
pdfs ${\{ f(\mathbf{y}; \mathbf{x}) \}}_{\mathbf{x} \in \mathcal{X}}$ 
%% of the statistical model 
are well defined with respect to $\mu_{\mathcal{E}}$ and condition \eqref{equ_f_y_x_0_not_vanish} is satisfied.
The Neyman-Fisher factorization theorem \cite{LC,HalmosRadonNikodymSuffStat,kay,IbragimovBook} then states that
the statistic $\mathbf{z} = \mathbf{t}(\mathbf{y})$ is sufficient for $\mathcal{E}= \scalarestproblem$ if and only if $f(\mathbf{y}; \mathbf{x})$ can be factored as 
%is 
%a sufficient statistic for $\mathcal{E}$ if and only if the statistical model $f(\mathbf{y}; \mathbf{x})$ can be factored as 
\begin{equation} 
\label{equ_suff_stat_factor} 
f(\mathbf{y}; \mathbf{x}) \ist=\ist  h(\mathbf{t}(\mathbf{y}); \mathbf{x}) \, k(\mathbf{y}) \,,
\end{equation}  
where $h(\cdot\ist\ist;\mathbf{x})$ and $k(\cdot)$ are nonnegative functions and the function $k(\cdot)$ does not depend on $\mathbf{x}$. 
Relation \eqref{equ_suff_stat_factor} has to be satisfied for every $\mathbf{y} \!\in\! \mathbb{R}^M\rmv$ except for a set of measure zero with respect to 
the dominating measure $\mu_{\mathcal{E}}$. 

%% \begin{theorem}
%% \label{thm_sufficient_stat_factor}
%% \end{theorem}
%% \begin{proof}
%% \cite{LC,kay,IbragimovBook} 
%% \end{proof} 

%% \newpage %%%%%%%%%

The probability measure on $\mathbb{R}^{K}$ (equipped with the system of $K$-dimensional Borel sets, cf.\ \cite[Section 10]{BillingsleyProbMeasure})
that is induced by the random vector $\mathbf{z} = \mathbf{t}(\mathbf{y})$ is obtained as $\mu^{\mathbf{z}}_{\mathbf{x}} = \mu^{\mathbf{y}}_{\mathbf{x}} \mathbf{t}^{-1}$ \cite{IbragimovBook,HalmosRadonNikodymSuffStat}.
According to Section \ref{SecRKHS-basics}, under condition \eqref{equ_f_y_x_0_not_vanish},
%% assumptions on $\mathcal{E}$ stated there, 
% the existence of the Radon-Nikoydym derivative of the measures ${\{ \mu^{\mathbf{y}}_{\mathbf{x}}Ê\}}_{\mathbf{x} \in \mathcal{X}}$ 
%with respect to the measure $\mu^{\mathbf{y}}_{\mathbf{x}_{0}}$ is guaranteed, and it follows 
%from the Radon-Nikodym theorem \cite{HalmosRadonNikodymSuffStat,BillingsleyProbMeasure} that 
%each measure of ${\{ \mu^{\mathbf{y}}_{\mathbf{x}}Ê\}}_{\mathbf{x} \in \mathcal{X}}$ possesses a density, i.e., the Radon Nikodym derivative, w.r.t. the measure $\mu^{\mathbf{y}}_{\mathbf{x}_{0}}$. 
%Therefore, according to \cite[]{}, 
the measure $\mu^{\mathbf{y}}_{\mathbf{x}_{0}}$ dominates the measures ${\{ \mu^{\mathbf{y}}_{\mathbf{x}}Ê\}}_{\mathbf{x} \in \mathcal{X}}$. This, in turn, 
implies via \cite[Lemma 4]{HalmosRadonNikodymSuffStat} that the measure $\mu^{\mathbf{z}}_{\mathbf{x}_{0}}$ dominates the measures 
${\{ \mu^{\mathbf{z}}_{\mathbf{x}}Ê\}}_{\mathbf{x} \in \mathcal{X}}\ist$, and therefore that, for each $\mathbf{x} \rmv\in\rmv \mathcal{X}$, 
there exists a pdf $f(\mathbf{z}; \mathbf{x})$ with respect to the measure $\mu^{\mathbf{z}}_{\mathbf{x}_{0}}$.
This pdf is given by the following result. (Note that we do not assume condition \eqref{equ_corr_likelihood_finite}.)

\vspace{1mm}

\begin{lemma}
\label{lem_likelihood_factorization_factor_suff_stat}
Consider an estimation problem $\mathcal{E}=\scalarestproblem$ 
\textcolor{red}{satisfying \eqref{equ_f_y_x_0_not_vanish}, i.e.,}
%% satisfying the assumptions of Section \ref{SecRKHS-basics} but not necessarily }. In particular, we assume that 
%for which there exists a dominating measure $\mu_{\mathcal{E}}$, 
which is such that the Radon-Nikodym derivative of $\mu_{\mathbf{x}}^{\mathbf{y}}$ 
with respect to $\mu_{\mathbf{x}_{0}}^{\mathbf{y}}$ is well defined and given by the likelihood ratio \textcolor{red}{$\llr (\mathbf{y}, \mathbf{x})$}.
%???whose statistical model $f(\mathbf{y}; \mathbf{x})$ satisfies \eqref{equ_corr_likelihood_finite}.
Furthermore consider a sufficient statistic $\mathbf{z} = \mathbf{t}(\mathbf{y})$ for $\mathcal{E}$.
Then, the pdf of $\mathbf{z}$ with respect to the dominating measure $\mu^{\mathbf{z}}_{\mathbf{x}_{0}}$ is given by 
\begin{equation}
\label{equ_pdf_suff_stat_factor_neyman_fisher}
f(\mathbf{z}; \mathbf{x}) \ist=\ist \frac{h(\mathbf{z}; \mathbf{x})}{h(\mathbf{z}; \mathbf{x}_{0})} \,,
\end{equation}
where the function $h(\mathbf{z}; \mathbf{x})$ is obtained from the factorization \eqref{equ_suff_stat_factor}.
\end{lemma}

\vspace{1mm}

\emph{Proof}:\,
%% \begin{proof}
The pdf $f(\mathbf{z}; \mathbf{x})$ of $\mathbf{z}$ with respect to $\mu_{\mathbf{x}_{0}}^{\mathbf{z}}$ is defined by the relation 
\be
\label{equ_def_prob-A}
\expect_{\mathbf{x}_{0}} \big\{ I_{\mathcal{A}}(\mathbf{z}) \ist f(\mathbf{z}; \mathbf{x}) \big\} \ist=\ist \Prob_{\!\mathbf{x}} \{  \mathbf{z} \rmv\in\rmv \mathcal{A} \} \,,
\ee
which has to be satisfied for every measurable set $\mathcal{A}Ê\subseteq \mathbb{R}^{K}$ \cite{BillingsleyProbMeasure}. 
Denoting the pre-image of $\mathcal{A}$ under the mapping $\mathbf{t}(\cdot)$ by 
$\mathbf{t}^{-1} (\mathcal{A}) \triangleq \big\{ \mathbf{y} \big| \mathbf{t}(\mathbf{y}) \rmv\in\rmv \mathcal{A} \big \} \subseteq \mathbb{R}^{M}\rmv$, we 
\vspace{1mm}
have 
%% that 
\begin{align}
%% \expect_{\mathbf{x}} \big\{ I_{\mathcal{A}}(\mathbf{z}) \ist f(\mathbf{z}; \mathbf{x}) \big\} 
%% &\eq 
\expect_{\mathbf{x}_{0}} \bigg\{ I_{\mathcal{A}}(\mathbf{z}) \ist\ist \frac{h(\mathbf{z}; \mathbf{x})}{h(\mathbf{z}; \mathbf{x}_{0})} \bigg\}  
%% \nonumber \\[1.5mm]Ê
&  \stackrel{(a)}{\,=\,} 
\expect_{\mathbf{x}_{0}} \bigg\{ I_{\mathcal{A}}(\mathbf{t}(\mathbf{y})) \ist\ist \frac{h(\mathbf{t}(\mathbf{y}); \mathbf{x})}{h(\mathbf{t}(\mathbf{y}); \mathbf{x}_{0})} \bigg\}  \nonumber \\[1.5mm]
&\eq \expect_{\mathbf{x}_{0}} \bigg\{ I_{\mathbf{t}^{-1}(\mathcal{A}) }(\mathbf{y}) \ist\ist \frac{h(\mathbf{t}(\mathbf{y}); \mathbf{x})}{h(\mathbf{t}(\mathbf{y}); \mathbf{x}_{0})} \bigg\}  \nonumber \\[1.5mm]
& \textcolor{red}{ \stackrel{\eqref{equ_suff_stat_factor},\eqref{equ_def_likelihood}}{\,=\,} 
\expect_{\mathbf{x}_{0}} \bigg\{ I_{\mathbf{t}^{-1}(\mathcal{A}) }(\mathbf{y}) \ist\ist \llr (\mathbf{y}, \mathbf{x})  \bigg\} }\nonumber \\[1.5mm]
& \stackrel{(b)}{\eq} \Prob_{\!\mathbf{x}} \{ \mathbf{y} \rmv\in\rmv \mathbf{t}^{-1} (\mathcal{A}) \}\nonumber \\[1.5mm]
&\eq\Prob_{\!\mathbf{x}} \{ \mathbf{z} \rmv\in\rmv\rmv \mathcal{A} \} \,, 
\label{equ_E-I-h}
\end{align}
where step $(a)$ follows from \cite[Theorem 16.12]{BillingsleyProbMeasure} and $(b)$ is due to the fact that the Radon-Nikodym derivative of 
$\mu_{\mathbf{x}}^{\mathbf{y}}$ with respect to $\mu_{\mathbf{x}_{0}}^{\mathbf{y}}$ is 
given by \textcolor{red}{$\llr (\mathbf{y}, \mathbf{x})$ (cf.\ \eqref{equ_def_likelihood})}, as explained in Section \ref{SecRKHS-basics}.
Comparing \eqref{equ_E-I-h} with \eqref{equ_def_prob-A}, we conclude that $\frac{h(\mathbf{z}; \mathbf{x})}{h(\mathbf{z}; \mathbf{x}_{0})} = f(\mathbf{z}; \mathbf{x})$ 
up to differences on a set of measure zero (with respect to $\mu_{\mathbf{x}_{0}}^{\mathbf{z}}$).
Note that because we require $\mathbf{t}(\cdot)$ to be a measurable mapping,
%%  from $\mathbb{R}^{M}$ into $\mathbb{R}^{K}\rmv$, 
it is guaranteed that the set $\mathbf{t}^{-1}(\mathcal{A}) = \big\{ \mathbf{y} \big| \mathbf{t}(\mathbf{y}) \rmv\in\rmv \mathcal{A} \big \}$ is measurable for any 
measurable set 
\vspace{3mm}
$\mathcal{A} \subseteq \mathbb{R}^{K}\rmv$. 
%???\text{WHERE did we need assumption }\eqref{equ_corr_likelihood_finite}???
 \hfill $\Box$
%and the definition of the expectation 
%as an integral \cite[Section 21]{BillingsleyProbMeasure}. 
%% \end{proof}

Consider next an estimation problem $\mathcal{E}=\scalarestproblem$ satisfying \eqref{equ_corr_likelihood_finite}, so that the kernel 
$R_{\mathcal{E}, \mathbf{x}_{0}}(\cdot\ist\ist,\cdot)$ exists according to \eqref{equ_def_kernel_est_problem}.
Let $\mathbf{z} = \mathbf{t}(\mathbf{y})$ be a sufficient statistic. We can then define 
the modified estimation problem $\mathcal{E}' \triangleqÊ\big(\mathcal{X}, f(\mathbf{z}; \mathbf{x}), g(\cdot)\big)$, which is based on the observation $\mathbf{z}$ and 
whose statistical model is given by the pdf $f(\mathbf{z}; \mathbf{x})$ (cf.\ \eqref{equ_pdf_suff_stat_factor_neyman_fisher}). 
The following theorem states that the RKHS associated with $\mathcal{E}'$ equals the RKHS associated with $\mathcal{E}$.

\vspace{1mm}

%% \textcolor{red}{
\begin{theorem}
\label{thm_suff_stac_est_problem_RKHS_same}
Consider an estimation problem $\mathcal{E}= \scalarestproblem$ satisfying \eqref{equ_corr_likelihood_finite} and a reference parameter vector $\mathbf{x}_{0} \!\in\! \mathcal{X}$. 
For a sufficient statistic $\mathbf{z} = \mathbf{t}(\mathbf{y})$, consider
%% we define 
the modified estimation problem $\mathcal{E}' = \big(\mathcal{X}, f(\mathbf{z}; \mathbf{x}), g(\cdot) \big)$. 
%with dominating measure $\mu_{\mathcal{E}'}=\mu^{\mathbf{z}}_{\mathbf{x}_{0}}$. 
Then, 
%% the modified estimation problem 
$\mathcal{E}'$ also satisfies \eqref{equ_corr_likelihood_finite} and furthermore
%% the RKHS $\mathcal{H}_{\mathcal{E}'\!,\mathbf{x}_{0}}$ coincides with the RKHS $\mathcal{H}_{\mathcal{E}\rmv,\mathbf{x}_{0}}$,
%%  of the original estimation problem $\mathcal{E}$, 
%% and 
$R_{\mathcal{E}'\!, \mathbf{x}_{0}}(\cdot \ist\ist, \cdot) = R_{\mathcal{E}, \mathbf{x}_{0}}(\cdot \ist\ist, \cdot)$ and 
$\mathcal{H}_{\mathcal{E}'\!,\mathbf{x}_{0}} = \mathcal{H}_{\mathcal{E}\rmv,\mathbf{x}_{0}}$. 
\vspace{1mm}
\end{theorem}
%% }

\emph{Proof}:\,
%The kernel $R_{\mathcal{E}'\!, \mathbf{x}_{0}}(\cdot \ist\ist, \cdot)$ is obtained as
We 
\vspace{-2mm}
have
%% Observe that 
%% Invoking Lemma \ref{lem_likelihood_factorization_factor_suff_stat} yields 
\begin{align}
\label{equ_coincide_likelihood_corr_suff_stat}
 R_{\mathcal{E}, \mathbf{x}_{0}}(\mathbf{x}_{1}, \mathbf{x}_{2}) & \stackrel{\eqref{equ_def_kernel_est_problem}}{\eq} 
 \expect_{\mathbf{x}_{0}} 
  \textcolor{red}{\big\{\llr (\mathbf{y}, \mathbf{x}_{1}) \llr (\mathbf{y}, \mathbf{x}_{2})  \big\} }  \nonumber \\[.5mm]
   &  \stackrel{\eqref{equ_suff_stat_factor},\textcolor{red}{\eqref{equ_def_likelihood}}}{\eq}
 \expect_{\mathbf{x}_{0}} \bigg\{  \frac{h(\mathbf{t}(\mathbf{y}); \mathbf{x}_{1}) \ist h(\mathbf{t}(\mathbf{y}); \mathbf{x}_{2})}{h^2(\mathbf{t}(\mathbf{y}); \mathbf{x}_{0})}\bigg\}  \nonumber \\[1.5mm]
  & \stackrel{(a)}{\eq} \expect_{\mathbf{x}_{0}} 
  \bigg\{  \frac{h(\mathbf{z}; \mathbf{x}_{1}) \ist h(\mathbf{z}; \mathbf{x}_{2})}{h^2(\mathbf{z}; \mathbf{x}_{0})}\bigg\}  \nonumber \\[1.5mm]
  & \stackrel{\eqref{equ_pdf_suff_stat_factor_neyman_fisher}}{\eq} 
  \expect_{\mathbf{x}_{0}} \bigg\{  \frac{f(\mathbf{z}; \mathbf{x}_{1}) \ist f(\mathbf{z}; \mathbf{x}_{2})}{f^2(\mathbf{z}; \mathbf{x}_{0})}\bigg\} \nonumber \\[1.5mm]Ê
 & \stackrel{\eqref{equ_def_kernel_est_problem}}{\eq}  
R_{\mathcal{E}'\!, \mathbf{x}_{0}}(\mathbf{x}_{1}, \mathbf{x}_{2})\,,
\end{align}
where, as before, step $(a)$ follows from \cite[Theorem 16.12]{BillingsleyProbMeasure}. % and the definition of the expectation 
%as an integral \cite[Section 21]{BillingsleyProbMeasure}. 
From \eqref{equ_coincide_likelihood_corr_suff_stat}, we conclude that if $\mathcal{E}$ satisfies \eqref{equ_corr_likelihood_finite} then so does $\mathcal{E}'\rmv$. 
Moreover, from $R_{\mathcal{E}'\!, \mathbf{x}_{0}}(\cdot \ist\ist, \cdot) = R_{\mathcal{E}, \mathbf{x}_{0}}(\cdot \ist\ist, \cdot)$ in \eqref{equ_coincide_likelihood_corr_suff_stat},
it follows that $\mathcal{H}_{\mathcal{E}'\!,\mathbf{x}_{0}} =\mathcal{H}(R_{\mathcal{E}'\!,\mathbf{x}_{0}})$ equals
$\mathcal{H}_{\mathcal{E}\rmv,\mathbf{x}_{0}} =\mathcal{H}(R_{\mathcal{E}\rmv,\mathbf{x}_{0}})$.
%!!!By following the chain of equalities in \eqref{equ_coincide_likelihood_corr_suff_stat} in the reverse direction, the fulfillment of condition \eqref{equ_corr_likelihood_finite} for the modified estimation problem $\mathcal{E}Õ$ is verified.!!!
\vspace{3mm}
\hfill $\Box$

Intuitively, one might expect that the RKHS associated with a sufficient statistic should be typically ``smaller'' or ``simpler'' than the RKHS associated with
the original observation, since in general the sufficient statistic is a compressed and ``more concise'' version of the observation. 
However, Theorem \ref{thm_suff_stac_est_problem_RKHS_same} states that the RKHS remains unchanged by this compression. 
One possible interpretation of this fact is
%% an interesting interpretation of Theorem \ref{thm_suff_stac_est_problem_RKHS_same} is 
that the RKHS description of an estimation problem is already
%% intrinsically 
``maximally efficient'' in the sense that it cannot be reduced or simplified by using a compressed (yet sufficiently informative) observation.

%%%%%%%%%%%%%%%%%%%%%%%%%%%%%%%%%%%%%%%%%%%%%%%%%%%%%%%%%%
\section{MVE for the Exponential Family} %%RKHS-based Analysis of MVE for Exponential Families
\label{sec_rkhs_exp_fam}
%%%%%%%%%%%%%%%%%%%%%%%%%%%%%%%%%%%%%%%%%%%%%%%%%%%%%%%%%%

An important class of estimation problems 
%% that are regular in the sense of Definition \ref{post_regular_cond_CRB} 
is defined by statistical models belonging to an exponential family. 
%% It is somehow natural to consider this specific class 
Such models are of considerable interest in the context of MVE because,
%%  One reason is that, 
under mild conditions, the existence of a UMV estimator is guaranteed. 
Furthermore, 
%% Another important fact is that 
any estimation problem 
that admits the existence of an \emph{efficient estimator}, i.e., an estimator whose variance achieves the CRB,
%%  (cf.\ Section \ref{sec_CRB}), 
must be necessarily based on an exponential family \cite[Theorem 5.12]{LC}. 
%% These estimation problems are regular in the sense of Definition \ref{post_regular_cond_CRB}.
%% A prominent example is the linear-Gaussian model.
%% , ???the sparse linear-Gaussian model, 
%% (cf. Section \ref{sec_basic_concepts}) 
%% and the sparse covariance model studied in a RKHS context in \cite{AlexSebICASSP2011}. 
In this section, we will characterize the RKHS for this class and use it to derive lower variance bounds.
%%  of estimation problems.
%%  with observations conforming to the exponential family.
%% Specific members of this class are, e.g., the LGM (cf. Section \ref{sec_basic_concepts}) and the \emph{sparse covariance model} (SCM), 
%% which has been studied in a RKHS context in \cite{AlexSebICASSP2011}. 

%% \newpage %%%%%%%%%

\subsection{\textcolor{red}{Review of the} Exponential Family}
\label{sec_exp_fam}
%%%%%%%%%%%%%%%%%%%%%%%%%%%%%%%%%%%%%%%%%%%%%%%%%%%%%%%%%%

%cf. "Lecture Notes from KTH filed under "$Proof_Complete_Suff_stat_is_minimal$"": 
%
%\begin{theorem} 
%If the natural parameter space $\Omega$ of an exponential family contains
%an open set in $\mathbb{R}^{k}$, then $T(X)$ is a complete sufficient statistic.
%\end{theorem}
%\begin{proof}
%\cite[Theorem 6.2.25]{CasellaBergerStatInf} 
%\end{proof}

An exponential family is defined as the following parametrized set of pdfs ${\{ f(\mathbf{y}; \mathbf{x}) \}}_{\mathbf{x} \in \mathcal{X}}$ (with respect to the Lebesgue measure on $\mathbb{R}^{M}$)
%% that is given by
%% whose general expression is 
\cite{LC,FundmentExpFamBrown,GraphModExpFamVarInfWainJor}:
\[ 
%% \label{equ_exponential_family}
f(\mathbf{y};\mathbf{x}) \eq \exp \rmv\rmv\big( {\bm \phi}^{T}\rmv(\mathbf{y}) \ist \mathbf{u}(\mathbf{x}) \rmv-\rmv A(\mathbf{x})\big) \, h(\mathbf{y}) \,, 
\] 
with the \emph{sufficient statistic} ${\bm \phi}(\cdot):\mathbb{R}^{M} \!\rightarrow \mathbb{R}^{P}$, 
the \emph{parameter function} $\mathbf{u}(\cdot): \mathbb{R}^{N} \!\rightarrow \mathbb{R}^{P}$, 
the \emph{cumulant function} $A(\cdot): \mathbb{R}^{N} \!\rightarrow \mathbb{R}$, 
and the \emph{weight function} $h(\cdot): \mathbb{R}^{M} \!\rightarrow \mathbb{R}$. 
Many well-known statistical models are
%% can be represented as a 
special instances of an exponential family \cite{GraphModExpFamVarInfWainJor}. 
Without loss of generality, we can restrict ourselves to an exponential family in \emph{canonical form} \cite{LC}, for which $P \rmv=\rmv N$ and $\mathbf{u}(\mathbf{x}) = \mathbf{x}$, i.e.,
\begin{equation} 
\label{equ_exponential_family_natural_parametrization}
f^{(A)}(\mathbf{y};\mathbf{x}) \eq \exp \rmv\rmv\big( {\bm \phi}^{T}\rmv(\mathbf{y}) \ist\ist \mathbf{x} \rmv-\rmv A(\mathbf{x})\big) \, h(\mathbf{y}) \,. 
\end{equation} 
Here, the superscript $^{(A)}$ emphasizes the importance of the cumulant function $A(\cdot)$ in the characterization of an exponential family.
In what follows, we assume that the parameter space is chosen as $\mathcal{X} \rmv\subseteq\rmv \mathcal{N}$, where $\mathcal{N} \rmv\rmv\subseteq \mathbb{R}^{N}$ is 
the \emph{natural parameter space} defined as
\[
%% \label{equ_def_natural_parameter_space}
\mathcal{N} \,\triangleq\, \bigg\{ \mathbf{x} \rmv\in\rmv \mathbb{R}^{N} \ist\bigg|  
  \int_{\mathbb{R}^{M}} \! \exp \rmv\rmv\big( {\bm \phi}^{T}\rmv(\mathbf{y}) \ist\ist \mathbf{x} \big) \ist h(\mathbf{y}) \, d \mathbf{y} < \infty \bigg\} \,.
\] 
%% As the notation suggests, an exponential family in canonical form is effectively determined or parametrized by the 
%% cumulant function $A(\cdot): \mathcal{N} \rightarrow \mathbb{R}$.
%% 
% of the form \eqref{equ_exponential_family} with the additional restriction that the parameter function is the identity map. Thus, the statistical model is then given as
%\begin{equation}
%\label{equ_exp_family_parameter_map_identity}
%f(\mathbf{y};\mathbf{x}) = \exp( {\bm \phi}(\mathbf{y})^{T} \mathbf{x} - A(\mathbf{x})) h(\mathbf{y}).
%\end{equation}  
%% Similarly, an exponential family is also completely specified by the sufficient statistic ${\bm \phi}(\mathbf{y})$ and the weight function $h(\mathbf{y})$, since 
From the normalization constraint $\int_{\mathbb{R}^{M}} f^{(A)}(\mathbf{y};\mathbf{x}) \ist\ist d\mathbf{y} =1$, 
it follows that the cumulant function $A(\cdot)$ is determined by the sufficient statistic ${\bm \phi}(\cdot)$ and the weight function $h(\cdot)$ 
\vspace*{.5mm}
as
\[
%% \label{equ_def_cumulant_function}
A(\mathbf{x}) \eq \log \rmv\rmv \bigg( \int_{\mathbb{R}^{M}} \!
  \exp \rmv\rmv\big( {\bm \phi}^{T}\rmv(\mathbf{y}) \ist\ist \mathbf{x} \big) \ist h(\mathbf{y}) \, d \mathbf{y} \rmv \bigg) \ist\ist , \quad\; \mathbf{x} \rmv\in\rmv \mathcal{N} \,.
\vspace{.5mm}
\]
%The interior of the natural parameter space will be denoted as $\mathcal{N}^{o}$, i.e., $\mathcal{N}^{o} \triangleq \big\{ \mathbf{x} \in \mathcal{N} | \exists r>0: \mathcal{B}(\mathbf{x},r) \subseteq \mathcal{N} \big\}$.

The \emph{moment-generating function} of $f^{(A)}(\mathbf{y};\mathbf{x})$
%% $\lambda(\mathbf{x})$ of the exponential family 
is defined 
\vspace{.5mm}
as 
\be 
\label{equ_def_momgen_function}
\lambda(\mathbf{x}) \,\triangleq\, \exp(A(\mathbf{x})) \,=
%% \stackrel{\eqref{equ_def_cumulant_function}}{=} 
\int_{\mathbb{R}^{M}} \! \exp \rmv\rmv\big( {\bm \phi}^{T}\rmv(\mathbf{y}) \ist\ist \mathbf{x} \big) \ist h(\mathbf{y}) \, d \mathbf{y}  \,, \quad\; \mathbf{x} \rmv\in\rmv \mathcal{N}\,.
\vspace{-2mm}
\ee 
Note that 
\begin{equation}
\label{equ_natural_parameter_space_momgen}
\mathcal{N} \ist=\ist \big\{ \mathbf{x} \rmv\in\rmv \mathbb{R}^{N} \ist\big|  
  \lambda(\mathbf{x}) < \infty \big\} \,.
\end{equation} 
Assuming a random vector $\mathbf{y} \sim f^{(A)}(\mathbf{y};\mathbf{x})$, it is known \cite[Theorem 2.2]{FundmentExpFamBrown}, \cite[Proposition 3.1]{GraphModExpFamVarInfWainJor}
that for any $\mathbf{x} \rmv\in\rmv \mathcal{X}^{\text{o}}$ and $\mathbf{p} \rmv\in\rmv \mathbb{Z}_{+}^{N}$, the moments 
$\expect_{\mathbf{x}} \big\{ {\bm \phi}^{\mathbf{p}}(\mathbf{y}) \big\}$ exist, i.e., $\expect_{\mathbf{x}} \big\{ {\bm \phi}^{\mathbf{p}}(\mathbf{y}) \big\} < \infty$, 
and they can be calculated from the partial derivatives of $\lambda(\mathbf{x})$ according to
\begin{equation}
\label{equ_part_derivates_relation_moments_exp_family}
\expect_{\mathbf{x}} \big\{ {\bm \phi}^{\mathbf{p}}(\mathbf{y}) \big\} \ist=\ist \frac{1}{\lambda(\mathbf{x})} \frac{\partial^{\mathbf{p}} \lambda(\mathbf{x})}{\partial \mathbf{x}^{\mathbf{p}}} \,.
\end{equation}
Thus, the partial derivatives $ \frac{\partial^{\mathbf{p}} \lambda(\mathbf{x})}{\partial \mathbf{x}^{\mathbf{p}}}$ exist for any 
$\mathbf{x} \!\in\! \mathcal{X}^{\text{o}}$ and $\mathbf{p} \!\in\! \mathbb{Z}_{+}^{N}$, and for any choice of the sufficient statistic $\phi(\cdot)$ and the
weight function $h(\cdot)$. Moreover, they depend continuously on $\mathbf{x} \rmv\in\rmv \mathcal{X}^{\text{o}}$ \cite{FundmentExpFamBrown,GraphModExpFamVarInfWainJor}. 
%???DON'T WE NEED SOME SMOOTHNESS ASSUMPTIONS ABOUT ${\bm \phi}(\mathbf{y})$ 
%OR $\lambda(\mathbf{x})$ FOR SUCH A STATEMENT???

\subsection{RKHS Associated with an Exponential Family Based MVP}
\label{sec_rkhs_exp_fam_mvp}
%%%%%%%%%%%%%%%%%%%%%%%%%%%%%%%%%%%%%%%%%%%%%%%%%%%%%%%%%%

Consider an estimation problem $\mathcal{E}^{(A)}Ê\triangleq \big(\mathcal{X},f^{(A)}(\mathbf{y};\mathbf{x}),g(\cdot) \big)$ with an exponential family statistical model 
${\{ f^{(A)}(\mathbf{y}; \mathbf{x}) \}}_{\mathbf{x} \in \mathcal{X}}$ as defined in \eqref{equ_exponential_family_natural_parametrization},
%%  with associated natural parameter space $\mathcal{N}$ (cf.\ \eqref{equ_def_natural_parameter_space}), 
and a fixed $\mathbf{x}_{0} \rmv\in\rmv \mathcal{X}$. 
%% For the estimation problem $\mathcal{E}^{(A)}$, we have the identity 
%% \pagebreak %%%%%%
Consider further the RKHS $\mathcal{H}_{\mathcal{E}^{(A))}\rmv\rmv,\mathbf{x}_{0}}$. Its kernel is obtained as 
\begin{align} 
R_{\mathcal{E}^{(A)}\rmv\rmv,\mathbf{x}_{0}}(\mathbf{x}_{1}, \mathbf{x}_{2}) & \ist\stackrel{\eqref{equ_def_kernel_est_problem}}{=}\,
  \expect_{\mathbf{x}_{0}} \bigg\{ \frac{f^{(A)}(\mathbf{y}; \mathbf{x}_{1}) \ist f^{(A)}(\mathbf{y}; \mathbf{x}_{2})}{(f^{(A)}(\mathbf{y}; \mathbf{x}_{0}))^2} \bigg \}
  \label{equ_R_E_A} 
 \\[2mm]
& \ist\stackrel{\eqref{equ_exponential_family_natural_parametrization}}{=}\, 
 \expect_{\mathbf{x}_{0}} \Bigg\{ \frac{ \exp \rmv\rmv\big( {\bm \phi}^{T}\rmv(\mathbf{y}) \ist\ist \mathbf{x}_1 \rmv-\rmv A(\mathbf{x}_1)\big)  
 \exp \rmv\rmv\big( {\bm \phi}^{T}\rmv(\mathbf{y}) \ist\ist \mathbf{x}_2 \rmv-\rmv A(\mathbf{x}_2)\big) }
   { \exp \rmv\rmv\big( 2 \big[ {\bm \phi}^{T}\rmv(\mathbf{y}) \ist\ist \mathbf{x}_0 \rmv-\rmv A(\mathbf{x}_0) \big]\big)} \Bigg\}  \nonumber \\[1.5mm]
 & \ist\eq  \expect_{\mathbf{x}_{0}} \big\{ \exp \rmv\rmv\big( {\bm \phi}^{T}\rmv(\mathbf{y}) \ist\ist (\mathbf{x}_{1} \rmv\rmv+\rmv \mathbf{x}_{2} \rmv-\rmv 2 \ist \mathbf{x}_{0}) 
   - A(\mathbf{x}_{1}) - A(\mathbf{x}_{2}) + 2\ist A(\mathbf{x}_{0})  \big) \big\}  \nonumber \\[2mm] 
 & \ist\stackrel{\eqref{equ_exponential_family_natural_parametrization}}{=}\,  
   \exp \rmv\rmv\big( A(\mathbf{x}_{1}) - A(\mathbf{x}_{2}) + 2\ist A(\mathbf{x}_{0}) \big) \nonumber \\[0mm] 
 &\rule{10mm}{0mm} \times \int_{\mathbb{R}^{M}} \! \exp \rmv\rmv\big( {\bm \phi}^{T}\rmv(\mathbf{y}) \ist\ist (\mathbf{x}_{1} \rmv\rmv+\rmv \mathbf{x}_{2} \rmv-\rmv 2 \ist \mathbf{x}_{0})\big)
   \exp \rmv\rmv\big( {\bm \phi}^{T}\rmv(\mathbf{y}) \ist\ist \mathbf{x}_{0} \rmv-\rmv A(\mathbf{x}_{0}) \big) \ist h(\mathbf{y}) \, d \mathbf{y}
   \nonumber \\[1mm] 
 &\eq\, \exp \rmv\rmv\big( \!-\rmv\rmv A(\mathbf{x}_{1}) \rmv-\rmv A(\mathbf{x}_{2}) + A(\mathbf{x}_{0}) \big) 
   \int_{\mathbb{R}^{M}} \! \exp \rmv\rmv\big( {\bm \phi}^{T}\rmv(\mathbf{y}) \ist\ist (\mathbf{x}_{1} \rmv\rmv+\rmv \mathbf{x}_{2} \rmv-\rmv \mathbf{x}_{0})\big) 
   \ist h(\mathbf{y}) \, d \mathbf{y} \nonumber \\[1mm]
%% & \stackrel{\eqref{equ_def_cumulant_function}}{=}\,  \exp \rmv\rmv\big( \!-\rmv\rmv A(\mathbf{x}_{1}) \rmv-\rmv A(\mathbf{x}_{2}) + A(\mathbf{x}_{0}) \big)
%%   \exp \rmv\rmv\big( A(\mathbf{x}_{1} \rmv\rmv+\rmv \mathbf{x}_{2} \rmv-\rmv \mathbf{x}_{0}) \big) \nonumber \\[2mm]
%% %%   \label{equ_kernel_exp_fam_canonical_form}
& \stackrel{\eqref{equ_def_momgen_function}}{=}\, \frac{\lambda(\mathbf{x}_{1} \rmv\rmv+\rmv \mathbf{x}_{2} \rmv-\rmv \mathbf{x}_{0}) \, \lambda(\mathbf{x}_{0})}
  {\lambda(\mathbf{x}_{1}) \, \lambda(\mathbf{x}_{2})} \,. 
  \label{equ_kernel_exp_family_moment_function} 
%% &\!\!\!\stackrel{\eqref{equ_sum_x_1_x_2_minus_x_0_in_N},\eqref{equ_def_natural_parameter_space}}< \infty \,.
\end{align}
Because \eqref{equ_R_E_A} and \eqref{equ_kernel_exp_family_moment_function} are equal, we see that condition \eqref{equ_corr_likelihood_finite} is satisfied, i.e.,
%% \label{equ_corr_likelihood_finite_exp}
$\expect_{\mathbf{x}_{0}} \Big\{ \frac{f^{(A)}(\mathbf{y}; \mathbf{x}_{1}) \ist f^{(A)}(\mathbf{y}; \mathbf{x}_{2})}{(f^{(A)}(\mathbf{y}; \mathbf{x}_{0}))^2} \Big\} < \infty$ 
for all $\mathbf{x}_{1}, \mathbf{x}_{2} \in \mathcal{X}$, if and only if 
$\frac{\lambda(\mathbf{x}_{1} + \mathbf{x}_{2} - \mathbf{x}_{0}) \, \lambda(\mathbf{x}_{0})} {\lambda(\mathbf{x}_{1}) \, \lambda(\mathbf{x}_{2})} < \infty$ 
for all $\mathbf{x}_{1}, \mathbf{x}_{2} \in \mathcal{X}$. 
Since $\mathbf{x}_{0}
%% , \mathbf{x}_{1},Ê\mathbf{x}_{2} 
\in \mathcal{X} \subseteq \mathcal{N}$, we have 
$\lambda(\mathbf{x}_{0}) \rmv<\rmv \infty$.
%% , $\lambda(\mathbf{x}_{1}) \rmv<\rmv \infty$, and $\lambda(\mathbf{x}_{2}) \rmv<\rmv \infty$.
Furthermore, $\lambda(\mathbf{x}) \not= 0$ for all $\mathbf{x} \rmv\in\rmv \mathcal{X}$.
Therefore, \eqref{equ_corr_likelihood_finite} is satisfied if and only if $\lambda(\mathbf{x}_{1} + \mathbf{x}_{2} - \mathbf{x}_{0}) \rmv<\rmv \infty$.
We conclude that for an estimation problem whose statistical model belongs to an exponential family, condition \eqref{equ_corr_likelihood_finite} 
is equivalent 
\vspace{-2.5mm}
to 
%% In the following, we consider an exponential family statistical model ${\{ f^{(A)}(\mathbf{y}; \mathbf{x}) \}}_{\mathbf{x} \in \mathcal{X}}$ 
%% and ???a fixed $\mathbf{x}_{0} \rmv\in\rmv ???\mathcal{X}$, and we assume that the parameter set $\mathcal{X}$ 
\begin{equation}
\label{equ_sum_x_1_x_2_minus_x_0_in_N} 
\mathbf{x}_{1},Ê\mathbf{x}_{2} \in \mathcal{X} \;\;\Rightarrow\;\; \mathbf{x}_{1} \rmv\rmv+\rmv \mathbf{x}_{2} \rmv-\rmv \mathbf{x}_{0} \in \mathcal{N} \,.
\vspace{.5mm}
\end{equation}
%% Together with \eqref{equ_natural_parameter_space_momgen} and \eqref{equ_sum_x_1_x_2_minus_x_0_in_N}, this expression shows that
%% condition \eqref{equ_corr_likelihood_finite} is satisfied, i.e.,
Furthermore, 
%% it follows from \eqref{equ_kernel_exp_family_moment_function} that the partial derivatives of $R_{\mathcal{E}^{(A)}\rmv\rmv,\mathbf{x}_{0}}(\cdot\ist\ist,\cdot)$
%in \eqref{equ_part-diff} are given by
%\[
%\frac{\partial^{\mathbf{p}_{1}} \partial^{\mathbf{p}_{2}} R(\mathbf{x}_{1}, \mathbf{x}_{2})}{\partial \mathbf{x}_{1}^{\mathbf{p}_{1}} \partial \mathbf{x}_{2}^{\mathbf{p}_{2}} } 
%  \bigg|_{\mathbf{x}_{1}=\mathbf{x}_{2}= \mathbf{x}_0} =\, ??? \,.
%\]
from \eqref{equ_kernel_exp_family_moment_function} and the fact that the partial derivatives $ \frac{\partial^{\mathbf{p}} \lambda(\mathbf{x})}{\partial \mathbf{x}^{\mathbf{p}}}$ exist for any 
$\mathbf{x} \rmv\in\rmv \mathcal{X}^{\text{o}}$ and $\mathbf{p} \rmv\in\rmv \mathbb{Z}_{+}^{N}$ and depend continuously 
on $\mathbf{x} \rmv\in\rmv \mathcal{X}^{\text{o}}\rmv$, we can conclude that the RKHS $\mathcal{H}_{\mathcal{E}^{(A)}\rmv\rmv,\mathbf{x}_{0}}$ 
is differentiable up to any order. We summarize this finding in 

\vspace{1mm}

\begin{lemma}
\label{lem_diff_RKHS_exp_family} 
Consider an estimation problem $\mathcal{E}^{(A)}=\big(\mathcal{X},f^{(A)}(\mathbf{y};\mathbf{x}),g(\cdot)\big)$ associated with an exponential family (cf.\ \eqref{equ_exponential_family_natural_parametrization}) with natural parameter space $\mathcal{N}$. The parameter set $\mathcal{X}$ is assumed to satisfy condition \eqref{equ_sum_x_1_x_2_minus_x_0_in_N} for some reference parameter vector $\mathbf{x}_{0} \rmv\in\rmv \mathcal{X}$. 
Then, the kernel $R_{\mathcal{E}^{(A)}\rmv\rmv,\mathbf{x}_{0}}(\mathbf{x}_{1}, \mathbf{x}_{2})$ 
%% (cf.\ \eqref{equ_def_kernel_est_problem}) 
and the RKHS $\mathcal{H}_{\mathcal{E}^{(A)}\rmv\rmv,\mathbf{x}_{0}}$ 
%associated 
%to $\mathcal{E}^{(\mathcal{A})}$ and the reference parameter vector $\mathbf{x}_{0}$ 
%% \textcolor{red}{
are differentiable up to any order $m$.
\vspace{1mm}
\end{lemma}

Next, by combining Lemma \ref{lem_diff_RKHS_exp_family} with \textcolor{red}{\eqref{equ_der_reproduction_prop}}, we will 
%% allow us to 
derive simple 
%% closed-form 
lower bounds on the variance of estimators with a prescribed bias function.
%???WHERE/WHY DO WE NEED DIFFERENTIABILITY???

\subsection{Variance Bounds for the Exponential Family}
\label{sec_lower_bound_exp_fam}
%%%%%%%%%%%%%%%%%%%%%%%%%%%%%%%%%%%%%%%%%%%%%%%%%%%%%%%%%%

If $\mathcal{X}^{\text{o}}$ is nonempty, the sufficient statistic ${\bm \phi}(\cdot)$ is a \emph{complete} sufficient statistic 
for the estimation problem $\mathcal{E}^{(A)}\rmv$, and thus there exists a UMV estimator $\hat{g}_{\text{\tiny{UMV}}}(\cdot)$ for any valid bias function $c(\cdot)$ \cite[p.\ 42]{LC}.
This UMV estimator is given by the conditional 
expectation\footnote{The %%%%%%%%%%
conditional expectation in \eqref{equ_constr_UMV_cond_expecation} can be taken with respect to the measure $\mu^{\mathbf{y}}_{\mathbf{x}}$ for an arbitrary 
$\mathbf{x} \rmv\in\rmv \mathcal{X}$. Indeed, since ${\bm \phi}(\cdot)$ is a sufficient statistic, 
%% the conditional expectation 
$\expect_{\mathbf{x}}\{ \hat{g}_{0}(\mathbf{y}) | \bm{ \phi }(\mathbf{y}) \}$ yields the same result for every 
$\mathbf{x} \rmv\in\rmv \mathcal{X}$.} %%%%%%%%%%%%
\begin{equation} 
\label{equ_constr_UMV_cond_expecation}
\hat{g}_{\text{\tiny{UMV}}}(\mathbf{y}) =  \expect_{\mathbf{x}}\{ \hat{g}_{0}(\mathbf{y}) | \ist\bm{ \phi }(\mathbf{y}) \} \,,
\end{equation}
where $\hat{g}_{0}(\cdot)$ is any estimator with bias function $c(\cdot)$, i.e., $b(\hat{g}_{0}(\cdot);\mathbf{x}_{0}) = c(\mathbf{x})$ for all $\mathbf{x} \!\in\! \mathcal{X}$. 
The minimum achievable variance $M(c(\cdot),\mathbf{x}_{0})$ is then equal to the variance of $\hat{g}_{\text{\tiny{UMV}}}(\cdot)$ at $\mathbf{x}_{0}$, i.e., 
$M(c(\cdot),\mathbf{x}_{0}) = v(\hat{g}_{\text{\tiny{UMV}}}(\cdot);\mathbf{x}_{0})$ 
%% (cf.\ 
\cite[p.\ 89]{LC}.
However, it may be difficult to actually construct the UMV estimator via \eqref{equ_constr_UMV_cond_expecation} and to calculate its variance. 
In fact, it may be already a difficult task to find an estimator $\hat{g}_{0}(\cdot)$ whose bias function equals $c(\cdot)$. 
Therefore, it is still of interest to find simple closed-form lower bounds on the variance of any estimator with bias $c(\cdot)$.

\vspace{1mm}

\begin{theorem}
\label{thm_lower_bound_exp_family}
Consider an estimation problem $\mathcal{E}^{(A)}=\big(\mathcal{X},f^{(A)}(\mathbf{y};\mathbf{x}),g(\cdot)\big)$ with parameter set $\mathcal{X}$ satisfying \eqref{equ_sum_x_1_x_2_minus_x_0_in_N} and a finite set of multi-indices 
${\{ \mathbf{p}_{l} \}}_{l\in [L]} \subseteq \mathbb{Z}_{+}^{N}$. Then, at any \emph{$\mathbf{x}_{0} \!\in\! \mathcal{X}^{\text{o}}\!$}, the variance 
%% $v(\hat{g}(\cdot);\mathbf{x}_{0})$ 
%% at $\mathbf{x}_{0}$ 
of any estimator $\hat{g}(\cdot)$ with mean function $\gamma(\mathbf{x}) = \expect_{\mathbf{x}} \{ \hat{g}(\mathbf{y}) \}$ and finite variance at $\mathbf{x}_{0}$ is lower bounded as
\begin{equation} 
\label{equ_variance_bound_exp_families}
v(\hat{g}(\cdot);\mathbf{x}_{0}) \,\geq\, \mathbf{n}^{T}\rmv(\mathbf{x}_{0}) \ist\ist\mathbf{S}^{\dagger}(\mathbf{x}_{0}) \ist\ist\mathbf{n}(\mathbf{x}_{0}) \ist-\ist \gamma^{2}(\mathbf{x}_{0}) \,, 
\end{equation}  
where the vector $\mathbf{n}(\mathbf{x}_{0}) \rmv\in\rmv \mathbb{R}^{L}\rmv$ and the matrix $\mathbf{S}(\mathbf{x}_{0}) \rmv\in\rmv \mathbb{R}^{L \times L}\rmv$ 
are given elementwise 
\vspace*{1mm}
by
\begin{align}
\big( \mathbf{n}(\mathbf{x}_{0}) \big)_l &\ist\ist\triangleq 
% \frac{\partial^{\mathbf{p}_{k}}[ \gamma(\mathbf{x})  \lambda(\mathbf{x})]}{\partial \mathbf{x}^{\mathbf{p}_{k}}}\big|_{\mathbf{x} = \mathbf{x}_{0}} = 
\sum_{\mathbf{p} \leq \mathbf{p}_{l}} \!\binom{\mathbf{p}_{l}}{\mathbf{p}} \,\expect_{\mathbf{x}_{0}} \rmv\big\{ {\bm \phi}^{\mathbf{p}_{l}-\mathbf{p}}(\mathbf{y}) \big\} 
  \ist\ist\frac{\partial^{\mathbf{p}} \gamma(\mathbf{x})}{\partial \mathbf{x}^{\mathbf{p}}} \bigg|_{\mathbf{x} = \mathbf{x}_{0}} 
  \label{equ_variance_bound_exp_families_n} \\[.5mm]
\big( \mathbf{S}(\mathbf{x}_{0}) \big)_{l,l'} &\ist\ist\triangleq\, \expect_{\mathbf{x}_{0}} \rmv\big\{ {\bm \phi}^{\mathbf{p}_{l} + \mathbf{p}_{l'}}(\mathbf{y}) \big\} \,,
\label{equ_variance_bound_exp_families_S}
\end{align}
respectively. Here, $\sum_{\mathbf{p} \leq \mathbf{p}_{l}}$ denotes the sum over all multi-indices $\mathbf{p} \in \mathbb{Z}_{+}^{N}$ such that 
$p_{k} \leq (\mathbf{p}_{l})_{k}$ for $k\in [N]$, 
and $\binom{\mathbf{p}_{l}}{\mathbf{p}} \triangleq \prod_{k=1}^{N} \binom{ (\mathbf{p}_{l})_{k}}{p_{k}}$.
%%  denotes a multivariable extension of the binomial coefficient.
\vspace{1mm}
\end{theorem}

A proof of this result is provided in Appendix \ref{app_proof_thm_lower_bound_exp_family}. This proof shows that the bound \eqref{equ_variance_bound_exp_families} is obtained by projecting an appropriately 
transformed version of the mean function $\gamma(\cdot)$ onto the finite-dimensional subspace $\mathcal{U} = \linspan \rmv\big\{ r^{(\mathbf{p}_{l})}_{\mathbf{x}_{0}}(\cdot) \big\}_{l\in [L]}$ 
of an appropriately defined RKHS $\mathcal{H}(R)$, with the functions $r^{(\mathbf{p}_{l})}_{\mathbf{x}_{0}}(\cdot)$
%%  \in \mathcal{H}(R)$ 
given by \eqref{equ_def_part_der_func}. 
If we increase the set $\big\{ r^{(\mathbf{p}_{l})}_{\mathbf{x}_{0}}(\cdot) \big\}_{l\in [L]}$ by adding further functions $r^{(\mathbf{p}')}_{\mathbf{x}_{0}}(\cdot)$ 
with multi-indices $\mathbf{p}' \!\notin\rmv\rmv {\{\mathbf{p}_{l} \}}_{l \in [L]}$, the subspace 
%% on which we project 
tends to become higher-dimensional and in turn the lower bound \eqref{equ_variance_bound_exp_families} becomes higher, i.e., tighter. 

The requirement of a finite variance $v(\hat{g}(\cdot);\mathbf{x}_{0})$ in Theorem \ref{thm_lower_bound_exp_family} implies via \eqref{equ_nec_suff_cond_validity_bias_function} that 
$\gamma(\cdot) \in \mathcal{H}_{\mathcal{E}^{(A)}\!,\mathbf{x}_{0}}$. 
This, in turn, guarantees via \eqref{equ_der_reproduction_prop}---which can be invoked since due to Lemma \ref{lem_diff_RKHS_exp_family}
the RKHS $\mathcal{H}_{\mathcal{E}^{(A)}\!,\mathbf{x}_{0}}$ is differentiable up to any order at $\mathbf{x}_{0}$---the existence of the partial derivatives $\frac{\partial^{\mathbf{p}} \gamma(\mathbf{x})}{\partial \mathbf{x}^{\mathbf{p}}} \big|_{\mathbf{x} = \mathbf{x}_{0}}$.
Note also that the bound \eqref{equ_variance_bound_exp_families} depends on the mean function $\gamma(\cdot)$ only via its local behavior
as given by the the partial derivatives of $\gamma(\cdot)$ at $\mathbf{x}_{0}$
%%  \rmv\in\rmv \mathcal{X}^{\text{o}}\!$ 
up to a suitable order. %Note that if the variance $v(\hat{g}(\cdot);\mathbf{x}_{0})$ is finite, the existence of the partial derivatives is guaranteed 
%Note that for the specific choice $\{ \mathbf{p}_{l} = \mathbf{e}_{l} \}_{l \in [N]}$, the bound in \eqref{equ_variance_bound_exp_families} reduces to the ordinary \CRBfull \cite{LC,kay}. 

Evaluating the bound \eqref{equ_variance_bound_exp_families} requires computation of the moments 
$\expect_{\mathbf{x}_{0}} \rmv\big\{ {\bm \phi}^{\mathbf{p}}(\mathbf{y}) \rmv\big\}$. This can be done
%%  efficiently using 
by means of message passing algorithms \cite{GraphModExpFamVarInfWainJor}. 

\textcolor{red}{For the choice $L=N$ and $\mathbf{p}_{l} = \mathbf{e}_{l}$, the bound \eqref{equ_variance_bound_exp_families} is closely related to the CRB
%% \CRBfull 
obtained for the estimation problem $\mathcal{E}^{(A)}$. In fact, the 
CRB for $\mathcal{E}^{(A)}$
%% an estimation problem based on an exponential family 
is obtained as \cite[Thm. 2.6.2]{LC} 
%% \vspace{-2mm}
\begin{equation}
\label{equ_CRB_exp_family}
v(\hat{g}(\cdot);\mathbf{x}_{0}) \,\geq\, \mathbf{n}^{T}\rmv(\mathbf{x}_{0}) \ist\ist\mathbf{J}^{\dagger}(\mathbf{x}_{0}) \ist\ist\mathbf{n}(\mathbf{x}_{0}) \,,
\vspace{-.5mm}
\end{equation}
with $\big( n(\mathbf{x}_{0}) \big)_l = \frac{\partial \gamma(\mathbf{x})}{\partial x_{l}} \big|_{\mathbf{x} = \mathbf{x}_{0}}$ 
and the Fisher information matrix given by
%% $\mathbf{J}(\mathbf{x}_{0})$ given elementwise as
\[
% \mathbf{J}(\mathbf{x}_{0}) =  \expect_{\mathbf{x}_{0}} \rmv\big\{ \big(\phi_{l}(\mathbf{y}) - \expect_{\mathbf{x}_{0}} \rmv\big\{ \phi_{l}(\mathbf{y}) \big\}\big)
%\big(\phi_{l'}(\mathbf{y}) - \expect_{\mathbf{x}_{0}} \rmv\big\{ \phi_{l'}(\mathbf{y}) \big\}\big) \big \}, 
\mathbf{J}(\mathbf{x}_{0}) \,=\, \expect_{\mathbf{x}_{0}} \rmv\big\{ \big({\bm \phi}(\mathbf{y}) - \expect_{\mathbf{x}_{0}} \{ {\bm \phi}(\mathbf{y}) \} \big)
\transp{\big({\bm \phi}(\mathbf{y}) - \expect_{\mathbf{x}_{0}} \{ {\bm \phi}(\mathbf{y}) \} \big)} \big \} \,, 
\]
i.e., 
%% nothing but 
the covariance matrix of the sufficient statistic vector ${\bm \phi}(\mathbf{y})$. 
On the other hand, evaluating the bound \eqref{equ_variance_bound_exp_families} for 
%% the choice 
$L = N$ and $\mathbf{p}_{l} = \mathbf{e}_{l}$ and assuming without loss of generality that $\gamma(\mathbf{x}_{0}) =0$, we obtain
\begin{equation}
\label{equ_var_bound_exp_families_compare_CRB}
v(\hat{g}(\cdot);\mathbf{x}_{0}) \,\geq\, \mathbf{n}^{T}\rmv(\mathbf{x}_{0}) \ist\ist\mathbf{S}^{\dagger}(\mathbf{x}_{0}) \ist\ist\mathbf{n}(\mathbf{x}_{0}) \,,
\vspace{-2mm}
\end{equation}
with $\mathbf{n}(\mathbf{x}_{0})$ as before
%% $\big( n(\mathbf{x}_{0}) \big)_l = \frac{\partial \gamma(\mathbf{x})}{\partial x_{l}} \big|_{\mathbf{x} = \mathbf{x}_{0}}$ 
\vspace{-1mm}
and 
\[
 \mathbf{S}(\mathbf{x}_{0}) \,=\, \expect_{\mathbf{x}_{0}} \rmv\big\{{\bm \phi}(\mathbf{y}) \transp{\bm \phi}(\mathbf{y})\big \} \,.
\]
Thus, the only difference 
%% between the CRB in \eqref{equ_CRB_exp_family} and the bound in \eqref{equ_var_bound_exp_families_compare_CRB} 
is that the CRB in \eqref{equ_CRB_exp_family} involves the covariance matrix of the sufficient statistic ${\bm \phi}(\mathbf{y})$ 
whereas the bound in \eqref{equ_var_bound_exp_families_compare_CRB} involves the correlation matrix of ${\bm \phi}(\mathbf{y})$.
}

\vspace{-1mm}

\subsection{Reducing the Parameter Set}
\label{sec_param_set_reduction_exp_fam}
%%%%%%%%%%%%%%%%%%%%%%%%%%%%%%%%%%%%%%%%%%%%%%%%%%%%%%%%%%

Using the RKHS framework, we will now show that, under mild conditions, the minimum achievable variance $M(c(\cdot),\mathbf{x}_{0})$ for an exponential family type
estimation problem $\mathcal{E}^{(A)}=\big(\mathcal{X},f^{(A)}(\mathbf{y};\mathbf{x}),g(\cdot)\big)$ is invariant to reductions of the parameter set $\mathcal{X}$. 
Consider two estimation problems $\mathcal{E}=\scalarestproblem$ and $\mathcal{E}'\rmv\rmv =\big(\mathcal{X}'\rmv\rmv,f(\mathbf{y}; \mathbf{x}),g(\cdot)\big|_{\mathcal{X}'}\big)$---for 
now, not necessarily of the exponential family type---that differ only in their parameter sets $\mathcal{X}$ and $\mathcal{X}'\rmv\rmv$.
More specifically, $\mathcal{E}'\rmv$ is obtained 
from $\mathcal{E}$ by reducing the parameter set, i.e.,  $\mathcal{X}' \!\subseteq\rmv \mathcal{X}$. 
For these two estimation problems,
%%  $\mathcal{E}$ and $\mathcal{E}'\rmv$, 
we consider corresponding MVPs at a specific parameter vector $\mathbf{x}_0
%% \mathbf{x}' 
\!\in\! \mathcal{X}'$
%%  \subseteq \mathcal{X}$ 
and for a certain prescribed bias $c(\cdot)$. More precisely, $c(\cdot)$ is the prescribed bias 
for $\mathcal{E}$ on the set $\mathcal{X}$, while the prescribed bias for $\mathcal{E}'\rmv$ is the restriction of $c(\cdot)$ to $\mathcal{X}'\rmv$, 
%% which will be denoted as 
$c(\cdot)\big|_{\mathcal{X}'}$. We will denote the minimum achievable variances of the MVPs corresponding to $\mathcal{E}$ and $\mathcal{E}'$ 
by $M(c(\cdot), \mathbf{x}_{0})$ and $M'\big(c(\cdot)\big|_{\mathcal{X}'}, \mathbf{x}_{0}\big)$, respectively. 
From 
%% expression 
\eqref{equ_barankin_bound1}, it follows that $M'\big(c(\cdot)\big|_{\mathcal{X}'}, \mathbf{x}_{0}\big) \le M(c(\cdot), \mathbf{x}_{0})$,
%% the Barankin bound can only decrease, 
since taking the supremum over a reduced set can never result in an increase of the supremum.

The effect that a reduction of the parameter set $\mathcal{X}$ has on the minimum achievable variance can be analyzed conveniently within the RKHS framework.
\textcolor{red}{This is based on the following result \cite{aronszajn1950}:}
%\begin{theorem}[\!\!\cite{aronszajn1950}]
%\label{abc_thm_reducing_domain_RKHS}
Consider an RKHS $\mathcal{H}(R_{1})$ of functions $f(\cdot) \rmv : \mathcal{D}_{1}  \rmv\rmv\rightarrow\rmv \mathbb{R}$,  
%% associated 
with 
%% the 
kernel $R_{1}(\cdot\ist\ist,\cdot) \rmv: \mathcal{D}_{1} \rmv\rmv\times\rmv \mathcal{D}_{1} \rmv\rmv\rightarrow\rmv \mathbb{R}$.
%% , and note that the elements of $\mathcal{H}(R_{1})$ are functions with domain $\mathcal{D}_{1}$. 
Let $\mathcal{D}_{2} \rmv\subseteq\rmv \mathcal{D}_{1}$.
Then, the set of functions $\big\{ \tilde{f}(\cdot) \triangleq f(\cdot) \big|_{\mathcal{D}_{2}} \big|\, f(\cdot) \in \mathcal{H}(R_{1})\big\}$
%% _{f(\cdot) \in \mathcal{H}(R_{1})}$ 
that is obtained by restricting each function $f(\cdot) \rmv\in\rmv \mathcal{H}(R_{1})$ to 
the subdomain $\mathcal{D}_{2}$
%%  \rmv\subseteq\rmv \mathcal{D}_{1}$ 
coincides with the RKHS $\mathcal{H}(R_{2})$ 
%% associated with the 
whose kernel $R_{2}(\cdot\ist\ist,\cdot) \rmv: \mathcal{D}_{2} \rmv\times\rmv \mathcal{D}_{2} \rmv\rightarrow \mathbb{R}$ 
is the restriction of the kernel $R_{1}(\cdot\ist\ist,\cdot) \rmv: \mathcal{D}_{1} \rmv\rmv\times\rmv \mathcal{D}_{1} \rmv\rmv\rightarrow\rmv \mathbb{R}$ 
to the subdomain $\mathcal{D}_{2} \rmv\times\rmv \mathcal{D}_{2}$, i.e., 
\begin{equation}
\label{equ_rkhs_reducing_domain_equivalence_set_functions}
 \textcolor{red}{\mathcal{H}(R_{2}) \,=\, \big\{ \tilde{f}(\cdot) \triangleq f(\cdot) \big|_{\mathcal{D}_{2}} \big|\, f(\cdot) \in \mathcal{H}(R_{1})\big\} \,, \quad\mbox{with}\;\, 
R_{2} (\cdot\ist\ist, \cdot) \triangleq R_{1}(\cdot\ist\ist, \cdot)\big|_{\mathcal{D}_{2} \times \mathcal{D}_{2}}. }
\end{equation}
Furthermore, the norm of an element $\tilde{f}(\cdot) \rmv\in\rmv \mathcal{H}(R_{2})$ is equal to the minimum of the norms of all functions 
$f(\cdot) \rmv\in\rmv \mathcal{H}(R_{1})$ 
that coincide with $\tilde{f}(\cdot)$ on 
%% the restricted domain 
$\mathcal{D}_{2}$, i.e.,
%% we have the relation 
\begin{equation} 
\label{equ_thm_reducing_domain_RKHS}
{\| \tilde{f}(\cdot) \|}_{\mathcal{H}(R_{2})} \,= \min_{\substack{ \rule{0mm}{2.5mm}f(\cdot) \ist\in\ist \mathcal{H}(R_{1}) \\ f(\cdot) \big|_{\mathcal{D}_{2}} \!=\ist \tilde{f}(\cdot)}} {\| f(\cdot) \|}_{\mathcal{H}(R_{1})} \,.
% \mbox{s.t.} & \quad f(\mathbf{x}) = g(\mathbf{x}) \quad \forall \mathbf{x} \in \mathcal{D}_{2}. \nonumber 
\vspace{1.5mm}
\end{equation}
%% Thus, the norm of any element $f(\cdot) \in \mathcal{H}(R_{2})$ is equal to the minimal norm of a function $g(\cdot) \in \mathcal{H}(R_{1})$ 
%% that coincides with $f(\cdot)$ on the restricted domain $\mathcal{D}_{2}$. 
%\end{theorem} 

%% \begin{proof}
%% \cite{aronszajn1950}
%% \end{proof}

%% \noindent
Consider an arbitrary but fixed $f(\cdot) \in \mathcal{H}(R_{1})$, and let $\tilde{f}(\cdot) \triangleq f(\cdot) \big|_{\mathcal{D}_{2}}$.
Because $\tilde{f}(\cdot) \in \mathcal{H}(R_{2})$, we can calculate ${\| \tilde{f}(\cdot) \|}_{\mathcal{H}(R_{2})}$. From \eqref{equ_thm_reducing_domain_RKHS}, we obtain
for ${\| \tilde{f}(\cdot) \|}_{\mathcal{H}(R_{2})} = \big\| f(\cdot) \big|_{\mathcal{D}_{2}} \big\|_{\mathcal{H}(R_{2})}$ the inequality
\begin{equation} 
\label{equ_thm_reducing_domain_RKHS_1}
\big\| f(\cdot) \big|_{\mathcal{D}_{2}} \big\|_{\mathcal{H}(R_{2})} \,\le\, {\| f(\cdot) \|}_{\mathcal{H}(R_{1})} \,.
% \mbox{s.t.} & \quad f(\mathbf{x}) = g(\mathbf{x}) \quad \forall \mathbf{x} \in \mathcal{D}_{2}. \nonumber 
%% \vspace{1.5mm}
\end{equation}
This inequality holds for all $f(\cdot) \rmv\in\rmv \mathcal{H}(R_{1})$.
%%  and $$ such that $$.
%% Based on Theorem \ref{abc_thm_reducing_domain_RKHS}, we can 

Let us now return to the MVPs corresponding to $\mathcal{E}$ and $\mathcal{E}'\rmv$.
From \eqref{equ_thm_reducing_domain_RKHS_1} with $\mathcal{D}_{1} = \mathcal{X}$, $\mathcal{D}_{2} = \mathcal{X}'\rmv$, 
$\mathcal{H}(R_{1}) = \mathcal{H}_{\mathcal{E},\mathbf{x}_{0}}$, and $\mathcal{H}(R_{2}) = \mathcal{H}_{\mathcal{E}'\rmv\rmv,\mathbf{x}_{0}}$, 
we can conclude that, for any 
\vspace{1mm}
$\mathbf{x}_{0} \rmv\in\rmv \mathcal{X}'\rmv$,
\begin{equation}
\label{equ_min_var_problem_reduc_para_set_RKHS}
M'\big(c(\cdot)\big|_{\mathcal{X}'}, \mathbf{x}_{0}\big) \,\stackrel{\eqref{equ_min_achiev_var_sqared_norm}}{=}\, 
  \big\| \gamma(\cdot)\big|_{\mathcal{X}'} \big\|^{2}_{\mathcal{H}_{\mathcal{E}'\rmv\rmv,\mathbf{x}_{0}}} \rmv\!- \gamma^{2}(\mathbf{x}_{0})  
\,\stackrel{\eqref{equ_thm_reducing_domain_RKHS_1}}{\leq}\, {\| \gamma(\cdot) \|}^{2}_{\mathcal{H}_{\mathcal{E}\rmv,\mathbf{x}_{0}}} \rmv\!- \gamma^{2}(\mathbf{x}_{0}) 
\eq M(c(\cdot), \mathbf{x}_{0}) \,. 
\vspace{1mm}
\end{equation} 
Here, we also used the fact that $\gamma(\cdot)\big|_{\mathcal{X}'} =c(\cdot)\big|_{\mathcal{X}'} +g(\cdot)\big|_{\mathcal{X}'} $.
The inequality in \eqref{equ_min_var_problem_reduc_para_set_RKHS} means that a reduction of the parameter set $\mathcal{X}$ can never result in a deterioration of
the achievable performance, i.e., in a higher minimum achievable variance. 
Besides this rather intuitive fact, \textcolor{red}{the result \eqref{equ_rkhs_reducing_domain_equivalence_set_functions}} has the following consequence:
Consider an estimation problem $\mathcal{E}=\scalarestproblem$ whose statistical model ${\{ f(\mathbf{y}; \mathbf{x}) \}}_{\mathbf{x} \in \mathcal{X}}$ 
satisfies \eqref{equ_corr_likelihood_finite} at some $\mathbf{x}_{0} \!\in\! \mathcal{X}$ and moreover is contained in a ``larger'' model 
$\{ f(\mathbf{y}; \mathbf{x}) \}_{\mathbf{x} \in \tilde{\mathcal{X}}}$ with $\tilde{\mathcal{X}} \rmv\supseteq\rmv \mathcal{X}$. 
%% In what follows, we consider a fixed parameter vector $\mathbf{x}_{0} \!\in\! \mathcal{X}$ at which the condition \eqref{equ_corr_likelihood_finite} 
%% and the validity of bias functions is examined. 
If the larger model $\{ f(\mathbf{y}; \mathbf{x}) \}_{\mathbf{x} \in \tilde{\mathcal{X}}}$ also satisfies \eqref{equ_corr_likelihood_finite}, 
it follows from \textcolor{red}{\eqref{equ_rkhs_reducing_domain_equivalence_set_functions}} that a prescribed bias function $c(\cdot) \!: \mathcal{X} \!\to\rmv \mathbb{R}$ 
can only be valid for $\mathcal{E}$ at $\mathbf{x}_{0}$ if it is the restriction of a 
%% valid bias 
function $c'(\cdot) \!: \tilde{\mathcal{X}} \!\to\rmv \mathbb{R}$ that is a valid bias function for the estimation problem 
$\tilde{\mathcal{E}}\rmv\rmv=\big(\tilde{\mathcal{X}}\rmv\rmv,f(\mathbf{y}; \mathbf{x}),g(\cdot)\big)$ at $\mathbf{x}_{0}$. 
%% that differs from $\mathcal{E}$ only by its parameter set $\tilde{\mathcal{X}}\supseteq \mathcal{X}$. 
This holds true since every valid bias function for $\mathcal{E}$ at $\mathbf{x}_{0}$ is an element of the RKHS $\mathcal{H}_{\mathcal{E}\rmv,\mathbf{x}_{0}}$, 
which by \textcolor{red}{\eqref{equ_rkhs_reducing_domain_equivalence_set_functions}} consists precisely of the restrictions of the elements of the RKHS $\mathcal{H}_{\tilde{\mathcal{E}}\rmv,\mathbf{x}_{0}}$, 
which by \eqref{equ_nec_suff_cond_validity_bias_function} consists
%% is made up 
precisely of the mean functions that are valid for $\tilde{\mathcal{E}}$ at $\mathbf{x}_{0}$.% (see the remark made immediately after Theorem \ref{equ_main_thm_RKHS_MVE}). 
%% \newpage %%%%%%%%

For the remainder of this section, we restrict our discussion to estimation problems $\mathcal{E}^{(\mathcal{A})}=\big(\mathcal{X},f^{(A)}(\mathbf{y};\mathbf{x}),g(\cdot)\big)$
whose statistical model is an exponential family model. The next result characterizes the analytic properties of the mean functions $\gamma(\cdot)$ that belong to an 
RKHS $\mathcal{H}_{\mathcal{E}^{(\mathcal{A})}\rmv\rmv,\mathbf{x}_{0}}$. A proof is provided in Appendix \ref{app_proof_thm_any_valid_bias_func_analytic}.

\vspace{1mm}

\begin{lemma} 
\label{thm_any_valid_bias_func_analytic} 
Consider an estimation problem $\mathcal{E}^{(\mathcal{A})}=\big(\mathcal{X},f^{(A)}(\mathbf{y};\mathbf{x}),g(\cdot)\big)$ with an open parameter set 
$\mathcal{X} \rmv\rmv\subseteq\rmv \mathcal{N}$
%% $\mathcal{X}$ 
satisfying \eqref{equ_sum_x_1_x_2_minus_x_0_in_N} for some $\mathbf{x}_{0} \rmv\in\rmv \mathcal{X}$.
Let $\gamma(\cdot) \in \mathcal{H}_{\mathcal{E}^{(\mathcal{A})}\rmv\rmv,\mathbf{x}_{0}}$ be such that the partial derivatives 
$\frac{\partial^{\mathbf{p}} \gamma(\mathbf{x})}{\partial \mathbf{x}^{\mathbf{p}}} \big|_{\mathbf{x} = \mathbf{x}_{0}}$ 
vanish for every multi-index $\mathbf{p} \in \mathbb{Z}_{+}^{N}$. 
Then $\gamma(\mathbf{x})=0$ for all $\mathbf{x} \rmv\in\rmv \mathcal{X}$. 
\vspace{1mm}
\end{lemma} 

Note that since $\mathcal{H}_{\mathcal{E}^{(\mathcal{A})}\!,\mathbf{x}_{0}}$ is differentiable at $\mathbf{x}_{0}$ up to any order 
(see Lemma \ref{lem_diff_RKHS_exp_family}), it contains 
the function set $\big\{ r^{(\mathbf{p})}_{\mathbf{x}_0} (\mathbf{x}) \big\}_{\mathbf{p} \in \mathbb{Z}^{N}_{+}}$ defined in \textcolor{red}{\eqref{equ_def_part_der_func}}. 
Moreover, by \eqref{equ_der_reproduction_prop}, for any $f(\cdot) \in \mathcal{H}_{\mathcal{E}^{(\mathcal{A})}\!,\mathbf{x}_{0}}$ and any $\mathbf{p} \rmv\in\rmv \mathbb{Z}_{+}^{N}$,
there is $\big\langle r^{(\mathbf{p})}_{\mathbf{x}_0}(\cdot) , f(\cdot) \big\rangle_{\mathcal{H}_{\mathcal{E}^{(\mathcal{A})}\!,\mathbf{x}_{0}}} 
\!\!= \frac{\partial^{\mathbf{p}} f(\mathbf{x})}{\partial \mathbf{x}^{\mathbf{p}}} \big|_{\mathbf{x} = \mathbf{x}_{0}}$. Hence, under the assumptions of Lemma \ref{thm_any_valid_bias_func_analytic}, we have that if a function $f(\cdot) \in \mathcal{H}_{\mathcal{E}^{(\mathcal{A})}\!,\mathbf{x}_{0}}$
satisfies $\big\langle r^{(\mathbf{p})}_{\mathbf{x}_0}(\cdot) , f(\cdot) \big\rangle_{\mathcal{H}_{\mathcal{E}^{(\mathcal{A})}\!,\mathbf{x}_{0}}} \!\!=0$
for all $\mathbf{p} \rmv\in\rmv \mathbb{Z}_{+}^{N}$, then $f(\cdot) \equiv 0$. Thus, in this case, the set 
$\big\{ r^{(\mathbf{p})}_{\mathbf{x}_0} (\mathbf{x}) \big\}_{\mathbf{p} \in \mathbb{Z}^{N}_{+}}$ is complete for the RKHS 
$\mathcal{H}_{\mathcal{E}^{(\mathcal{A})}\!,\mathbf{x}_{0}}$.

Upon combining \textcolor{red}{\eqref{equ_rkhs_reducing_domain_equivalence_set_functions} and \eqref{equ_thm_reducing_domain_RKHS}} with Lemma \ref{thm_any_valid_bias_func_analytic}, we arrive at the second main result of this section: 

\vspace{1mm}

\begin{theorem}
\label{thm_par_set_reduction_est_problem_exp_family}
Consider an estimation problem $\mathcal{E}^{(\mathcal{A})} \rmv=\big(\mathcal{X},f^{(A)}(\mathbf{y};\mathbf{x}),g(\cdot)\big)$ with an open parameter set 
$\mathcal{X} \rmv\subseteq\rmv \mathcal{N}$ 
satisfying \eqref{equ_sum_x_1_x_2_minus_x_0_in_N} for some $\mathbf{x}_{0} \!\in\! \mathcal{X}$, and 
%% Consider furthermore 
a prescribed bias function $c(\cdot)$ that is valid for $\mathcal{E}^{(A)}$ at $\mathbf{x}_{0}$.
%%  \rmv\in\rmv \mathcal{X}$ 
%% and such that the corresponding mean function $\gamma(\cdot) =c(\cdot)+g(\cdot)$ 
%% is an element of $\mathcal{H}_{\mathcal{E}^{(\mathcal{A})}\rmv\rmv,\mathbf{x}_{0}}$.
%% Finally consider 
Furthermore consider a reduced parameter set $\mathcal{X}_{1} \rmv\subseteq\rmv \mathcal{X}$ such that \emph{$\mathbf{x}_{0} \rmv\in\rmv \mathcal{X}_{1}^{\text{o}}\rmv$}.  
Let $\mathcal{E}_{1}^{(A)} \!\triangleq\rmv \big(\mathcal{X}_{1},f^{(A)}(\mathbf{y}; \mathbf{x}); g(\cdot) \big)$
denote the estimation problem that is obtained from $\mathcal{E}^{(A)}$ by reducing the parameter set 
%% from $\mathcal{X}$ 
to $\mathcal{X}_{1}$,
%%  \subseteq \mathcal{X}$, 
and let $c_1(\cdot) \triangleq c(\cdot)\big|_{\mathcal{X}_{1}}$. Then, the minimum achievable variance for the restricted estimation problem $\mathcal{E}_{1}^{(A)}$ and 
the restricted bias function $c_1(\cdot)$, denoted by $M_{1}(c_1(\cdot),\mathbf{x}_{0})$, is equal to the minimum achievable variance for the original estimation problem 
$\mathcal{E}^{(A)}$ and the original bias function $c(\cdot)$, i.e.,
\[
M_{1}(c_1(\cdot),\mathbf{x}_{0}) \ist=\ist M(c(\cdot), \mathbf{x}_{0}) \,.
\vspace{1.5mm}
\]
\end{theorem}

A proof of this theorem is provided in Appendix \ref{app_proof_thm_par_set_reduction_est_problem_exp_family}.
Note that the requirement $\mathbf{x}_{0} \!\in\! \mathcal{X}_{1}^{\text{o}}$ of the theorem 
%% Theorem \ref{thm_par_set_reduction_est_problem_exp_family} 
implies that the reduced parameter set $\mathcal{X}_{1}$ must contain a neighborhood of $\mathbf{x}_{0}$, i.e., 
an open ball $\mathcal{B}(\mathbf{x}_{0},r)$ with some radius $r >0$. 
The main message of the theorem 
%% Theorem \ref{thm_par_set_reduction_est_problem_exp_family} 
is that, for an estimation problem based on an exponential family, parameter set reductions have no effect on the minimum achievable variance at $\mathbf{x}_{0}$ 
as long as the reduced parameter set contains a neighborhood of $\mathbf{x}_{0}$.

%% \newpage %%%%%%%

%%%%%%%%%%%%%%%%%%%%%%%%%%%%%%%%%%%%%%%%%%%%%%%%%%%%%%%%%%
\section{Conclusion}
%%%%%%%%%%%%%%%%%%%%%%%%%%%%%%%%%%%%%%%%%%%%%%%%%%%%%%%%%%

The mathematical framework of reproducing kernel Hilbert spaces (RKHS) provides powerful tools for the analysis of minimum variance estimation (MVE)
problems.
%% We considered the reproducing kernel Hilbert spaces (RKHS) approach to minimum variance estimation (MVE).
%%  within the classical (frequentist) estimation framework. 
Building upon the theoretical foundation developed in the seminal papers \cite{Parzen59} and \cite{Duttweiler73b}, we
derived novel results concerning the RKHS-based analysis of lower variance bounds for MVE, of sufficient statistics, and of MVE problems conforming to an exponential family of distributions.
More specifically, we presented an RKHS-based geometric interpretation of several well-known lower bounds on the estimator variance. 
We showed that each of these bounds is related to the orthogonal projection onto an associated subspace of the RKHS.
In particular, the subspace associated with the \CRBfull is based on
%% connected with 
the strong structural properties of a \emph{differentiable} RKHS.
For a wide class of estimation problems, we proved that the minimum achievable variance, which is the tightest possible lower bound on the estimator variance
(Barankin bound), 
is a lower semi-continuous function of the parameter vector. In some cases, 
%% This finding gives a criterion for judging wether 
this fact can be used to show that a given lower bound on the estimator variance is not maximally tight.
Furthermore, we proved that the RKHS associated with an estimation problem 
remains unchanged if the observation is replaced by a sufficient statistic. 

Finally, we specialized the RKHS description
%% approach 
to estimation problems whose observation conforms to an exponential family of distributions. We showed that the kernel of the 
%% associated 
RKHS has a particularly simple expression in terms of the moment-generating function of the 
%% underlying 
exponential family, and the RKHS itself is differentiable up to any order.
%% associated with estimation problems based on exponential families have the property of being differentiable. 
Using this differentiability, we derived novel closed-form lower bounds on the estimator variance.  
%% Moreover, based on the tight relationship between the RKHS and the moment-generating function, we 
We also showed that reducing the parameter set 
has no effect on the minimum achievable variance at a given reference parameter vector $\mathbf{x}_{0}$ if the reduced parameter set contains a neighborhood of 
$\mathbf{x}_{0}$. 
% These bounds are illustrated by means of numerical studies.

Promising
%% interesting 
directions for future work include the practical implementation of message passing algorithms for the efficient 
%% (possibly approximate) evaluation 
computation of the lower variance bounds 
%% on the estimator variance 
for exponential families derived in Section \ref{sec_lower_bound_exp_fam}.
%It would also be worthwhile to further develop the RKHS approach to estimation problems based on an exponential family with the goal of obtaining closed-form expressions 
%for the minimum achievable variance. 
Furthermore, in view of the close relations between exponential families and probabilistic graphical models \cite{GraphModExpFamVarInfWainJor}, it would be 
interesting to explore the relations between the graph-theoretic properties of the graph associated with
%% underlying 
an exponential family and the properties of the RKHS associated with that exponential family.
%It the context of graphical models \cite{}, 

%???THIS IS JUST A SUMMARY. ANY ``HIGHER-LEVEL CONCLUSIONS''? SOME CONCLUDING WORDS ON THE SIGNIFICANCE AND RELEVANCE OF OUR RESUTS? 
%SUGGESTIONS FOR FUTURE RESEARCH? (E.G., DEVELOP NEW, IMPROVED BOUNDS BY DEFINING NEW SUBSPACES ONTO WHICH TO PROJECT...)
%% \appendices
\appendices

\section{Proof of Theorem \ref{thm_lower_semi_cont_varying_kernel}}
 \label{app_proof_thm_lower_semi_cont_RKHS}
%%%%%%%%%%%%%%%%%%%%%%%%%%%%%%%%%%%%%%%%%%%%%%%%%%%%%%%%%%

\vspace{1mm}

%??? The finiteness of $M(c(\cdot),\mathbf{x})$ for every $\mathbf{x} \in \mathcal{X}$ follows...  ???THE THEOREM DIDN'T TALK ABOUT FINITENESS??? 

We first note that our assumption that the prescribed bias function $c(\cdot)$ is valid for $\mathcal{E}$ at every $\mathbf{x} \rmv\in\rmv \mathcal{C}$
has two consequences. First, $M(c(\cdot),\mathbf{x}) < \infty$ for every $\mathbf{x} \rmv\in\rmv \mathcal{C}$ (cf.\ our definition of the validity of a bias function in Section \ref{SecMVE});
second, due to \textcolor{red}{\eqref{equ_nec_suff_cond_validity_bias_function}}, the prescribed mean function $\gamma(\cdot) = c(\cdot) + g(\cdot)$ belongs to $\mathcal{H}_{\mathcal{E}\rmv,\mathbf{x}}$ 
for every $\mathbf{x} \rmv\in\rmv \mathcal{C}$.

Following \cite{Parzen59}, we define the \emph{linear span of a kernel function 
$R(\cdot\ist\ist,\cdot) \!: \mathcal{X} \rmv\rmv\times\rmv\rmv \mathcal{X} \rmv\rmv\rightarrow \mathbb{R}$}, denoted by $\mathcal{L}(R)$, as the
set of all functions $f(\cdot) \rmv: \mathcal{X} \rmv\rmv\rightarrow \mathbb{R}$ that are finite linear combinations of the form
\begin{equation}
\label{equ_def_linear_span_kernel}
f(\cdot) \ist\ist=\ist \sum_{l \in [L]} a_{l} \ist R(\cdot\ist,\mathbf{x}_{l}) \,, \quad\; \text{with} \;\,\mathbf{x}_{l}Ê\rmv\in\rmv \mathcal{X} \ist, \; a_{l} \rmv\in\rmv \mathbb{R} \,, \; L \rmv\in\rmv \mathbb{N} \,. 
\end{equation}
%% the associated linear span $\mathcal{L}(R)$ as the set of all finite linear combinations of functions of the form 
%% $f(\cdot) = R(\cdot, \mathbf{x})$ for some $\mathbf{x} \in \mathcal{X}$, i.e., 
%% \begin{equation}
%% \label{equ_def_linear_span_kernel}
%% \mathcal{L}(R) \triangleq  \bigg\{ f(\cdot): \mathcal{X} \rightarrow \mathbb{R} \big| f(\cdot) = \sum_{l = 1}^{L} a_{l} R(\cdot,\mathbf{x}_{l}) 
%% \mbox{, where } \mathbf{x}_{l}Ê\in \mathcal{X} \mbox{, } a_{l} \in \mathbb{R} \mbox{, and } L \in \mathbb{N} \bigg \}. 
%% \end{equation}
The linear span $\mathcal{L}(R)$ can be used to 
%% approximate 
express the norm of any function $h(\cdot) \rmv\in\rmv \mathcal{H}(R)$
%%  arbitrarily well, in the sense of the following result:
according to
%% \begin{lemma} 
%% Consider a kernel function $R(\cdot,\cdot): \mathcal{X} \times \mathcal{X} \rightarrow \mathbb{R}$ and the corresponding 
%% RKHS $\mathcal{H}(R)$ and linear span $\mathcal{L}(R)$, respectively. 
%% For every function $f(\cdot) \in \mathcal{H}(R)$, it holds that 
\begin{equation}
\label{equ_approx_norm_inner_prod_linear_span_RKHS}
 {\| h(\cdot) \|}_{\mathcal{H}(R)}^{2} \,=\!  \sup_{ \substack{\rule[-1.2mm]{0mm}{3.7mm} f(\cdot) \ist\in\ist \mathcal{L}(R) \\  {\| f(\cdot) \|}^{2}_{\mathcal{H}(R)} >\ist 0}} 
\!\rmv \frac{ {\langle  h(\cdot), f(\cdot) \rangle}_{\mathcal{H}(R)}^{2} }{ {\| f(\cdot) \|}^{2}_{\mathcal{H}(R)} } \,.
\end{equation}
%% \end{lemma} 
%% \begin{proof}
%% \cite[Theorem 3.1.2]{JungPHD} combined with \cite[Theorem 3.2.2]{JungPHD}. 
%% \end{proof}
This expression can be shown by combining \cite[Theorem 3.1.2]{JungPHD} and \cite[Theorem 3.2.2]{JungPHD}.
%% With this result, we can 
We can now develop the minimum achievable variance $M(c(\cdot), \mathbf{x})$ as follows:
\begin{align*}
M(c(\cdot),\mathbf{x}) & \,\stackrel{\eqref{equ_min_achiev_var_sqared_norm}}{=}\, {\| \gamma(\cdot) \|}_{\mathcal{H}_{\mathcal{E}\rmv,\mathbf{x}}}^{2} \!\rmv-\ist \gamma^{2}(\mathbf{x})  \nonumber \\[1mm]
   & \,\stackrel{\eqref{equ_approx_norm_inner_prod_linear_span_RKHS}}{=} \sup_{ \substack{\rule[-1.2mm]{0mm}{3.7mm} f(\cdot) \ist\in\ist \mathcal{L}(R_{\mathcal{E}\rmv,\mathbf{x}}) \\  {\| f(\cdot) \|}^{2}_{\mathcal{H}_{\mathcal{E}\rmv,\mathbf{x}}} >\ist 0}} 
\!\rmv \frac{ {\langle \gamma(\cdot), f(\cdot) \rangle}_{\mathcal{H}_{\mathcal{E}\rmv,\mathbf{x}}}^{2} }{ {\| f(\cdot) \|}^{2}_{\mathcal{H}_{\mathcal{E}\rmv,\mathbf{x}}} }
 \ist-\ist \gamma^{2}(\mathbf{x}) \,.
%% \nonumber \\[3mm]
\end{align*}
Using \eqref{equ_def_linear_span_kernel} and letting $\mathcal{D} \triangleq \{\mathbf{x}_{1},\ldots,\mathbf{x}_{L}\}$, 
$\mathbf{a} \triangleq ( a_{1} \cdots\ist a_{L})^T\!$, and
$\mathcal{A}_{\mathcal{D}} \ist\triangleq \big\{ \mathbf{a} \rmv\rmv\in\rmv \mathbb{R}^{L} \big|  \sum_{l,l'  \in [L]} a_{l} \ist a_{l'} 
R_{\mathcal{E}\rmv,\mathbf{x}}(\mathbf{x}_{l}, \mathbf{x}_{l'}) >0 \big\}$,
we obtain 
\vspace{-1mm}
further
\be
M(c(\cdot),\mathbf{x}) \,= \sup_{\mathcal{D} \subseteq \mathcal{X}, \ist\ist L\in \mathbb{N}, \ist\ist \mathbf{a} \in \mathcal{A}_{\mathcal{D}}} 
  \! h_{\mathcal{D},\mathbf{a}}(\mathbf{x}) \,.
\label{equ_proof_cont_kernel_finite_approx_sup_4}
\vspace{1mm}
\ee
%% (note that $L \rmv=\rmv | \mathcal{D}|$, and 
Here, our notation $\sup_{\mathcal{D} \subseteq \mathcal{X}, \ist\ist L\in \mathbb{N}, \ist\ist \mathbf{a} \in \mathcal{A}_{\mathcal{D}}}$ 
indicates that the supremum is taken not only with respect to the elements $\mathbf{x}_{l} $ of $\mathcal{D}$ but also with respect to
the size of $\mathcal{D}$,
%%  (number of elements), 
$L \rmv=\rmv | \mathcal{D}|$, and the function $h_{\mathcal{D},\mathbf{a}}(\cdot) \rmv: \mathcal{X} \rmv\rmv\rightarrow \mathbb{R}$ is given by
\begin{align*}
%% \label{equ_proof_cont_kernel_finite_approx_sup_1}
h_{\mathcal{D}, \mathbf{a}}(\mathbf{x}) 
& \,\triangleq\, \frac{ \big\langle \gamma(\cdot) \ist, \sum_{l \in [L]} a_{l} \ist R_{\mathcal{E}\rmv,\mathbf{x}}(\cdot\ist,\mathbf{x}_{l}) \big\rangle_{\mathcal{H}_{\mathcal{E}\rmv,\mathbf{x}}}^{2} }{
\big\| \sum_{l \in [L]} a_{l} \ist R_{\mathcal{E}\rmv,\mathbf{x}}(\cdot\ist,\mathbf{x}_{l}) \big\|_{\mathcal{H}_{\mathcal{E}\rmv,\mathbf{x}}}^{2} }  \ist-\ist \gamma^{2}(\mathbf{x})  \nonumber \\[1mm]
%%   \sum_{l,l'  \in [L]} a_{l} \ist a_{l'} R_{\mathcal{E}\rmv,\mathbf{x}}(\mathbf{x}_{l}, \mathbf{x}_{l'})} \ist-\ist \gamma^{2}(\mathbf{x}) \nonumber \\[1mm]
& \eq \frac{ \big( \sum_{l \in [L]} a_{l} \ist\ist \big\langle  \gamma(\cdot) \ist, R_{\mathcal{E}\rmv,\mathbf{x}}(\cdot\ist,\mathbf{x}_{l}) \big\rangle_{\mathcal{H}_{\mathcal{E}\rmv,\mathbf{x}}} \big)^{2} }{ 
  \sum_{l,l'  \in [L]} a_{l} \ist a_{l'} \ist\ist \big\langle R_{\mathcal{E}\rmv,\mathbf{x}}(\cdot\ist,\mathbf{x}_{l}) \ist\ist R_{\mathcal{E}\rmv,\mathbf{x}}(\cdot\ist,\mathbf{x}_{l'})
  \big\rangle_{\mathcal{H}_{\mathcal{E}\rmv,\mathbf{x}}} }  \ist-\ist \gamma^{2}(\mathbf{x})  \nonumber \\[1mm]
& \,\stackrel{\eqref{equ_reproducing_property}}{=}\, \frac{ \big( \sum_{l \in [L]} a_{l} \ist \gamma(\mathbf{x}_{l}) \big)^{2} }{ 
  \sum_{l,l'  \in [L]} a_{l} \ist a_{l'} R_{\mathcal{E}\rmv,\mathbf{x}}(\mathbf{x}_{l}, \mathbf{x}_{l'})} \ist-\ist \gamma^{2}(\mathbf{x}) \,. \nonumber \\[-9mm]
&\nonumber
\end{align*}
%The step $(a)$ follows from Theorem \ref{thm_main_facts_RKHS_MVE}, step $(b)$ is due to Theorem \ref{thm_approx_norm_dense_set} since the linear space $\mathcal{L}(R_{\mathcal{M}})$ 
%(cf.\ Definition \ref{def_linear_space_kernel}) is dense in the RKHS $\mathcal{H}(\mathcal{M})$ by Theorem \ref{thm_constr_RKHS_closure_linear_span}. 
%The step $(c)$ follows from the fact that any function $f(\cdot)$ belonging to the linear space $\mathcal{L}(R_{\mathcal{M}})$ can be written as a finite linear combination $f(\cdot) =  \sum_{l \in [L]} a_{l} R(\cdot,\mathbf{x}_{l})$ with 
%suitable coefficients $a_{l} \in \mathbb{R}$ and points $\mathcal{D} = \{ \mathbf{x}_{1}, \ldots, \mathbf{x}_{L} \} \subseteq \mathcal{X}$, which can be used in combination with the reproducing property \eqref{equ_reproduction_property} to 
%express the squared norm of $f(\cdot)$ as $ \| f(\cdot) \|^{2}_{\mathcal{H}(\mathcal{M})} =  \sum_{l,l'  \in [L]} a_{l} a_{l'} R_{\mathcal{M}(\mathbf{x})}(\mathbf{x}_{l}, \mathbf{x}_{l'})$. %Finally, step $(c)$ follows from the reproducing property \eqref{equ_reproduction_property}.   
For any finite set $\mathcal{D} = \{\mathbf{x}_{1},\ldots,\mathbf{x}_{L}\} \subseteq \mathcal{X}$ and any $\mathbf{a} \in \mathcal{A}_{\mathcal{D}}$, 
it follows from our assumptions of continuity of $R_{\mathcal{E}\rmv,\mathbf{x}}(\cdot\ist\ist,\cdot)$ with respect to $\mathbf{x}$ on $\mathcal{C}$ (see \eqref{equ_def_cont_varying_kernel}) 
and 
%% our assumption of 
continuity of $\gamma(\mathbf{x})$ on $\mathcal{C}$ that the function $h_{\mathcal{D},\mathbf{a}}(\mathbf{x})$ is continuous in a neighborhood around 
any point $\mathbf{x}_{0} \rmv\in\rmv \mathcal{C}$. 
%Using elementary linear algebra \cite{golub96} and the continuity of the kernel (cf. \eqref{equ_def_cont_varying_kernel}) we have that for any point $\mathbf{x}_{0} \in \mathcal{X}$ and every set $\mathcal{D}'$ such that $\rank(\mathbf{R}_{\mathbf{x}_{0}, \mathcal{D}'}) =|\mathcal{D}'|$ 
Thus, for any 
%% point 
$\mathbf{x}_{0} \rmv\in\rmv \mathcal{C}$, there exists a radius $\delta_{0} \rmv>\rmv 0$ such that 
%% the function 
$h_{\mathcal{D},\mathbf{a}}(\mathbf{x})$ is continuous on $\mathcal{B}(\mathbf{x}_{0},\delta_{0}) \subseteq \mathcal{C}$. 

 %the function $h_{\mathcal{D}'}(\mathbf{x})$ is continuos at $\mathbf{x}_{0}$ for every set $\mathcal{D}'$ such that $\rank(\mathbf{R}_{\mathbf{x}_{0}, \mathcal{D}'}) =|\mathcal{D}'|$.
We will now show that the function $M(c(\cdot),\mathbf{x})$ given by \eqref{equ_proof_cont_kernel_finite_approx_sup_4} is
%% must be 
lower semi-continuous at every $\mathbf{x}_{0} \rmv\in\rmv \mathcal{C}$, i.e., for any $\mathbf{x}_{0} \rmv\in\rmv \mathcal{C}$ and $\varepsilon >0$, 
we can find a radius $r>0$ such that 
\be
M(c(\cdot),\mathbf{x}) \,\geq\, M(c(\cdot),\mathbf{x}_{0}) \ist-\ist \varepsilon \,, \qquad \text{for all} \;\; \mathbf{x} \in \mathcal{B}(\mathbf{x}_{0},r) \,.
\label{equ_MM}
\ee
%Suppose $M(c(\cdot),\mathbf{x})$ is not lower semi-continuous at a specific parameter vector $\mathbf{x}_{0} \rmv\in\rmv \mathcal{C}$, i.e., 
%$\liminf_{\mathbf{x} \rightarrow \mathbf{x}_{0}} M(c(\cdot),\mathbf{x})  \leq M(c(\cdot),\mathbf{x}_{0}) - \varepsilon_{0}$ with some
%% a specific 
%$\varepsilon_{0} \rmv>\rmv 0$. 
%This implies that for any radius $r \rmv>\rmv 0$, there exists at least one parameter vector $\mathbf{x}' \!\in \mathcal{B}(\mathbf{x}_{0},r)$ 
%such that $M(c(\cdot),\mathbf{x}')  < M(c(\cdot),\mathbf{x}_{0})  - \varepsilon_{0}/2$, i.e, 
%\begin{equation}
%\label{equ_proof_lower_semi_equ_1}
%\text{for all} \; r \rmv>\rmv 0: \qquad
%M(c(\cdot),\mathbf{x}') \ist<\ist M(c(\cdot),\mathbf{x}_{0})  - \frac{\varepsilon_{0}}{2} \,, \quad \text{for some} \;\, \mathbf{x}' \! \in \mathcal{B}(\mathbf{x}_{0},r) \,.
%% ???\forall r>0: \exists \mathbf{x}' \in \mathcal{X} \cap \mathcal{B}(\mathbf{x}_{0}, r)Ê\Rightarrow M(c(\cdot),\mathbf{x}')  < M(c(\cdot),\mathbf{x}_{0})  - \frac{\varepsilon_{0}}{2} \,. 
%\end{equation} 
Due to \eqref{equ_proof_cont_kernel_finite_approx_sup_4}, there must be a finite subset $\mathcal{D}_{0} \rmv\subseteq\rmv \mathcal{X}$ 
and a vector $\mathbf{a}_{0} \rmv\in\rmv \mathcal{A}_{\mathcal{D}_{0}}$ such 
that\footnote{Indeed, %%%%%%%%%%%
if \eqref{equ_proof_lower_semi_equ_2} were not true, we would have 
$h_{\mathcal{D},\mathbf{a}}(\mathbf{x}_0) < M(c(\cdot),\mathbf{x}_{0})  - \varepsilon/2$ 
for every choice of $\mathcal{D}$ and $\mathbf{a}$. This, in turn, would imply that $\sup_{\mathcal{D} \subseteq \mathcal{X}, \ist\ist L\in \mathbb{N}, \ist\ist \mathbf{a} \in \mathcal{A}_{\mathcal{D}}} 
  \! h_{\mathcal{D},\mathbf{a}}(\mathbf{x}_{0}) \leq M(c(\cdot),\mathbf{x}_{0})  - \varepsilon/2 < M(c(\cdot),\mathbf{x}_{0})$, yielding the contradiction $M(c(\cdot),\mathbf{x}_{0}) Ê\stackrel{\eqref{equ_proof_cont_kernel_finite_approx_sup_4}}{=} 
  \sup_{\mathcal{D} \subseteq \mathcal{X}, \ist\ist L\in \mathbb{N}, \ist\ist \mathbf{a} \in \mathcal{A}_{\mathcal{D}}} h_{\mathcal{D},\mathbf{a}}(\mathbf{x}_{0}) < M(c(\cdot),\mathbf{x}_{0})$.
 } %%%%%%%%%
\begin{equation} 
\label{equ_proof_lower_semi_equ_2}
h_{\mathcal{D}_{0},\mathbf{a}_{0}}(\mathbf{x}_{0}) \,\geq\,  M(c(\cdot),\mathbf{x}_{0})  - \frac{\varepsilon}{2} \,,
\end{equation} 
for any given $\varepsilon> 0$.
Furthermore, since $h_{\mathcal{D}_{0},\mathbf{a}_{0}}(\mathbf{x})$ is continuous on $\mathcal{B}(\mathbf{x}_{0},\delta_{0})$ as shown above, 
there is a 
%% small 
radius $r_{0} \rmv>\rmv 0$ (with $r_{0} \rmv<\rmv \delta_{0}$) such that 
%% for any $\mathbf{x} \in \mathcal{B}(\mathbf{x}_{0},r_{0})$ we have 
\begin{equation} 
\label{equ_proof_lower_semi_equ_3}
h_{\mathcal{D}_{0},\mathbf{a}_{0}}(\mathbf{x}) \,\geq\, h_{\mathcal{D}_{0},\mathbf{a}_{0}}(\mathbf{x}_{0}) - \frac{\varepsilon}{2} \,, \qquad
  \text{for all} \;\; \mathbf{x} \in \mathcal{B}(\mathbf{x}_{0},r_{0}) \,.
\end{equation}  
By combining this inequality with \eqref{equ_proof_lower_semi_equ_2},
%%  and \eqref{equ_proof_lower_semi_equ_3}, 
it follows that there is a radius $r \rmv>\rmv 0$ (with $r \rmv<\rmv \delta_{0}$) such that for any $\mathbf{x} \in \mathcal{B}(\mathbf{x}_{0},r)$ we have
\begin{equation} 
\label{equ_proof_lower_semi_equ_4}
h_{\mathcal{D}_{0},\mathbf{a}_{0}}(\mathbf{x}) \,\stackrel{\eqref{equ_proof_lower_semi_equ_3}}{\geq}\, h_{\mathcal{D}_{0},\mathbf{a}_{0}}(\mathbf{x}_{0}) - \frac{\varepsilon}{2} \,\stackrel{\eqref{equ_proof_lower_semi_equ_2}}{\geq}\, 
%% M(c(\cdot),\mathbf{x}_{0})  - \frac{\varepsilon_{0}}{4} - \frac{\varepsilon_{0}}{4} \eq 
M(c(\cdot),\mathbf{x}_{0}) - \varepsilon \,,
\vspace{-3mm}
\end{equation}  
%% This lower bound on $h_{\mathcal{D}_{0},\mathbf{a}_{0}}(\mathbf{x})$ implies also a lower bound on the supremum in \eqref{equ_proof_cont_kernel_finite_approx_sup_4}, 
%% i.e., for those $\mathbf{x}$ for which the bound in \eqref{equ_proof_lower_semi_equ_4} is in force, we have 
%% simultaneously the lower bound 
and 
\vspace{.5mm}
further
\[ 
%% \label{equ_proof_lower_semi_equ_6}
M(c(\cdot),\mathbf{x}) \,\stackrel{\eqref{equ_proof_cont_kernel_finite_approx_sup_4}}{=}\! 
  \sup_{\mathcal{D} \subseteq \mathcal{X}, \ist\ist L\in \mathbb{N}, \ist\ist \mathbf{a} \in \mathcal{A}_{\mathcal{D}}} \! h_{\mathcal{D},\mathbf{a}}(\mathbf{x})
 \,\geq\, h_{\mathcal{D}_{0},\mathbf{a}_{0}}(\mathbf{x}) \stackrel{\eqref{equ_proof_lower_semi_equ_4}}{\geq} M(c(\cdot),\mathbf{x}_{0}) -  \varepsilon \,.
\vspace{1mm}
\]  
%% Putting together the pieces, we have that 
Thus, for any given $\varepsilon >0$, there is a radius $r \rmv>\rmv 0$ (with $r \rmv<\rmv \delta_{0}$) such that
%% \[ 
%% \label{equ_proof_lower_semi_equ_5}
%% M(c(\cdot),\mathbf{x}) \,\geq\, M(c(\cdot),\mathbf{x}_{0}) - \frac{\varepsilon_{0}}{2} \,, \qquad
%%   \text{for all} \;\; \mathbf{x} \in \mathcal{B}(\mathbf{x}_{0},r_{0}) \,.
%% \]  
$M(c(\cdot),\mathbf{x}) \geq M(c(\cdot),\mathbf{x}_{0}) - \varepsilon$ for all $\mathbf{x} \!\in\! \mathcal{B}(\mathbf{x}_{0},r)$, i.e.,
\eqref{equ_MM} has been proved.
%Because this contradicts \eqref{equ_proof_lower_semi_equ_1}, we have shown that $M(c(\cdot),\mathbf{x})$ is lower semi-continuous at all
%%  every point 
%$\mathbf{x} \rmv\in\rmv \mathcal{C}$. 
%Thus we have by \eqref{equ_proof_cont_kernel_finite_approx_sup_4} that the minimum achievable variance $L_{\mathcal{M}(\mathbf{x}_{0})}$ (viewed as a function of $\mathbf{x}_{0}$) is the point-wise minimum of continuos functions and 
%therefore lower semi-continous Zitat ``Counterexamples in Analyis''. 

%%%%%%%%%%%%%%%%%%%%%%%%%%%%%%%%%%%%%%%%%%%%%%%%%%%%%%%%%%
%\section*{Appendix C:\, Proof of Theorem \ref{thm_lower_bound_exp_family}}
\section{Proof of Theorem \ref{thm_lower_bound_exp_family}}
\label{app_proof_thm_lower_bound_exp_family}
%%%%%%%%%%%%%%%%%%%%%%%%%%%%%%%%%%%%%%%%%%%%%%%%%%%%%%%%%%

\vspace{1mm}

%% \emph{Proof}:\,
%% \begin{proof}
The bound \eqref{equ_variance_bound_exp_families} in Theorem \ref{thm_lower_bound_exp_family} is derived by using an
%%  specific 
isometry between the RKHS $\mathcal{H}_{\mathcal{E}^{(\mathcal{A})}\rmv\rmv,\mathbf{x}_{0}}$ 
and the RKHS $\mathcal{H}(R)$ that is defined by the kernel 
\begin{equation}
\label{equ_proof_exp_family_kernel_lambda_isometry}
R(\cdot\ist\ist,\cdot) \rmv:\ist \mathcal{X} \!\times\! \mathcal{X} \rmv\rmv\to \mathbb{R} \,, \qquad 
  R(\mathbf{x}_{1},\mathbf{x}_{2}) \ist=\ist \frac{\lambda(\mathbf{x}_{1} \rmv\rmv+\rmv \mathbf{x}_{2} \rmv-\rmv \mathbf{x}_{0})}{\lambda(\mathbf{x}_{0})} \,.
\end{equation}
It is easily verified that $R(\cdot\ist\ist,\cdot)$ and, thus, $\mathcal{H}(R)$ are differentiable up to any order.
Invoking \cite[Theorem 3.3.4]{JungPHD}, it can be verified that the two RKHSs $\mathcal{H}_{\mathcal{E}^{(\mathcal{A})}\rmv\rmv,\mathbf{x}_{0}}$ and $\mathcal{H}(R)$ 
are isometric and a specific congruence $\mathsf{J} \rmv:\ist \mathcal{H}_{\mathcal{E}^{(\mathcal{A})}\rmv\rmv,\mathbf{x}_{0}} \!\rightarrow \mathcal{H}(R)$ 
%% between these two RKHSs 
is given by 
\begin{equation} 
\label{equ_def_congruence_exp_families_kernel_RKHS}
\mathsf{J} [f(\cdot)] \ist=\ist \frac{\lambda(\mathbf{x})}{\lambda(\mathbf{x}_{0})} \, f(\mathbf{x}) \,.
%% \triangleq g(\cdot)\mbox{, where } g(\mathbf{x}) = \frac{\lambda(\mathbf{x})}{\lambda(\mathbf{x}_{0})} f(\mathbf{x}).
\end{equation}  
%In what follows,  of the moment generating function $\lambda(\mathbf{x})$ (cf.\ Theorem \ref{thm_part_der_moment_gen}), the kernel $R$ and in turn the RKHS $\mathcal{H}(R)$ is differentiable. 
Similarly to the bound  \eqref{equ_lower_bound_variance_projection}, we can then obtain a lower bound on 
%% the estimator variance 
$v(\hat{g}(\cdot);\mathbf{x}_{0})$ via an orthogonal projection onto a subspace of $\mathcal{H}(R)$. 
Indeed, with $c(\cdot) = \gamma(\cdot) - g(\cdot)$
%% \triangleq b(\hat{g}(\cdot); \mathbf{x})$ 
denoting the 
%% actual 
bias function of the estimator $\hat{g}(\cdot)$, we have 
\begin{align}
v(\hat{g}(\cdot); \mathbf{x}_{0}) & \,\ist\stackrel{\eqref{equ_lower_vound_variance_trivial_min_achiev_var}}{\geq}\, M(c(\cdot), \mathbf{x}_{0})\nonumber \\[.5mm]
&\ist\stackrel{\eqref{equ_min_achiev_var_sqared_norm}}{=}\ist {\| \gamma(\cdot) \|}^{2}_{\mathcal{H}_{\mathcal{E}^{(\mathcal{A})}\rmv\rmv,\mathbf{x}_{0}}}
  \!- \gamma^{2}(\mathbf{x}_{0}) \nonumber \\[.5mm]
& \ist\ist\stackrel{(a)}{=}\ist\ist \big\| \mathsf{J} [\gamma(\cdot)] \big\|^{2}_{\mathcal{H}(R)}
  - \gamma^{2}(\mathbf{x}_{0}) \nonumber \\[.5mm]
&\,\ist\geq\,\ist \big\| \big( \mathsf{J} [\gamma(\cdot)] \big)_{\mathcal{U}} \big\|^{2}_{\mathcal{H}(R)}
  - \gamma^{2}(\mathbf{x}_{0}) \,, 
\label{equ_lower_bound_variance_projection_isometry_exp_family}
\end{align}
for an arbitrary subspace 
%% $c(\mathbf{x}) = b(\hat{g}(\cdot), \mathbf{x})$ denotes the actual bias function of the estimator $\hat{g}(\cdot)$ and 
$\mathcal{U} \rmv\subseteq\rmv \mathcal{H}(R)$.
%%  is an arbitrary subspace of $$ and $\mathsf{P}_{\mathcal{U}}$ denotes the orthogonal projection operator on the subspace . 
Here, step $(a)$ is due to the fact that $\mathsf{J}$ is a congruence, and ${(\ist\cdot\ist)}_\mathcal{U}$ denotes orthogonal projection onto $\mathcal{U}$.
%% , and step $(b)$ follows from the projection theorem in Hilbert spaces \cite[Chapter 4]{RudinBook}. 
The bound \eqref{equ_variance_bound_exp_families} is obtained from \eqref{equ_lower_bound_variance_projection_isometry_exp_family} 
by choosing the subspace as
$\mathcal{U} \triangleq \linspan \rmv\rmv \big \{  r^{(\mathbf{p}_{l})}_{\mathbf{x}_{0}}(\cdot)  \big\}_{l\in [L]}$, with the functions 
$r^{(\mathbf{p}_{l})}_{\mathbf{x}_{0}}(\cdot) \in \mathcal{H}(R)$ as defined in \eqref{equ_def_part_der_func}, i.e.,
$r^{(\mathbf{p}_{l})}_{\mathbf{x}_{0}}(\mathbf{x}) = \frac{\partial^{\mathbf{p}_l} R(\mathbf{x}, \mathbf{x}_{2})}
{\partial \mathbf{x}_{2}^{\mathbf{p}_l}} \big|_{\mathbf{x}_{2} = \mathbf{x}_0}$. 

Let us denote the image of $\gamma(\cdot)$ under the isometry $\mathsf{J}$ by $\tilde{\gamma}(\cdot) \triangleq \mathsf{J} [\gamma(\cdot)]$. 
According to \eqref{equ_def_congruence_exp_families_kernel_RKHS},
%%  =  \mathsf{J} [\gamma(\cdot)]$, which according to \eqref{equ_def_congruence_exp_families_kernel_RKHS}, is given pointwise as 
\begin{equation}
\label{equ_proof_lower_bound_exp_fam_pointwise_1}
\tilde{\gamma}(\mathbf{x}) \ist=\ist \frac{\lambda(\mathbf{x})}{\lambda(\mathbf{x}_{0})} \, \gamma(\mathbf{x}) \,. 
%The function $v_{0}(\cdot)$ is orthogonal to every $ r^{(\mathbf{p}_{l})}_{\mathbf{x}_{0}}(\cdot)$, for $l \in \{1,\ldots,L\}$, since 
\end{equation}
Furthermore, the variance bound \eqref{equ_lower_bound_variance_projection_isometry_exp_family} reads
\[
v(\hat{g}(\cdot); \mathbf{x}_{0}) \,\geq\, {\| \tilde{\gamma}_{\mathcal{U}}(\cdot) \|}^{2}_{\mathcal{H}(R)}
  - \gamma^{2}(\mathbf{x}_{0}) \,. 
%% \label{equ_lower_bound_variance_projection_isometry_exp_family_1}
\]
Using \eqref{equ_projection_finite_dim_subspace}, we obtain further
\be
v(\hat{g}(\cdot); \mathbf{x}_{0}) \,\geq\, \mathbf{n}^{T}\rmv(\mathbf{x}_{0}) \ist\ist\mathbf{S}^{\dagger}(\mathbf{x}_{0}) \ist\ist\mathbf{n}(\mathbf{x}_{0}) 
  \ist-\ist \gamma^{2}(\mathbf{x}_{0}) \,,
\label{equ_lower_bound_variance_projection_isometry_exp_family_2}
\ee
where, according to \eqref{equ_n_S_0}, the entries of $\mathbf{n}(\mathbf{x}_{0})$ and $\mathbf{S}(\mathbf{x}_{0})$ are calculated as follows:
\begin{align}
\big( \mathbf{n}(\mathbf{x}_{0}) \big)_l &\stackrel{\eqref{equ_n_S_0}}{=}\ist \big\langle \tilde{\gamma}(\cdot), r^{(\mathbf{p}_{l})}_{\mathbf{x}_{0}}(\cdot) \big\rangle_{\mathcal{H}(R)}\nonumber \\[1.5mm]
&\stackrel{\eqref{equ_der_reproduction_prop}}{=}\ist \frac{\partial^{\mathbf{p}_{l}} \tilde{\gamma}(\mathbf{x})}{\partial \mathbf{x}^{\mathbf{p}_{l}}} 
  \bigg|_{\mathbf{x} = \mathbf{x}_{0}} \nonumber \\[2mm]
&\stackrel{\eqref{equ_proof_lower_bound_exp_fam_pointwise_1}}{=}\ist
\frac{1}{\lambda(\mathbf{x}_{0})} \frac{\partial^{\mathbf{p}_{l}} [\lambda(\mathbf{x})\gamma(\mathbf{x})]}{\partial \mathbf{x}^{\mathbf{p}_{l}}} 
  \bigg|_{\mathbf{x} = \mathbf{x}_{0}}  \nonumber \\[2mm]
& \ist\stackrel{(a)}{=}\ist \frac{1}{\lambda(\mathbf{x}_{0})} \sum_{\mathbf{p} \leq \mathbf{p}_{l}} \!\binom{\mathbf{p}_{l}}{\mathbf{p}} \ist\ist
  \frac{\partial^{\mathbf{p}_{l} - \mathbf{p}} \lambda(\mathbf{x})}{\partial \mathbf{x}^{\mathbf{p}_{l} - \mathbf{p}}}
  %% \bigg|_{\mathbf{x} = \mathbf{x}_{0}} \!
  \ist\frac{\partial^{\mathbf{p}} \gamma(\mathbf{x})}{\partial \mathbf{x}^{\mathbf{p}}} \bigg|_{\mathbf{x} = \mathbf{x}_{0}} \nonumber \\[2.5mm]
 & \stackrel{\eqref{equ_part_derivates_relation_moments_exp_family}}{=} \sum_{\mathbf{p} \leq \mathbf{p}_{l}} \!\binom{\mathbf{p}_{l}}{\mathbf{p}} \ist\ist
   \expect_{\mathbf{x}_{0}} \rmv\big\{ {\bm \phi}^{\mathbf{p}_{l}-\mathbf{p}}(\mathbf{y}) \big\} \ist\ist \frac{\partial^{\mathbf{p}} \gamma(\mathbf{x})}{\partial \mathbf{x}^{\mathbf{p}}} \bigg|_{\mathbf{x} = \mathbf{x}_{0}}
   %%  = m_{l},
\label{equ_expr_exp_coeffs_exp_family_proof}\\[-9mm]
\nonumber
\end{align} 
(here, $(a)$ is due to the generalized Leibniz rule for differentiation of a product of two functions \cite[p. 104]{RudinBookPrinciplesMatheAnalysis}), and 
%% we obtain 
\vspace{-3mm}
\begin{align}
\big( \mathbf{S}(\mathbf{x}_{0}) \big)_{l,l'} &\stackrel{\eqref{equ_n_S_0}}{=}\ist 
  \big\langle r^{(\mathbf{p}_{l})}_{\mathbf{x}_{0}}(\cdot), r^{(\mathbf{p}_{l'})}_{\mathbf{x}_{0}}(\cdot) \big \rangle_{\mathcal{H}(R)}\nonumber \\[1.5mm]
&\stackrel{\eqref{equ_der_reproduction_prop}}{=}\ist \frac{\partial^{\mathbf{p}_{l}} r^{(\mathbf{p}_{l'})}_{\mathbf{x}_{0}}(\mathbf{x})} {\partial \mathbf{x}^{\mathbf{p}_{l}}} 
  \bigg|_{\mathbf{x} = \mathbf{x}_{0}}\nonumber \\[1.5mm]
&\stackrel{\eqref{equ_def_part_der_func}}{=}\ist
\frac{\partial^{\mathbf{p}_{l}}} {\partial \mathbf{x}^{\mathbf{p}_{l}}} 
  \bigg\{ \frac{\partial^{\mathbf{p}_{l'}} \rmv R(\mathbf{x}, \mathbf{x}_{2})} {\partial \mathbf{x}_{2}^{\mathbf{p}_{l'}}} \bigg|_{\mathbf{x}_{2} = \mathbf{x}_0} \bigg\} 
    \bigg|_{\mathbf{x} = \mathbf{x}_{0}}  \nonumber \\[1.5mm] 
&\stackrel{\eqref{equ_proof_exp_family_kernel_lambda_isometry}}{=}\ist
\frac{1}{\lambda(\mathbf{x}_{0})} \ist\ist \frac{\partial^{\mathbf{p}_{l}+\mathbf{p}_{l'}} \lambda(\mathbf{x})}{\partial \mathbf{x}^{\mathbf{p}_{l}+\mathbf{p}_{l'}}} 
  \bigg|_{\mathbf{x} = \mathbf{x}_{0}}  \nonumber \\[1.5mm] 
& \stackrel{\eqref{equ_part_derivates_relation_moments_exp_family}}{=}\ist 
\expect_{\mathbf{x}_{0}} \rmv\big\{ {\bm \phi}^{\mathbf{p}_{l} + \mathbf{p}_{l'}}(\mathbf{y}) \rmv\big\} \,. 
%% = F_{k,l}.
\label{equ_exp_gramming_exp_family_proof}
\end{align}
Note that the application of \eqref{equ_der_reproduction_prop} was based on the differentiability of $\mathcal{H}(R)$.
Comparing \eqref{equ_lower_bound_variance_projection_isometry_exp_family_2}, \eqref{equ_expr_exp_coeffs_exp_family_proof}, and \eqref{equ_exp_gramming_exp_family_proof}
with \eqref{equ_variance_bound_exp_families}, \eqref{equ_variance_bound_exp_families_n}, and \eqref{equ_variance_bound_exp_families_S}, respectively,
we conclude that the theorem is proved.
%% \hfill $\Box$
%% \end{proof}

%%%%%%%%%%%%%%%%%%%%%%%%%%%%%%%%%%%%%%%%%%%%%%%%%%%%%%%%%%
%\section*{Appendix D:\, Proof of Lemma \ref{thm_any_valid_bias_func_analytic}}
\section{Proof of Lemma \ref{thm_any_valid_bias_func_analytic}}
 \label{app_proof_thm_any_valid_bias_func_analytic}
%%%%%%%%%%%%%%%%%%%%%%%%%%%%%%%%%%%%%%%%%%%%%%%%%%%%%%%%%%

\vspace{1mm}

%% \emph{Proof}:\,
%We first summarize three relevant properties of the parameter space $\mathcal{X}$: (i) $\mathcal{X} \rmv\subseteq\rmv \mathcal{N}$; (ii) $\mathcal{X}$ is an open set, i.e., 
%$\mathcal{X} \rmv\rmv=\rmv \mathcal{X}^{\text{o}}\rmv$; and (iii) $\mathcal{X}$ is a connected set ???because it satisfies \eqref{equ_sum_x_1_x_2_minus_x_0_in_N}.
%Since the bias function $c(\cdot)$ is assumed valid, there must exist an estimator $\hat{g}(\cdot)$ with finite variance at $\mathbf{x}_{0}$ and whose bias function and mean function
%equal $c(\cdot)$ and $\gamma(\cdot) = c(\cdot) + g(\cdot)$, respectively; ???furthermore, $\gamma(\mathbf{x}) = \expect_{\mathbf{x}} \{ \hat{g}(\cdot) \} < \infty$
%for all $\mathbf{x} \rmv\in\rmv \mathcal{X}$.
For $\mathcal{E}^{(\mathcal{A})}=\big(\mathcal{X},f^{(A)}(\mathbf{y};\mathbf{x}),g(\cdot)\big)$ and $\mathbf{x}_{0} \rmv\in\rmv \mathcal{X}$, 
consider a function $\gamma(\cdot) \!: \mathcal{X} \!\rightarrow\rmv \mathbb{R}$ belonging to the RKHS $\mathcal{H}_{\mathcal{E}^{(\mathcal{A})}\rmv\rmv,\mathbf{x}_{0}}$. 
By \textcolor{red}{\eqref{equ_nec_suff_cond_validity_bias_function}}, the function $c(\cdot) = \gamma(\cdot) - g(\cdot)$ is a valid bias function for 
$\mathcal{E}^{(\mathcal{A})}=\big(\mathcal{X},f^{(A)}(\mathbf{y};\mathbf{x}),g(\cdot)\big)$ at $\mathbf{x}_{0}$; furthermore, the LMV estimator at $\mathbf{x}_{0}$ exists 
and is given by $\hat{g}^{(\mathbf{x}_{0})}(\cdot) = \mathsf{J}[\gamma(\cdot)]$. 
Trivially, this estimator has the finite variance $v\big(\hat{g}^{(\mathbf{x}_{0})}(\cdot);\mathbf{x}_{0}\big) = M(c(\cdot),\mathbf{x}_{0})$ at $\mathbf{x}_{0}$ 
and its mean function equals $\gamma(\cdot)$, i.e., $\expect_{\mathbf{x}} \big \{\hat{g}^{(\mathbf{x}_{0})}(\mathbf{y}) \} = \gamma(\mathbf{x})$ for all 
$\mathbf{x} \!\in\! \mathcal{X}$. Hence, the mean power $\expect_{\mathbf{x}} \big\{ \big( \hat{g}^{(\mathbf{x}_{0})}(\mathbf{y}) \big)^{\rmv 2} \big\}$
%% of this estimator 
is finite at $\mathbf{x}_{0}$, since 
\begin{equation}
\label{equ_stoch_power_LMV_finite_at_x0}
\expect_{\mathbf{x}_{0}} \rmv\big\{ \big( \hat{g}^{(\mathbf{x}_{0})}(\mathbf{y}) \big)^{\rmv 2} \big\} 
  \eq v\big(\hat{g}^{(\mathbf{x}_{0})}(\mathbf{y});\mathbf{x}_{0}\big)  \ist+\ist  \big( \expect_{\mathbf{x}_{0}} \rmv\big\{ \hat{g}^{(\mathbf{x}_{0})}(\mathbf{y}) \big \} \big)^{\rmv 2} 
  \eq M(c(\cdot),\mathbf{x}_{0}) \ist+\ist \gamma^2(\mathbf{x}_{0}) \,<\, \infty \,. 
\end{equation}

Now, for any exponential family based estimation problem $\mathcal{E}^{(\mathcal{A})}=\big(\mathcal{X},f^{(A)}(\mathbf{y};\mathbf{x}),g(\cdot)\big)$, 
it follows from \cite[Theorem 2.7]{FundmentExpFamBrown} that the mean function $\expect_{\mathbf{x}} \{ \hat{g}(\cdot) \}$ of any estimator $\hat{g}(\cdot)$ 
%% $\gamma(\cdot)$ 
is
%% has to be 
analytic\footnote{Following %%%%%%%%%%
\cite[Definition 2.2.1]{KranzPrimerAnalytic}, we call a real-valued function $f(\cdot) \!: \mathcal{U} \rightarrow \mathbb{R}$ defined on some open domain $\mathcal{U}Ê\subseteq \mathbb{R}^{N}\rmv$ \emph{analytic} if 
for
%% at 
every point $\mathbf{x}_{c} \!\in\rmv \mathcal{U}$ there exists a power series $\sum_{\mathbf{p} \in \mathbb{Z}_{+}^{N}} a_{\mathbf{p}} (\mathbf{x} - \mathbf{x}_{c})^{\mathbf{p}}$ converging to $f(\mathbf{x})$ for every $\mathbf{x}$ in some neighborhood of $\mathbf{x}_{c}$. Note that the coefficients $a_{\mathbf{p}}$ may vary with 
$\mathbf{x}_{c}$.} %%%%%%%%%%%%%%  
on the interior $\mathcal{T}^{\text{o}}$ of the set 
$\mathcal{T} \triangleq \big\{ \mathbf{x} \!\in\! \mathcal{N} \ist\big|\ist \expect_{\mathbf{x}} \{ | \hat{g}(\mathbf{y}) | \} \rmv\rmv<\rmv\rmv \infty \big\}$.
Furthermore, $\mathcal{T}$ can be shown to be a convex set \cite[Corollary 2.6]{FundmentExpFamBrown}. 
In particular, the mean function $\gamma(\mathbf{x})$ of the LMV estimator $\hat{g}^{(\mathbf{x}_{0})}(\cdot)$ is analytic on the interior $\mathcal{T}_{0}^{\text{o}}$ 
of the convex set $\mathcal{T}_{0} \triangleq \big\{ \mathbf{x} \!\in\! \mathcal{N} \ist\big|\ist \expect_{\mathbf{x}} \{ |\hat{g}^{(\mathbf{x}_{0})}(\mathbf{y})| \} \rmv\rmv<\rmv\rmv \infty \big\}$. 
We will now verify that $\mathcal{X} \!\subseteq\rmv\rmv \mathcal{T}_{0}$. 
%First we have that $\mathbf{x} \in \mathcal{X}$ implies that $\mathbf{x} \in \mathcal{N}$, since $\mathcal{T} \subsetq \mathcal{N}$ by assumption. 
%By a reasoning similar to the proof of Theorem \ref{thm_regul_cond_CRB_imply_diff} in Appendix A, again
\textcolor{red}{Using} the Hilbert space $\mathcal{H}_{\mathbf{x}_{0}} \triangleq \big\{ t(\mathbf{y}) \ist\big|\ist \expect_{\mathbf{x}_{0}} \{ t^{2}(\mathbf{y}) \} \!<\! \infty \big\}$ 
and associated inner product ${\langle t_{1}(\mathbf{y}) , t_{2}(\mathbf{y}) \rangle}_{\text{RV}} =
%% _{\mathcal{H}_{\mathbf{x}_{0}}} \!\triangleq 
\expect_{\mathbf{x}_{0}} \{ t_{1}(\mathbf{y}) \ist t_{2}(\mathbf{y}) \}$,
we obtain for an arbitrary 
\vspace{1.5mm}
$\mathbf{x} \rmv\in\rmv \mathcal{X} \rmv\subseteq\rmv \mathcal{N}$
\begin{align*}
\expect_{\mathbf{x}} \{ |\hat{g}_{0}(\mathbf{y})| \} 
&\eq \expect_{\mathbf{x}_{0}} \rmv\bigg\{ \big| \hat{g}^{(\mathbf{x}_{0})}(\mathbf{y}) \big| \, \frac{f(\mathbf{y}; \mathbf{x})}{f(\mathbf{y}; \mathbf{x}_{0})} \bigg\}  \nonumber \\ 
& %% \stackrel{\eqref{equ_def_finite_power_hilbert_space},\eqref{equ_def_inner_prod_finite_power_hilbert_space}}{=} 
  \eq \big\langle \big| \hat{g}^{(\mathbf{x}_{0})}(\mathbf{y}) \big|, \rho(\mathbf{y}, \mathbf{x}) \big\rangle_{\text{RV}} 
 \nonumber \\[.5mm]
 & \stackrel{(a)}{\leq} \sqrt{ \big\langle \big| \hat{g}^{(\mathbf{x}_{0})}(\mathbf{y}) \big|, \big| \hat{g}^{(\mathbf{x}_{0})}(\mathbf{y}) \big| \big\rangle_{\text{RV}}
   \,\big\langle \rho(\mathbf{y}, \mathbf{x}), \rho(\mathbf{y}, \mathbf{x}) \big\rangle_{\text{RV}} }\nonumber \\[.5mm]Ê
 &  %% \stackrel{\eqref{equ_def_finite_power_hilbert_space},\eqref{equ_def_inner_prod_finite_power_hilbert_space}}{=} 
   \eq \sqrt{ \ist\ist \expect_{\mathbf{x}_{0}} \rmv \big\{  \big( \hat{g}^{(\mathbf{x}_{0})}(\mathbf{y}) \big)^{\rmv 2} \big\} 
   \,\expect_{\mathbf{x}_{0}} \rmv\bigg\{ \rmv\rmv\bigg( \frac{f(\mathbf{y}; \mathbf{x})}{f(\mathbf{y}; \mathbf{x}_{0})} \bigg)^{\rmv\!2}\bigg\} } \nonumber \\[.5mm]Ê
 & \!\!\!\rmv\stackrel{\eqref{equ_stoch_power_LMV_finite_at_x0},\eqref{equ_corr_likelihood_finite}}{\leq}\! \infty \,,
\end{align*}
where $(a)$ follows from the Cauchy-Schwarz inequality in the Hilbert space $\mathcal{H}_{\mathbf{x}_{0}}$.
Thus, we have verified that $\mathcal{X} \!\subseteq\rmv\rmv \mathcal{T}_{0}$. Moreover, we 
\vspace{-2mm}
have 
\begin{equation} 
\label{equ_open_parameter_set_included_interior_tau}
\mathcal{X} \rmv\rmv\subseteq\rmv \mathcal{T}_{0}^{\text{o}}.
\end{equation} 
This is implied\footnote{Indeed, %%%%%%%%% 
assume that the open set $\mathcal{X} \rmv\rmv\subseteq\rmv \mathcal{T}_{0}$ contains a vector $\mathbf{x}' \!\in\rmv \mathcal{X}$ that does not belong to the interior 
$\mathcal{T}_{0}^{\text{o}}$.
%%  of $\mathcal{T}_{0}$. 
It follows that no single neighborhood of $\mathbf{x}'$ can be 
%% fully 
contained in $\mathcal{T}_{0}$ and, thus, 
no single neighborhood of $\mathbf{x}'$ can be contained in $\mathcal{X}$, since $\mathcal{X} \rmv\rmv\subseteq\rmv \mathcal{T}_{0}$. 
However, because $\mathbf{x}'$ belongs to the open set $\mathcal{X} \!=\! \mathcal{X}^{\text{o}}\rmv$, there must be at least one neighborhood of $\mathbf{x}'$ that is 
%% fully 
contained in $\mathcal{X}$. Thus, we arrived at a contradiction, which implies that 
%% there cannot be a vector $\mathbf{x}' \!\in\rmv \mathcal{X}$ that does not belong to $\mathcal{T}_{0}^{\text{o}}$, meaning that
every vector $\mathbf{x}' \!\in\rmv \mathcal{X}$ must belong to $\mathcal{T}_{0}^{\text{o}}$, or, equivalently, 
that 
$\mathcal{X} \rmv\rmv\subseteq\rmv \mathcal{T}_{0}^{\text{o}}$.} %%%%%%%%%%%%%
by $\mathcal{X} \rmv\rmv\subseteq\rmv \mathcal{T}_{0}$ together with 
%% $\mathcal{X}^{\text{o}} \rmv\rmv=\rmv \mathcal{X}$ (recall that $\mathcal{X}$ was assumed an open set).
the fact that (by assumption) $\mathcal{X}$ is an open set.

Let us now consider the restrictions 
\begin{equation} 
\label{equ_restr_gamma_R_x_1}
\gamma_{\mathcal{R}_{\mathbf{x}_{1}}}\!(a) \,\triangleq\, \gamma \big(a \mathbf{x}_{1} + (1 \!-\rmv\rmv a) \ist\mathbf{x}_{0}\big) \,, \qquad a \!\in\! {(-\varepsilon,1+\varepsilon)} \,,
\end{equation}
of $\gamma(\cdot)$ on line segments of the form $\mathcal{R}_{\mathbf{x}_{1}} \triangleq \big\{ a \mathbf{x}_{1} + (1 \!-\rmv\rmv a) \ist\mathbf{x}_{0} \ist\big|\ist a \!\in\! {(-\varepsilon,1+\varepsilon)} \big\}$,
%% 0 \!\leq\! a \!\leq\! 1 \}$, 
where $\mathbf{x}_{1} \!\in\! \mathcal{T}^{\text{o}}_{0}$ and $\varepsilon>0$. Here, $\varepsilon$ is chosen sufficiently small such that the 
%% two 
vectors $\mathbf{x}_{a} \triangleq \mathbf{x}_{0} - \varepsilon (\mathbf{x}_{1}- \mathbf{x}_{0})$ and 
$\mathbf{x}_{b} \triangleq \mathbf{x}_{1} + \varepsilon (\mathbf{x}_{1}- \mathbf{x}_{0})$ belong to $\mathcal{T}_{0}^{\text{o}}$, i.e., $\mathbf{x}_{a}, \mathbf{x}_{b} \in \mathcal{T}_{0}^{\text{o}}$. Such an $\varepsilon$ can always be found, since---due to \eqref{equ_open_parameter_set_included_interior_tau}---we have 
$\mathbf{x}_{0} \rmv\in\rmv \mathcal{T}_{0}^{\text{o}}$.
%% , i.e., $\mathbf{x}_{0}$ like $\mathbf{x}_{1}$ (by assumption) is an interior point of $\mathcal{T}_{0}$. %Note that, since the set $\mathcal{T}_{0}$ is convex, we have $\mathcal{R}_{\mathbf{x}_{1}} \!\rmv\subseteq\rmv\rmv \mathcal{T}_{0}$ for any $\mathbf{x}_{1} \!\in\! \mathcal{T}_{0}$.
As can be verified easily, any vector in $\mathcal{R}_{\mathbf{x}_{1}}$ is a convex combination of the vectors $\mathbf{x}_{a}$ and $\mathbf{x}_{b}$, which both belong to the interior $\mathcal{T}^{\text{o}}_{0}$ of the convex set $\mathcal{T}_{0}$. Therefore we have $\mathcal{R}_{\mathbf{x}_{1}} \!\rmv\subseteq\rmv\rmv \mathcal{T}_{0}^{\text{o}}$ for any $\mathbf{x}_{1} \!\in\! \mathcal{T}^{\text{o}}_{0}$, as the interior  $\mathcal{T}^{\text{o}}_{0}$ of the convex set $\mathcal{T}_{0}$ is itself a convex set 
\cite[Theorem 6.2]{RockafellarBook},\footnote{Strictly %%%%%%%%%%
speaking, \cite[Theorem 6.2]{RockafellarBook} states that 
the \emph{relative interior} of a convex set is a convex set. However, since we assume that $\mathcal{X}$ is open with non-empty interior and therefore, by \eqref{equ_open_parameter_set_included_interior_tau}, also $\mathcal{T}_{0}$ has a nonempty interior, the relative interior of $\mathcal{T}_{0}$ 
coincides with the interior of 
$\mathcal{T}_{0}$.} %%%%%%%%%%
i.e., the interior $\mathcal{T}^{\text{o}}_{0}$ contains any convex combination of its elements.

The function $\gamma_{\mathcal{R}_{\mathbf{x}_{1}}}\!(\cdot) \!: (-\varepsilon,1+\varepsilon) \rightarrow \mathbb{R}$ in \eqref{equ_restr_gamma_R_x_1} 
is the composition of the mean function $\gamma(\cdot) \!: \mathcal{X} \rmv\rmv\rightarrow \mathbb{R}$, which is analytic on $\mathcal{T}_{0}^{\text{o}} \subseteq \mathcal{X}$, 
% $\gamma(\cdot): \mathcal{T}_{0}^{\text{o}} \rightarrow \mathbb{R}$ 
with the vector-valued function $\mathbf{b}(\cdot) \!: {(-\varepsilon,1+\varepsilon)} \rightarrow \mathcal{T}_{0}^{\text{o}}$ given by 
$\mathbf{b}(a) =a \mathbf{x}_{1} + (1 \rmv- a) \ist\mathbf{x}_{0}$. Since each 
component $b_{l}(\cdot)$ of the function $\mathbf{b}(\cdot)$, whose domain is the open interval ${(-\varepsilon,1+\varepsilon)}$, is an analytic function, 
the function $\gamma_{\mathcal{R}_{\mathbf{x}_{1}}}\!(\cdot)$ is itself analytic \cite[Proposition 2.2.8]{KranzPrimerAnalytic}.  
 
Since the partial derivatives of $\gamma(\cdot)$ at $\mathbf{x}_{0}$, $\frac{\partial^{\mathbf{p}} \gamma(\mathbf{x})}{\partial \mathbf{x}^{\mathbf{p}}} \big|_{\mathbf{x} = \mathbf{x}_{0}}$,
are assumed to vanish for every 
%% multi-index 
$\mathbf{p} \rmv\in\rmv \mathbb{Z}_{+}^{N}$, the (ordinary) derivatives of arbitrary order of the scalar function $\gamma_{\mathcal{R}_{\mathbf{x}_{1}}}\!(a)$ 
vanish at $a=0$ (cf.\ \cite[Theorem 9.15]{RudinBookPrinciplesMatheAnalysis}). 
According to \cite[Corollary 1.2.5]{KranzPrimerAnalytic}, since 
%% the derivatives of the analytic function $\gamma_{\mathcal{R}_{\mathbf{x}_{1}}}\!(a)$ vanish for any order at $a=0$, 
$\gamma_{\mathcal{R}_{\mathbf{x}_{1}}}\!(a)$ is an analytic function, 
this implies that $\gamma_{\mathcal{R}_{\mathbf{x}_{1}}}\!(a)$ vanishes everywhere on its open domain ${(-\varepsilon,1+\varepsilon)}$. 
This, in turn, implies that $\gamma(\cdot)$ vanishes on every line segment $\mathcal{R}_{\mathbf{x}_{1}}$ with some $\mathbf{x}_{1} \!\in\! \mathcal{T}^{\text{o}}_{0}$ 
and, thus, $\gamma(\cdot)$ vanishes everywhere on $\mathcal{T}_{0}^{\text{o}}$.
%% $a \!\in\! [0,1]$,
%i.e., 
%% separately to $\gamma_{\mathcal{R}_{\mathbf{x}_{1}}}\!(a)$ 
%for each , it follows 
%% the function 
%$\gamma(\cdot)$ vanishes on some
%% an arbitrary 
%ball $\mathcal{B}(\mathbf{x}_{0},r) \rmv\subseteq\rmv \mathcal{X} \subseteq \mathcal{T}_{0}$,
%%  with $r \rmv>\rmv 0$, 
%that if the derivatives of $\gamma_{\mathcal{R}_{\mathbf{x}_{1}}}\!(a)$ all vanish at $a=0$, 
%% the function 
%$\gamma(\cdot)$ vanishes on every line segment $\mathcal{R}_{\mathbf{x}_{1}}$ and in turn everywhere 
%% is identically zero 
%on , if the partial . 
By \eqref{equ_open_parameter_set_included_interior_tau}, we finally conclude that $\gamma(\cdot)$ vanishes everywhere on 
\vspace{-2mm}
$\mathcal{X}$.
%% \vspace{2mm}
%% \hfill $\Box$

%%%%%%%%%%%%%%%%%%%%%%%%%%%%%%%%%%%%%%%%%%%%%%%%%%%%%%%%%%
%\section*{Appendix E:\, Proof of Theorem \ref{thm_par_set_reduction_est_problem_exp_family}}
\section{Proof of Theorem \ref{thm_par_set_reduction_est_problem_exp_family}}
\label{app_proof_thm_par_set_reduction_est_problem_exp_family}
%%%%%%%%%%%%%%%%%%%%%%%%%%%%%%%%%%%%%%%%%%%%%%%%%%%%%%%%%%

\vspace{1mm}

%% \emph{Proof}:\,
Because $c(\cdot)$ was assumed valid at $\mathbf{x}_{0}$, the corresponding mean function $\gamma(\cdot) =c(\cdot)+g(\cdot)$ 
is an element of $\mathcal{H}_{\mathcal{E}^{(\mathcal{A})}\rmv\rmv,\mathbf{x}_{0}}$ (see \textcolor{red}{\eqref{equ_nec_suff_cond_validity_bias_function}}).
Let $\gamma_1(\cdot) \triangleq \gamma(\cdot)\big|_{\mathcal{X}_{1}}$, and note that $\gamma_1(\cdot)$ is the mean function 
corresponding to the restricted bias function $c_1(\cdot)$, i.e., $\gamma_1(\cdot) = c_1(\cdot) + g(\cdot)\big|_{\mathcal{X}_{1}}$. 
We have $\gamma_1(\cdot) \in \mathcal{H}_{\mathcal{E}_1^{(\mathcal{A})}\!,\mathbf{x}_{0}}$ due to \textcolor{red}{\eqref{equ_nec_suff_cond_validity_bias_function}}, 
because $\gamma_1(\mathbf{x})$ is the 
mean function (evaluated for $\mathbf{x} \in \mathcal{X}_{1}$) of an estimator $\hat{g}(\cdot)$ that has finite variance at $\mathbf{x}_{0}$ and
whose bias function on $\mathcal{X}$ 
%% evaluated at $\mathbf{x} \in \mathcal{X}$ 
equals $c(\mathbf{x})$. (The existence of such an estimator $\hat{g}(\cdot)$ is guaranteed 
since $c(\cdot)$ was assumed valid at $\mathbf{x}_{0}$.)
For the minimum achievable variance for the restricted estimation problem, we obtain
\be 
\label{equ_proof_param_set_reduc_exp_fam_proof_1}
M_{1}(c_1(\cdot),\mathbf{x}_{0}) \,\stackrel{\eqref{equ_min_achiev_var_sqared_norm}}{=}\, {\| \gamma_1(\cdot) \|}^{2}_{\mathcal{H}_{\mathcal{E}_1^{(\mathcal{A})}\!,\mathbf{x}_{0}}}
  \!\!-\ist \gamma_1^{2}(\mathbf{x}_{0}) 
\,\stackrel{\eqref{equ_thm_reducing_domain_RKHS}}{=}\rmv \min_{\substack{\rule{0mm}{3mm}\gamma'(\cdot) \ist\in\ist \mathcal{H}_{\mathcal{E}^{(\mathcal{A})}\rmv\rmv,\mathbf{x}_{0}} \\ \gamma'(\cdot)\big|_{\mathcal{X}_{1}} \!=\ist\ist \gamma_1(\cdot)}} \!{\| \gamma'(\cdot) \|}^{2}_{\mathcal{H}_{\mathcal{E}^{(\mathcal{A})}\rmv\rmv,\mathbf{x}_{0}}} \!\!-\ist \gamma_1^{2}(\mathbf{x}_{0}) \,.
\ee 
However, the only function $\gamma'(\cdot) \in \mathcal{H}_{\mathcal{E}^{(\mathcal{A})}\rmv\rmv,\mathbf{x}_{0}}$ that satisfies 
$\gamma'(\cdot)\big|_{\mathcal{X}_{1}} \!= \gamma_1(\cdot)$ is the mean function $\gamma(\cdot)$.
%%  (which corresponds to the prescribed bias $c(\cdot)$).
This is a consequence of Lemma \ref{thm_any_valid_bias_func_analytic} and can be verified as follows.
Consider a function $\gamma'(\cdot) \in \mathcal{H}_{\mathcal{E}^{(\mathcal{A})}\rmv\rmv,\mathbf{x}_{0}}$ that satisfies $\gamma'(\cdot)\big|_{\mathcal{X}_{1}} \!= \gamma_1(\cdot)$. 
By the definition of $\gamma_1(\cdot)$, we also have $\gamma(\cdot)\big|_{\mathcal{X}_{1}} \!= \gamma_1(\cdot)$. 
Therefore, the difference $\gamma''(\cdot)
%% : \mathcal{X} \rightarrowÊ\mathbb{R} 
\triangleq \gamma'(\cdot) - \gamma(\cdot) \in \mathcal{H}_{\mathcal{E}^{(\mathcal{A})}\rmv\rmv,\mathbf{x}_{0}}$ satisfies 
$\gamma''(\cdot) \big|_{\mathcal{X}_{1}} \!= \gamma'(\cdot) \big|_{\mathcal{X}_{1}} \!- \gamma(\cdot) \big|_{\mathcal{X}_{1}} \!= \gamma_1(\cdot) - \gamma_1(\cdot) = 0$,
i.e., $\gamma''(\mathbf{x}) = 0$ for all $\mathbf{x} \rmv\in\rmv \mathcal{X}_{1}$. 
Since $\mathbf{x}_{0} \rmv\in\rmv \mathcal{X}_{1}^{\text{o}}$, 
%% we can find a radius $r$ such that $\mathcal{B}(\mathbf{x}_{0}, r) \rmv\subseteq\rmv \mathcal{X}_{1}$,
%% and therefore $\gamma''(\cdot)$ vanishes on the open ball $\mathcal{B}(\mathbf{x}_{0}, r)$. 
this implies that $\frac{\partial^{\mathbf{p}} \gamma''(\mathbf{x})}{\partial \mathbf{x}^{\mathbf{p}}} \big|_{\mathbf{x} = \mathbf{x}_{0}} = 0$ 
for all $\mathbf{p} \in \mathbb{Z}_{+}^{N}$. It then follows from Lemma \ref{thm_any_valid_bias_func_analytic} 
that $\gamma''(\mathbf{x}) = 0$ 
%% vanishes everywhere on 
%% the open domain 
for all $\mathbf{x} \rmv\in\rmv \mathcal{X}$ and, thus, $\gamma'(\mathbf{x}) = \gamma(\mathbf{x})$ for all $\mathbf{x} \rmv\in\rmv \mathcal{X}$.
This shows that $\gamma(\cdot)$ is the unique
%% only 
function satisfying $\gamma(\cdot)\big|_{\mathcal{X}_{1}} \!= \gamma_1(\cdot)$. Therefore, we have 
\[
\min_{\substack{\rule{0mm}{3mm}\gamma'(\cdot) \ist\in\ist \mathcal{H}_{\mathcal{E}^{(\mathcal{A})}\rmv\rmv,\mathbf{x}_{0}} \\ \gamma'(\cdot)\big|_{\mathcal{X}_{1}} \!=\ist\ist \gamma_1(\cdot)}} 
\!{\| \gamma'(\cdot) \|}^{2}_{\mathcal{H}_{\mathcal{E}^{(\mathcal{A})}\rmv\rmv,\mathbf{x}_{0}}}
\!\rmv\eq {\| \gamma(\cdot) \|}^{2}_{\mathcal{H}_{\mathcal{E}^{(\mathcal{A})}\rmv\rmv,\mathbf{x}_{0}}} \,,
\] 
and thus \eqref{equ_proof_param_set_reduc_exp_fam_proof_1} becomes
\[ 
M_{1}(c_1(\cdot),\mathbf{x}_{0}) \eq {\| \gamma(\cdot) \|}^{2}_{\mathcal{H}_{\mathcal{E}^{(\mathcal{A})}\rmv\rmv,\mathbf{x}_{0}}} \!\!-\ist \gamma_1^{2}(\mathbf{x}_{0})
\eq
%% \stackrel{\gamma(\mathbf{x}_{0}) = \gamma_1(\mathbf{x}_{0})}{=}
{\| \gamma(\cdot) \|}^{2}_{\mathcal{H}_{\mathcal{E}^{(\mathcal{A})}\rmv\rmv,\mathbf{x}_{0}}} \!\!-\ist \gamma^{2}(\mathbf{x}_{0})  
\,\stackrel{\eqref{equ_min_achiev_var_sqared_norm}}{=}\, M(c(\cdot),\mathbf{x}_{0}) \,.
\]
Here, the second equality is due to the fact that $\gamma_1(\mathbf{x}_{0}) = \gamma(\mathbf{x}_{0})$ 
\pagebreak %%%%%%%%%
(because $\mathbf{x}_{0} \rmv\in\rmv \mathcal{X}_{1}^{\text{o}}$).
%% \vspace{2mm}
%% \hfill $\Box$

%%%%%%%%%%%%%%%%%%%%%%%%%%%%%%%%%%%%%%%%
\bibliographystyle{IEEEtran}
\bibliography{/Users/ajung/Arbeit/LitAJ_ITC.bib,/Users/ajung/Arbeit/tf-zentral}
%\bibliography{LitAJ_ITC.bib,tf-zentral}
%%%%%%%%%%%%%%%%%%%%%%%%%%%%%%%%%%%%%%%%

%%%%%%%%%%%%%%%%%%%%%%%%%%%%%%%%%%%%%%%%
\end{document}